%% file: my_thesis.tex
\documentclass[11pt]{ucthesis}
\usepackage{latexsym, amsmath, amssymb, amsfonts, amscd}
\usepackage{amsthm}
\usepackage{comment}
\usepackage{t1enc}
\usepackage[mathscr]{eucal}
\usepackage{indentfirst}
\usepackage{graphicx, pb-diagram}
\usepackage{fancyhdr}
\usepackage{fancybox}
\usepackage{enumerate}
\usepackage[all]{xy} 

\newtheorem{Theorem}{Theorem}[section]
\newtheorem{Proposition}[Theorem]{Proposition} 
\newtheorem{Lemma}[Theorem]{Lemma}

\newtheorem{Corollary}[Theorem]{Corollary}
\newtheorem{Conjecture}[Theorem]{Conjecture}

\newtheorem{Question}{Question}[chapter]
\newtheorem{Comment}[Theorem]{Comment}
\newtheorem{Definition}[Theorem]{Definition}

\newcommand{\low}{\text{low}}

\newcommand{\wt}{\operatorname{wt}}
\newcommand{\wi}{\mathbf{i}}
\newcommand{\bigdot}{\cdot}
\newcommand{\bz}{\Bbb{Z}}
\newcommand{\bc}{\Bbb{C}}

\newcommand{\bn}{\Bbb{N}}
\newcommand{\g}{\mathfrak{g}}

\newcommand{\asl}{\widehat{\mbox{sl}}}
\newcommand{\agl}{\widehat{\mbox{gl}}}

\newcommand{\slt}{\mbox{sl}_3}
\newcommand{\gli}{\mbox{gl}_\infty}
\newcommand{\hgl}{\widehat{\mbox{gl}}}
\newcommand{\Aa}{\mathcal{A}}
\newcommand{\LL}{\mathcal{L}}

\newcommand{\End}{\mbox{End}}

\newcommand{\Flip}{\mbox{Flip}}

\newcommand{\Id}{\text{Id}}

\newcommand{\C}{\mathcal{C}}

\newcommand{\modq}{{(mod \, q^{-1} \LL)}}
\newcommand{\rev}{\mathrm{rev}}
\newcommand{\barr}{\text{bar}}

\newcommand{\Br}{\sigma^\text{br}}

\theoremstyle{definition}

\renewcommand{\theenumi}{\roman{enumi}}

\newcommand{\mathfig}[2]{{\hspace{-3pt}\begin{array}{c}%
  \raisebox{-2.5pt}{\includegraphics[width=#1\textwidth]{#2}}
\end{array}\hspace{-3pt}}}

\begin{document}

% Declarations for Front Matter

\title{Some results on the crystal commutor and $\asl_n$ crystals}
\author{Peter William Tingley}
\degreesemester{Spring}
\degreeyear{2008} 
\degree{Doctor of Philosophy} 
\chair{Professor Nicolai Reshetikhin}
\othermembers{Professor Mark Haiman \\
  Professor Kam-Biu Luk }
\numberofmembers{3}
 \prevdegrees{B Math (University of Waterloo, Canada) 2001 \\ M Sc. (Carleton University, Canada) 2002} 
\field{Mathematics} 
\campus{Berkeley}

\maketitle
\approvalpage
\copyrightpage

\begin{abstract}
%remove the next line 
\def\dsp{\def\baselinestretch{1.28}\large\normalsize} \dsp
We present a number of results concerning crystal bases and their combinatorial models. There are two parts, which are largely independent. 

The first part consists of a series of results concerning the crystal commutor of Henriques and Kamnitzer. We first describe the relationship between the crystal commutor and Drinfeld's unitarized R-matrix. We then give a new definition for the crystal commutor, which makes sense for any symmetrizable Kac-Moody algebra. We show that this new definition agrees with A. Henriques and J. Kamnitzer's definition in the finite type case, but we cannot prove our commutor remains a coboundary structure in the non-finite type cases. Next, we give a new formula for the R-matrix and hence for the unitarized R-matrix, which we hope will be useful in eventually proving that the crystal commutor is a coboundary structure in all cases.

In the second part, we define three combinatorial models for $\asl_n$ crystals. These are parameterized by partitions, configurations of beads on an "abacus", and cylindric plane partitions, respectively. Our models are reducible, but we can identify an irreducible subcrystal corresponding to any dominant integral highest weight $\Lambda$. Cylindric plane partitions actually parameterize a basis for an irreducible representation of $\widehat{\text{gl}}_n$. This allows us to calculate the partition function for a system of random cylindric plane partitions first studied by A. Borodin. There is a symmetry in our model which allows us calculate the same partition function using the Weyl character formula of either a level $\ell$ representation for $\agl_n$ or a level $n$ representation for $\agl_\ell$. Thus we observe a form of rank level duality, originally due to I. Frenkel.  Finally, we use an explicit bijection to relate our work to the Kyoto path model.   
%put in for copy to be signed
%\abstractsignature
\end{abstract}

\begin{frontmatter}

\begin{dedication}
\null\vfil
{\large
\begin{center}
To Emily.
\end{center}}
\vfil\null
\end{dedication}

\tableofcontents
\listoffigures
%\listoftables
\begin{acknowledgements}
Many people have helped along road to this thesis, for which I am extremely grateful. 
I would especially like to thank my advisor Nicolai Reshetikhin for his patience and for many useful insights, both mathematical and otherwise. I also thank Joel Kamnitzer who was a collaborator for much of Chapter 3, and who has been a great pleasure to work with. Finally I would like to thank all the participants in the student representation theory, geometry and combinatorics seminar at Berkeley between 2003 and 2007; your talks, questions, comments and general hecklings have been invaluable to me. 
\end{acknowledgements}
%\end{comment}
\end{frontmatter}

\chapter{Introduction}

Quantum groups and in particular quantized universal enveloping algebras have 
provided a rich field of study. They have been used in many applications ranging from 
representation theory to topology to statistical mechanics. In the present work we study the representation theory of the quantized universal enveloping algebra $U_q(\g)$ associated to a symmetrizable Kac-Moody algebra $\g$. 

We focus on two aspects of this theory. The first is the fact that representations of $U_q(\g)$ have the structure of a braided monoidal category. In particular, there is a natural system of isomorphisms $\Br_{V, W}: V \otimes W \rightarrow W \otimes V$, for all pairs of representations $V$ and $W$, which satisfies the braid relations:
\begin{equation} \label{braiding}
\Br_{1,2} \circ \Br_{2,3} \circ \Br_{1,2} = \Br_{2,3} \circ \Br_{1,2} \circ \Br_{2,3} 
\end{equation}
as isomorphisms from $U \otimes V \otimes W$ to $W \otimes V \otimes U$, where $\Br_{i, i+1}$ means apply the braiding to the $i$ and $i+1$st factor of the tensor product. This is one of the main ingredients used in constructing the well known quantum group knot invariants. 

The other is the existence of crystal bases for representations of $U_q(\g)$. These are remarkable bases discovered by Kashiwara (see \cite{Kashiwara:1991}). By allowing $q$ to approach $0$ (or $\infty$) in an appropriate way, one can obtain a combinatorial object called a crystal associated to each integrable highest weight representation. This crystal consists of an underlying set along with operators $e_i$ and $f_i$ which correspond to the Chevalley generators $E_i$ and $F_i$ of $U_q(\g)$. A crystal can be depicted as an edge-colored directed graph, where the underlying set of the crystal is represented as the vertices of the graph, and the operators are encoded with the edges. One can use these crystals to answer various questions about the original representations. For instance, there is a simple tensor product rule for crystals which allows one to calculate tensor product multiplicities for the original representations. 

One might hope that the braiding on the category of $U_q(\g)$ representations would descend to a braiding on the category on $U_q(\g)$ crystals, but this is not the case. One should instead use the related notion of a coboundary structure, which is still a system of ismorphisms $\{ \sigma_{V,W}: V \otimes W \rightarrow W \otimes V \}$, but satisfying a different set of axioms than a braiding. In \cite{cactus}, Henriques and Kamnitzer construct a coboundary structure on the category of $U_q(\g)$ crystals for all finite type $\g$.

We are also interested in the relationship between crystals and enumerative combinatorics. It is often possible to realize crystals in a purely combinatorial way. For instance, in the case of the crystal associated with an irreducible representation of $\text{sl}_n$, the vertices are in bijection with semi-standard Young tableau of a certain shape. Furthermore, the weight of a vertex corresponds to a natural weight on the corresponding semi-standard Young tableau. Among other things, this shows that the generating function for semi-standard Young tableau of a given shape, counted with their weight, is equal to the character of the associated representation. In this way one recovers the well known correspondence between Young tableau, representations of  $\text{sl}_n$, and the Schur symmetric functions.

The first part of this work (Chapter \ref{commutor_chapter}) is concerned with understanding the coboundary structure on the category of crystals, and its relationship with the braiding on the category of representations. We show how the coboundary structure on crystals can be realized as a combinatorial limit of Drinfeld's coboundary structure on representations. We then attempt to generalize Henriques and Kamnitzer's coboundary structure to include the cases when $\g$ is a symmetrizable Kac-Moody algebra, but not necessarily of finite type.  We do not achieve this goal, but we give a conjectural solution based on Kashiwara's $*$-involution, and show that it agrees with  Henriques and Kamnitzer's  structure in finite type cases. In our attempts to prove our conjecture we obtain a new realization of the standard $R$-matrix for any symmetrizable Kac-Moody algebra.

In the second part (Chapter \ref{abacus_chapter}) we restrict our attention to the crystal associated to an integrable highest weight representations of $\asl_n$. We present several ways to realize these crystals combinatorially. In one of these models, the vertices of the crystal are given by cylindric plane partitions, as defined by Borodin \cite{Borodin:2006}. This crystal is reducible, but studying its structure reveals that cylindric plane partitions naturally parameterize a basis for an irreducible $\agl_n$ representation. This allows us to give a formula for the generating of cylindric plane partitions in terms of the Weyl character formula. It also allows us to observe a form of rank-level duality. 

Chapters \ref{commutor_chapter} and \ref{abacus_chapter} each describe the current state of an ongoing project, and each ends with a discussion of the directions we hope to explore in the future. We now give a somewhat more detailed description of their contents.

\section{Introduction to Chapter \ref{commutor_chapter}}

Let $\g$ be symmetrizable Kac-Moody algebra, and $U_q(\g)$ the associated quantized universal enveloping algebra (see \cite{CP}). Consider the following three categories:

$\bullet$ $Rep(\g$):= The category of integrable highest weight representations of $\g$.

$\bullet$ $Rep(U_q(\g)$):= The category of type 1 integrable highest weight representations 

\hspace{0.9in} of $U_q(\g)$.

$\bullet$ $Cryst(\g$):= The category of crystals of representations in $Rep(U_q(\g))$.

\noindent These are all semi-simple monoidal categories, and are in many ways very similar. In particular:

\begin{enumerate}
\item The irreducible objects in each category are indexed by the set $\Lambda_+$ of dominant integral weights of $\g$. 

\item For each $\lambda, \mu \in \Lambda_+$, it is clear that $V_\lambda \otimes V_\mu$ is isomorphic to $\bigoplus_\nu \oplus^{c_{\lambda \mu}^\nu} V_\nu$ for certain uniquely defined non-negative integers $c_{\lambda \mu}^\nu$. These $c_{\lambda \mu}^\nu$ are identical in the three categories. 

\item \label{cat_prop_3}  For any objects $X$ and $Y$ in any one of these categories, $X \otimes Y$ is isomorphic to $Y \otimes X$.
\end{enumerate}
 
It is often useful to give a category satisfying (\ref{cat_prop_3}) above additional structure by fixing a system of isomorphisms $\sigma_{X,Y}: X \otimes Y \rightarrow Y \otimes X$ for each pair of objects $(X, Y)$. There are different possible choices of this system of isomorphisms $\sigma = \{ \sigma_{X,Y} \}$ in the different categories we are considering:

$\bullet Rep(\g):$ The system of isomorphisms $\sigma^{sym}$ defined by $\sigma^{sym}_{V,W} (v \otimes w) = w \otimes v$ for all $v \in V$ and $w \in W$ gives $Rep(\g)$ the structure of a symmetric tensor category (see Definition \ref{sym_cat_def}). This implies that for any permutation $p$ of the integers $1$ to $N$, and any objects $V_1, \ldots V_N$ in $Rep(\g)$, there is a canonical isomorphism
\begin{equation}
I_p: V_1 \otimes \cdots \otimes V_N \rightarrow V_{p(1)} \otimes \cdots \otimes V_{p(N)}.
\end{equation}
In other words, the action of the system of isomorphisms $\sigma^{sym}$ on an $N$-fold tensor product factors through the symmetric group $S_N$.

$\bullet Rep(U_q(\g)):$ We will consider two different choices for $\sigma$. The first and most common is the standard braiding $\sigma^{br}$. This gives $Rep(U_q(\g))$ the structure of a braided monoidal category (see Definition \ref{braided_cat_def}), which essentially means that the action of the system of isomorphisms $\sigma^{br}$ on a $N$-fold tensor product factors through the braid group. There is also a second choice $\sigma^{D}$ which gives $Rep(U_q(\g))$ a structure which Drinfeld \cite{Drinfeld:1990} called a coboundary category (see Definition \ref{coboundary_cat_def}). This means that the action of $\sigma^D$ on an $N$-fold tensor product factors through a certain group $J_N$, sometimes called the cactus group (see \cite{cactus}). There are in fact many possible coboundary structures on $Rep(U_q(\g))$. There is however no way to give $Rep(U_q(\g))$ the structure of a symmetric monoidal category.

$\bullet Cryst(\g):$ When $\g$ is of finite type, Henriques and Kamnitzer defined a system of isomorphisms $\sigma^{HK}$ (which they call the crystal commutor) so that $(Cryst(\g), \sigma^{HK})$ is a coboundary category. There is no choice which makes $Cryst(\g)$ into a symmetric or even a braided monoidal category. In the present work we present a generalization of $\sigma^{HK}$ to non-finite type cases. We conjecture that the result is still a coboundary category, but this has not been proven.

In this chapter we will study these different systems of isomorphisms $\sigma$, and how they are related\footnote{Added May 1, 2008: An introduction to these objects can be found in the recent expository paper \cite{savage2}}. We will mainly focus of the categories $Rep(U_q(\g))$ and $Cryst(\g)$.

\subsection{Constructing $\sigma$ using systems of endomorphisms} \label{autint_intro}

We will now review a method for constructing natural systems of isomorphisms $\sigma_{V,W}: V \otimes W \rightarrow W \otimes V$ for the category of representations of $U_q(\g)$. This construction was used by Henriques and Kamnitzer in \cite{cactus}, and was further developed in \cite{Rcommutor}, where it is used to construct both the braiding $\sigma^{br}$ and Drinfeld's coboundary structure $\sigma^D$. The data needed to construct $\sigma$ in this way is:

\begin{enumerate}

\item An algebra involution $C_\xi$ of $U_q(\g)$ which is also a co-algebra anti-involution.

\item \label{diagg} A natural system of invertible (vector space) endomorphisms $\xi_V$ of each representation $V$ of $U_q(\g)$ such that the following diagram commutes for all $V$:
\begin{equation*}
\xymatrix{
V  \ar@(dl,dr) \ar@/ /[rrr]^{\xi_V} &&& V \ar@(dl,dr) \\
U_q(\g)  \ar@/ /[rrr]^{\C_{\xi}} &&&   U_q(\g). \\
}
\end{equation*}
\end{enumerate}
It follows immediately from the definition of coalgebra anti-automorphism that  \begin{equation} \sigma^\xi := \Flip \circ (\xi_V^{-1} \otimes \xi_W^{-1}) \circ \xi_{V \otimes W} \end{equation}  is an isomorphism of $U_q(\g)$ representations from $V \otimes W$ to $W \otimes V$. 

We will normally denote the system $\{ \xi_V \}$ simply by $\xi$, and will denote the action of $\xi$ on the tensor product of two representations by $\Delta(\xi)$. This is justified since, as explained in \cite{Rcommutor}, $\xi$ in fact belongs to a completion of $U_q(\g)$, and the action of $\xi$ on $V \otimes W$ is calculated using the coproduct.
We say a system of (vector space) endomorphisms of $U_q(\g)$ representations $\xi = \{ \xi_V \}$ is compatible with an algebra automorphism $C_\xi$ if the diagram in (\ref{diagg}) above commutes for all $V$. 

\begin{Comment} \label{ximake}
To describe the data $( \C_\xi, \xi )$, it is sufficient to describe $C_\xi$, and the action of $\xi_{V_\lambda}$ on any one vector $v$ in each irreducible representation of $V_\lambda$. This is usually more convenient then describing $\xi_{V_\lambda}$ explicitly. Of course, the choice of $C_\xi$ imposes a restriction on the possibilities for $\xi_{V_\lambda} (v)$, so when we give a description of $\xi$ in this way we are always claiming that the action on our chosen vector in each $V_\lambda$ is compatible with $C_\xi$.
\end{Comment}

Following ideas of  Kirillov-Reshetikhin \cite{KR:1990} and Levendorskii-Soibelman \cite{LS:1991}, we use this construction to construct the standard braiding.

\begin{Definition}
Let $C_X$ be the algebra automorphism of $U_q(\g)$ given by 
\begin{equation}
\begin{cases}
C_X(E_i)= -F_{\theta(i)} \\
C_X(F_i)=-E_{\theta(i)} \\
C_X(K_i)= K_{\theta(i)}^{-1},
\end{cases}
\end{equation}
where $\theta$ is the involution of the set $I$ such that $\alpha_{\theta(i)}= -w_0(\alpha_i)$.
\end{Definition}

\begin{Theorem}[{\cite[Theorem 3]{KR:1990}, \cite[Theorem 1]{LS:1991}}]  \label{makeX_intro}
Assume $\g$ is of finite type. Then there is a system of endomorphisms $X_{V}$ on all finite dimensional representations of $U_q(\g)$ which is compatible with $C_X$, and such that the standard braiding is given by 
\begin{equation}
\sigma^{br}_{V,W}= \Flip \circ  (X_V^{-1} \otimes X_W^{-1}) \circ X_{V \otimes W}. 
\end{equation}
\end{Theorem}

The action of $X_{V_\lambda}$ on an irreducible representation $V_\lambda$ can be described explicitly, and is related to the element $T_{w_0}$ of the braid group. In particular, $X_{V_\lambda}$ takes the highest weight space of $V_\lambda$ to the lowest weight space. For this reason $X_{V_\lambda}$ only makes sense when $\g$ is of finite type. 

We note that the endomorphisms $X_V$ can be realized using the action of an element $X$ in a certain completion of $U_q(\g)$. With this in mind, Theorem \ref{makeX_intro} becomes
\begin{equation}
\sigma^{br}  = \Flip \circ  (X^{-1} \otimes X^{-1}) \Delta(X).
\end{equation}
Equivalently 
$R  =  (X^{-1} \otimes X^{-1}) \Delta(X),$ where $R$ is the universal $R$-matrix for $U_q(\g)$ (see Section \ref{coboundary_section}). It is sometimes more convinient to work with the $R$-matrix rather then the braiding $\sigma^{br}$.

We also make use of a similar construction of Drinfeld's coboundary structure $\sigma^D$, which can be obtained via a modification of the construction for $\sigma^{br}$:

\begin{Theorem}[{\cite[Theorem 7.5]{Rcommutor}}] \label{makeY_intro}
Assume $\g$ is of finite type. Then there is a system of endomorphisms $Y_{V}$ on all finite dimensional representations of $U_q(\g)$ which is compatible with $C_X$, and such that Drinfeld's coboundary structure is given by  
\begin{equation}
\sigma^D_{V,W}= \Flip \circ (Y_V^{-1} \otimes Y_W^{-1}) \circ Y_{V \otimes W}.
\end{equation} 
\end{Theorem}

Since both $X_{V_\lambda}$ and $Y_{V_\lambda}$ are compatible with $C_X$, they differ only by an overall scalar (in fact, a power of $q$). This scalar depends on $\lambda$.

\subsection{Crystals and tensor products (informal)} \label{crystal_section_intro}

Before going on, let us briefly discuss Kashiwara's crystal bases. Due to our choice of conventions, we are forced to use crystal bases at $q=\infty$ rather then at $q=0$.
Let $V$ be a finite dimensional representation of $U_q(\g)$, and let $\Aa_\infty$ be the algebra of rational functions in $q$ that do not have poles at $\infty$. Kashiwara defined operators $\tilde{E}_i$ and $\tilde{F_i}$ on $V$, which are modifications of the actions of the Chevalley generators $E_i$ and $F_i$ on $V$ (see Section \ref{crystal_section}). He showed that one could find an $A_\infty$ lattice $\LL \subset V_\lambda$ (a crystal lattice), and a basis $B$ for $\LL/ q^{-1} \LL$ (the crystal basis) satisfying the following properties:

\begin{enumerate}
\item $\LL $ and $ B $ are compatible with the weight decomposition of $ V $.
\item $\LL $ is invariant under the Kashiwara operators and $ B \cup 0 $ is invariant under their residues $ e_i := \tilde{E}_i^\modq, f_i := \tilde{F}_i^\modq : \LL/q^{-1}\LL \rightarrow \LL/q^{-1} \LL$.
\item For any $ b, b' \in B $, we have $ e_i b = b' $ if and only if $ f_i b' = b $.
\end{enumerate}

\begin{figure}
\begin{center}
\input{tensor_figure}
\end{center}

\caption[A tensor product of two $\text{sl}_2$ crystals]{The tensor product $B \otimes C$ of two $\text{sl}_2$ crystals. $B$ and $C$ are represented by the strings on the top and the left of the diagram respectively. The underlying set of $B \otimes C$ consists of $b \otimes c$ for all $b \in B$ and $c \in C$, and is represented in the diagram by the remaining vertices. These vertices are connected with arrows which move along the $k^{th}$ row from the top, then down the $k^{th}$ row from the right, for all possible $k$. \label{tensor_figure}}

\end{figure}
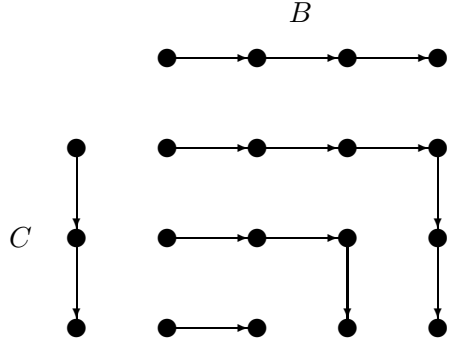

The set $B$, along with the operators $e_i$ and $f_i$ is referred to as the crystal of the representation $V$. It is often depicted as an edge-colored directed graph where the vertices consist of the elements of $B$, and for each $i$ in the index set for the Chevalley generators, the $i$ colored edges point from $b$ to $c$ if and only if $f_i(b)=c$. This crystal contains a lot of information about $V$. For instance, the connected components of $B$ naturally correspond to the irreducible components of $V$. Furthermore, crystal bases behave nicely under tensor products: if $(\LL_1, B_1)$ and $(\LL_2, B_2)$ are crystal bases for representations $V_1$ and $V_2$ respectively, then $(\LL_1 \otimes \LL_2, B_1 \times B_2)$ is a crystal basis for $V_1 \otimes V_2$. Furthermore, the action of $e_i$ and $f_i$ on the tensor product basis can be calculated as follows:

 $(\LL_1, B_1)$ and $(\LL_2, B_2)$ will remain crystal bases of $V_1$ and $V_2$ considered as representations of the copy of $\text{sl}_2$ generated by $E_i$ and $F_i$. We will calculate $e_i$ and $f_i$ on the tensor product crystal only considering its structure as a $\text{sl}_2$ crystal, so it is sufficient to give the definition for a tensor product of two irreducible $\text{sl}_2$ crystals. These consist of finite, linear directed graphs, where the operator $f_i$ moves one step in the direction of the arrows and $e_i$ moves one step against the arrows, and if you fall of the end of the graph, the element gets sent to zero. The tensor product of two such crystals is given in Figure \ref{tensor_figure}.

A highest weight element in a crystal is an element $b$ which is sent to zero by all the operators $e_i$. The crystal corresponding to the irreducible representation $V_\lambda$ has a unique highest weight element $b_\lambda$. Notice that the tensor product rule implies that all highest element in $B_\lambda \otimes B_\mu$ are of the form $b_\lambda \otimes c$ for various $c \in B_\mu$. This will be useful in our attempts to generalize Henriques and Kamnitzer's crystal commutor.

\subsection{Schutzenberger involution and the crystal commutor} \label{inv_commutor_intro}

The method described in Section \ref{autint_intro} cannot be used directly to construct the crystal commutor $\sigma^{HK}$ for the category $Cryst(\g)$, because $Cryst(\g)$ is not actually the category of representations of a Hopf algebra. However the construction used by Henriques and Kamnitzer \cite{cactus} proceeds along much the same lines.  They established the existence (and uniqueness) of a "Sch\"utzenberger involution" $ \xi_B : B \rightarrow B $, which satisfies the properties
\begin{equation}
\xi_B(e_i \cdot b) = f_{\theta(i)} \cdot \xi_B(b), \quad \xi_B(f_i \cdot b) = e_{\theta(i)} \cdot \xi_B(b),
\end{equation}
where $\theta$ is a certain diagram automorphism. This allows them to construct the crystal commutor as follows:

\begin{Theorem}[\cite{cactus}, Theorem 6] 
For any two $\g$-crystals $A$ and $B$, define the map $\sigma^{HK}_{A,B}: A \otimes B \rightarrow B \otimes A$ by
\begin{equation}
\begin{aligned} \label{commuter_definition_intro}
  \sigma^{HK}_{A,B} : A\otimes B &\rightarrow B \otimes A \\
 a \otimes b &\mapsto \mathrm{Flip} \circ (\xi_A \otimes \xi_B) \circ \xi_{A \otimes B} (a \otimes b).
 \end{aligned}
 \end{equation}
Then $\sigma^{HK}_{A,B}$ is always a crystal isomorphism, and the collection of all these isomorphisms gives $Cryst(\g)$ the structure of a coboundary category.
\end{Theorem}

\subsection{The crystal commutor as a limit of Drinfeld's coboundary structure}
Our first main theorem in Chapter \ref{commutor_chapter} essentially says that Drinfeld's coboundary structure on $Rep(U_q(\g))$ respects crystal lattices, and, in the crystal limit, one recovers Henriques and Kamnitzer's coboundary structure on crystal bases (up to some signs which are explained in Section \ref{ocb}). The standard braiding on the other hand does not respect crystal lattices. We feel this explains why the crystal commutor is a coboundary structure instead of a braiding, as one might have naively expected. We make this more precise is Section \ref{ocb}.

\begin{Theorem} \label{main1_intro} (see \cite{Rcommutor})
Let $(B_1, \LL_1)$ and $(B_2, \LL_2)$ be crystal bases for $V$ and $W$ respectively. Then 
\begin{enumerate}
\item $\sigma^D (\LL_1 \otimes \LL_2) =\LL_2 \otimes \LL_1$.

\item For any $b \in B_1$ and $c \in B_1$,  $\sigma^D (b_1 \otimes b_2) =\pm b_2 \otimes b_1$.
\end{enumerate}
\end{Theorem}

We prove Theorem \ref{main1_intro} by showing that the element $Y$ from Theorem \ref{makeY_intro} preserves the crystal lattice $\LL$, and acts on the crystal basis $B$ as Schutzenberger involution, up to some signs.

\subsection{A conjectural coboundary structure on crystals} \label{kcom_intro}

In \cite{commutorKashiwara} (joint work with Joel Kamnitzer) we show that the crystal commutor admits an alternative definition using Kashiwara's involution. This definition has the advantage that it makes sense even when $\g$ is not of finite type. Unfortunately, we do not know a proof that the crystal commutor will satisfy the conditions of a coboundary structure when $\g$ is not of finite type.

Our definition makes use of Kashiwara's $*$ involution, which is defined on $U_q(\g)$, and also on the "Verma crystal" $B_\infty.$ On $U_q(\g)$, this is defined on generators by
\begin{equation}
\begin{cases}
*E_i= E_i \\
*F_i=F_i \\
*K_i=K_i^{-1}.
\end{cases}
\end{equation}
One can easily check that $*$ is an algebra anti-automorphism, and that it is also a coalgebra automorphism. 

The Verma crystal $B_\infty$ is defined as follows: One can see that if $\lambda - \mu$ is a dominant integral weight, then there is an embedding of the crystal $B_\mu$ into the crystal $B_\lambda$ which is equivariant with respect to the operators $e_i$, and which takes the highest weight vector to the highest weight vector. These embedding turn irreducible $\g$ crystals into a directed system, and the limit of this system is $B_\infty$. 

Alternatively one can defined $B_\infty$ directly. There are analogues of the Kashiwara operators  $\tilde{E}_i$ and $\tilde{F}_i$ defined on $U_q^-(\g)$. Let $U_q^-(\g)$ be the sub-algebra of $U_q(\g)$ generated by all $F_i$, and let $\LL_\infty$ be the $\Aa_\infty$-lattice generated by all products $\tilde{F}_{i_n} \cdots \tilde{F}_{i_1}$. Then one can define $B_\infty$ to be the set consisting of all images of products $\tilde{F}_{i_n} \cdots \tilde{F}_{i_1}$ in $\LL_\infty/ q^{-1} \LL_\infty$. 

It was shown in \cite[Theorem 4]{Kashiwara:1991} that these two definitions of $B_\infty$ give the same crystal. However, some structure is visible in the second definition which is not in the first. Namely, $*$ involution naturally acts on $B_\infty$ according to the second definition. This action has the following property, which is essentially \cite[Prop 8.2]{Kashiwara:1991}:

\begin{Theorem}(Kashiwara) \label{*fact}
Let $c \in V_\mu \subset B_\infty$. Then $b_\lambda \otimes c$ is a highest weight element of $B_\lambda \otimes B_\mu$ if and only if $*c \in B_\lambda$, where we identify $B_\lambda$ and $B_\mu$ with their images in $B_\infty$. \qed
\end{Theorem}

Notice that in order to define a crystal isomorphism from $B_\lambda \otimes B_\mu$ to $B_\mu \otimes B_\lambda$, it is sufficient to specify the image of each highest weight element. By the tensor product rule for crystals (see Section \ref{crystal_section_intro} or \ref{crystal_section}), the highest weight vectors of $B_\lambda \otimes B_\mu$ are all of the form $b_\lambda \otimes c$ for various $c \in B_\mu$. In light of Theorem \ref{*fact}, we can define such an isomorphism as follows:

\begin{Definition}
For each pair of irreducible $\g$ crystals, define the crystal isomorphism $\sigma^{KT}_{V_\lambda, V_\mu}$ to take each highest weight element $b_\lambda \otimes c$ to $b_\mu \otimes *c$, where again we identify $B_\lambda$ and $B_\mu$ with their images in $B_\infty$.  Extend this by naturality to get a crystal isomorphism $\sigma^{KT}_{B,C}$ for any $\g$-crystals $B$ and $C$.
\end{Definition}

The following is the second main theorem in Chapter \ref{commutor_chapter}:

\begin{Theorem} \label{Main2_intro}
If $\g$ is of finite type, then $\sigma^{HK}=\sigma^{KT}$.
\end{Theorem}

Since $\sigma^{HK}$ is a coboundary structure we see that, if $\g$ is of finite type, $\sigma^{KT}$ is as well. We have not used anything specific to finite type algebras in defining $\sigma^{KT}$, so it seems natural to hope this remains a coboundary structure in other cases. Thus we conjecture:

\begin{Conjecture} \label{commutor_conjecture_intro}
For any symmetrizable Kac-Moody algebra $\g$, $\sigma^{KT}$, as defined above, is a coboundary structure on $Cryst(\g)$.\footnote{This has recently been proven by Savage \cite[Theorem 6.4]{Savage:2008}}
\end{Conjecture}

\subsection{An alternate construction of the braiding}

One source of difficulty in constructing a coboundary structure on $Cryst(\g)$ when $\g$ is not of finite type is that Theorems \ref{makeX_intro} and \ref{makeY_intro} do not make sense. The endomorphisms $X$ and $Y$ used in these theorems take the highest weight space of an irreducible representation to the lowest weight space, and hence are not well defined unless $\g$ is of finite type. We will now give an  analogue of Theorem \ref{makeX_intro} that is valid in all cases. We feel this is significant progress toward proving Conjecture \ref{commutor_conjecture_intro}, although it also raises some new difficulties.

As in Section \ref{autint_intro}, the data we need to construct a system of isomorphisms from $V \otimes W$ to $W \otimes V$ is an algebra automorphism $C_\Theta$ of $U_q(\g)$ which is also a coalgebra anti-automorphism, along with a corresponding system of automorphisms $\Theta_V$ of each representation. We now require that these make sense for any symmetrizable Kac-Moody algebra. There is a $C_\Theta$ with the right properties, although it is not linear over the base field $\bc(q)$, but is instead $\barr$-linear (i.e. is compatible with the automorphism of $\bc(q)$ which inverts $q$):
\begin{equation}
\begin{cases}
C_\Theta (E_i) =  K_i^{-1} E_{i} \\
C_\Theta (F_i) = F_{i}  K_i \\
C_\Theta (K_i) = K_{i}^{-1} \\
C_\Theta (q) = q^{-1}.
\end{cases}
\end{equation}
The fact that this involution is bar-linear instead of linear makes it somewhat more difficult to define the necessary endomorphisms $\Theta_{V_\lambda}$ of the irreducible representations $V_\lambda$. In fact, we need to work not in the category of representations itself, but instead in the category of representations with a chosen basis for the highest weight space of each isotypic component (for $V_\lambda$ that just means choosing a highest weight vector $v_\lambda$). Given a pair $(V_\lambda, v_\lambda)$, one fixes $\Theta_{(V_\lambda, v_\lambda)}$ by insisting that it send $v_\lambda$ to $q^{(-\lambda, \lambda)/2+(\lambda, \rho)} v_\lambda$, and that it be compatible with $\Theta$. 

One also needs to define a tensor product on this new category. That is, one must decide on a way to choose a basis for the highest weight space of each isotypic component of $V_\lambda \otimes V_\mu$, depending on the chosen highest weight vectors $v_\lambda \in V_\lambda$ and $v_\mu \in V_\mu$. This can be done in several ways. Once such a tensor product is chosen, one can define
\begin{equation}
\sigma^\Theta_{V,W}:= \Flip \circ \Theta_V \otimes \Theta_W \circ \Theta_{V \otimes W}^{-1}.
\end{equation}

It follows as in Section \ref{autint_intro} that $\sigma^\Theta_{V,W}$ will be an isomorphism of representations. With the tensor product we use one can show that it will not depend on the chosen bases for the highest weight spaces of $V$ and $W$, so it is well defined in the category of representations itself, without chosen highest weight vectors. The following is the third main theorem in Chapter \ref{commutor_chapter}.

\begin{Theorem} \label{RTheta_general_th_intro}
With the appropriate choice of normalization and tensor product, $\sigma^\Theta$ can be made to agree with either the standard braiding or Drinfeld's coboundary structure.
\end{Theorem}
We hope that Theorem \ref{RTheta_general_th_intro} will lead to a new proof of Conjecture \ref{commutor_conjecture_intro}.

\section{Introduction to Chapter \ref{abacus_chapter}}

The results in this chapter were motivated by the Hayashi realization for crystals of level one $\asl_n$ representations, originally developed by Misra and Miwa \cite{MM:1990} using work of Hayashi \cite{Hayashi:1990} (see also \cite{ariki:2000}, Chapter 10). In that realization, the underlying set of the crystal consists of partitions, and the operators $f_i$ act by adding a box to the associated Young diagram. We wondered if there was a similar realization for representations of arbitrary level $\ell$, where the operators $f_i$ would act by adding an $\ell$-ribbon (see Section \ref{lotsastuff}). It turns out that there is. To prove that our construction works, we need a second model, which is based on the abacus used by James and Kerber (\cite{JK:1981}, Chapter 2.7). The abacus model is reminiscent of  a "Dirac sea,"  and we think it is interesting in its own right. 

The crystals one obtains from the abacus model are not irreducible. However, one can pick out a "highest" irreducible sub-crystal, so we do have a model for the crystal of any irreducible integrable highest weight representation of $\asl_n$.  In finding this irreducible subcrystal, one is led to consider another subcrystal which is reducible. This turns out to allow a different combinatorial model, based on the cylindric plane partitions studied by Borodin \cite{Borodin:2006}, and has some interesting consequences.

\subsection{The partition and abacus models}

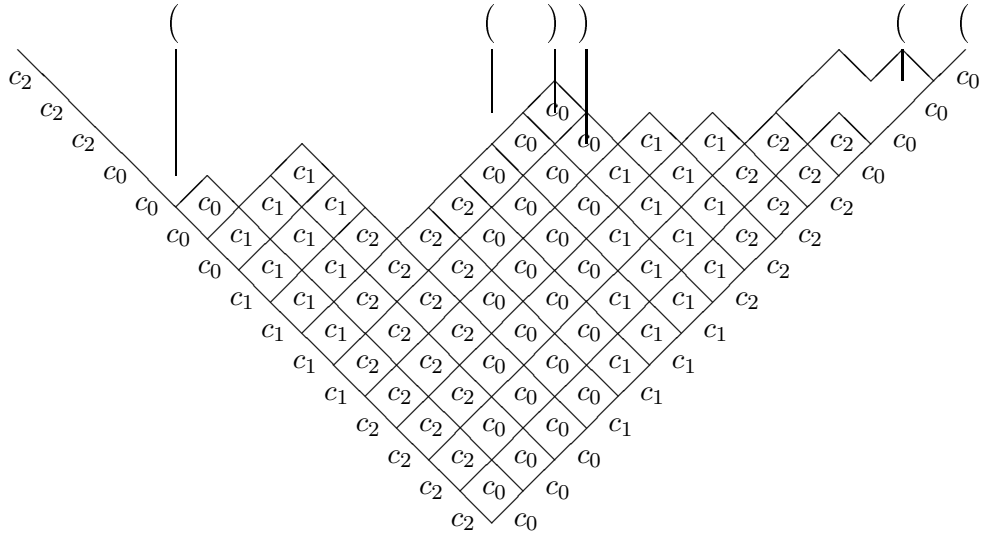
\begin{figure}
\begin{center}
\input{ribboncrystal}
\end{center}
\caption[Crystal moves on a partition]{Crystal moves on a partition. This example shows the action of $f_0$ on a level $4$ crystal for $\asl_3$. Color the boxes of a partition with $n$ colors $c_0, c_1, \ldots c_{n-1}$, where all boxes above position $k$ on the horizontal axis are colored $c_s$ for $\displaystyle s \equiv  \huge{ \lfloor} k /l \huge{\rfloor}$ modulo $n$. To act by $f_i$, place a $``("$ above the horizontal position $k$ if boxes in that position are colored $c_i$, and  you  can add an $\ell$-ribbon (see Section \ref{lotsastuff})  whose rightmost box is above $k$. Similarly, put a $``)"$ above each position $k$ where boxes are colored $c_i$ and you can remove an $\ell$-ribbon whose rightmost box is above $k$. $f_i$ acts by adding an $\ell$-ribbon whose rightmost box is below the first uncanceled $``("$ from the left, if possible, and sending that partition to $0$ if there is no uncanceled $``("$. $e_i$ acts by removing an $\ell$-ribbon whose rightmost box is below the first uncanceled $``)"$ from the right, if possible, and sending the partition to $0$ otherwise.  \label{ribboncrystal}}
\end{figure}

In our first model, the underlying set of the crystal consists of partitions. The crystal operators $f_i$ and $e_i$ act by adding and removing $\ell$-ribbons, where $\ell$ is the level of the crystal. These are described in Figure \ref{ribboncrystal}. In section \ref{highpart_section}, we prove that these rules define $\asl_n$ crystals, except possible when $n=2$. In that case, although we do not show that the whole model gives a crystal, we do prove Theorem \ref{highpart_intro}, so we still have a model for all irreducible integrable highest weight crystals. 

\setlength{\unitlength}{0.42cm}

\setlength{\unitlength}{0.5cm}
\begin{figure}
\begin{center}
\input{abacuscrystal}
\end{center}
\caption[Crystal moves on the abacus]{Crystal operators on the level $4$ abacus. Beads are represented by black circles, and empty spaces by white circles. Color the gaps between the columns of beads with $n$ colors, putting  $c_0$ at the origin, and $c_{[ i \mod n ]}$ in the $i^{th}$ gap, counting left to right. The operators $e_i$ and $f_i$ are then calculated as follows: Put a $``("$ every time a bead could move to the right across color $c_i$, and a $``)"$ every time a bead could move to the left across $c_i$. The brackets are ordered moving up each $c_i$ colored gap in turn from left to right. We group all the brackets corresponding to the same gap above that gap. $f_i$ moves the bead corresponding to the first uncanceled $``("$ from the left one place forward, if possible, and sends that element to 0 otherwise. Similarly, $e_i$ moves the bead corresponding the the first uncanceled $``)"$ from the right one space backwards, if possible, and sends the element to $0$ otherwise. \label{abacuscrystal}}
\end{figure}
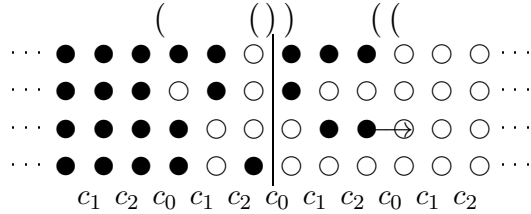

Our second model is based on the abacus used by James and Kerber (\cite{JK:1981}, Chapter 2.7). A level $\ell$ abacus configuration consists of $\ell$ strands of beads, with the possible positions of beads labeled by $\bz+1/2$, such that all but finitely many positions labeled by negative half integers of occupied, and all but finitely many position labeled by positive half integers are unoccupied. We color the abacus with $n$ colors, and define operators $e_i$ and $f_i$ depending on the coloring, as shown in Figure \ref{abacuscrystal}. The set of all abacus configurations becomes an $\asl_n$ crystal, at least when $n >2$. See Section \ref{abcryst} for more details. As explained by James and Kerber \cite{JK:1981}, there is a natural way to associate an abacus configuration to a given partitions, so these two models are closely related. 

The crystals one obtains from either of these models are reducible. However, using the abacus model, one can pick out a "highest" irreducible sub-crystal. To do this, we need to introduce a few definitions. These are made more precise in Section \ref{highpart_section}.

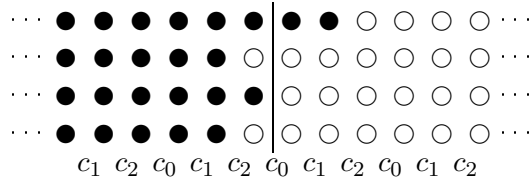
\begin{figure}
\begin{center}
\begin{picture}(20,4.5)
\input{compactification_intro}
\end{picture}
\end{center}
\caption[A compactification]{The compactification of the abacus configuration shown in Figure \ref{abacuscrystal}.  The configuration in Figure \ref{abacuscrystal} has weight 13 because you need to move a single bead one step to the left 13 times to reach the compactification. \label{compactification_intro} }
\end{figure}

\begin{Definition}
Let $\psi$ be an abacus configuration. The {\bf compactification} of $\psi$, denoted $\psi_{(0)}$, is the configuration obtained by pushing all the (black) beads to the left, using only finitely many moves, and not changing the row of any bead. For example, Figure \ref{compactification_intro} is the compactification of the configuration in Figure  \ref{abacuscrystal}.
\end{Definition}

\begin{Definition}
The {\bf weight} of an abacus configuration $\psi$ is the total number of times one bead needs to be pushed one step to the left before reaching the compactification of $\psi$.
\end{Definition}

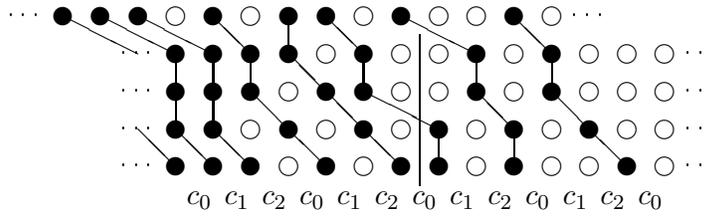
\begin{figure}
\begin{center}
\input{descendingabacus}
\end{center}
\caption[A descending abacus configuration]{A descending abacus configuration. An extra row has been added at the top, which is the same as the bottom row but shifted $n=3$ steps to the left. For every $k$, the $k^{th}$ black bead from the right in each row are connected with a line. The configuration is descending because all of these "descending strands" slope only up and to the left. 
This abacus configuration is not tight because each bead on the second strand can be shifted down one position, and they all remain to the right of the third strand. Here the first and fourth strands could also be shifted down, but the others cannot. If none of the strands can be shifted down in this way, the abacus configuration is called tight.
\label{descending_figure}}
\end{figure}

\begin{Definition}
Given a level $\ell$ abacus configuration $\psi$, put one more row at the top which is equal to the bottom row, but shifted $n$ steps to the left. For each $r \in \bn$, connect the $r$th bead in each row with a line, as in Figure \ref{descending_figure}. The configuration $\psi$ is called {\bf descending}  if all these lines slope only up and to the left.
\end{Definition}

\begin{Definition}
A descending abacus configuration is called {\bf tight} if one cannot "tighten" it by shifting all the beads in one descending strand down one position, staying entirely to the right of with the next strand. See Figure \ref{descending_figure}.
\end{Definition}

%\begin{figure}
%\begin{center}
%\begin{picture}(20,4)
%\input{abacus1_intro}
%\end{picture}
%\end{center}

%\caption[An abacus configuration]{A level $4$-abacus configuration for $\asl_3$. \label{abacus1_intro}}

%\end{figure}

Our first main result in Chapter \ref{abacus_chapter} is a construction of the crystal associated to any irreducible integrable highest weight $\asl_n$-module, of any level $\ell >0$:

\begin{Theorem} \label{highpart_intro}
The set of tight descending abacus configurations with a fixed compactification is an irreducible level $\ell$ crystal for $\asl_n$. Furthermore, every irreducible integrable highest weight $\asl_n$ crystal, of every level $\ell \geq 1$, occurs in this way. 
\end{Theorem}

\subsection{The cylindric plane partition model}
We now have a model for the crystal of any irreducible integrable highest weight representation of $\asl_n$. It is also useful to consider the larger crystal consisting of all descending abacus configurations with a given compactification, but which may not be tight. This decomposes as a union of infinitely many irreducible crystals, all of which are isomorphic. In fact, as discussed in Section \ref{rld}, it can be thought of as the crystal associated to an irreducible representation for $\agl_n$. It turns out that this set is naturally in bijection with the set of cylindric plane partitions, which we now define.

\begin{Definition} \label{cppd}
A cylindric plane partition of type $(n, \ell)$ is an array of non-negative integers $(\pi_{ij})$, defined for all sufficiently large $i, j \in \bz$, and satisfying:
\begin{enumerate}
\item \label{cppd1} If $\pi_{ij}$ is defined, then so is $\pi_{km}$ whenever $k \geq i$ and $m \geq j$.

\item \label{cppd2} $(\pi_{ij})$ is weakly decreasing in both $i$ and $j$. Furthermore, for all $i$, $\displaystyle \lim_{j \rightarrow \infty} \pi_{ij}=0$, and, for all $j$, $\displaystyle \lim_{i \rightarrow \infty} \pi_{ij}=0$. 

\item \label{cppd3} If $\pi_{ij}$ is defined, then $\pi_{ij} = \pi_{i+\ell, j-n}$.

The {\bf weight} of a cylindric plane partition is the sum of the entries over one period. The {\bf boundary} of a cylindric plane partition is the data of which $\pi_{ij}$ are defined.
\end{enumerate}
\end{Definition}

\begin{figure} \input{aplaincpp}
\caption[A cylindric plane partition]{A cylindric plane partition of type $(3,6)$. All empty squares are understood to contain $0$. The entries are weakly decreasing along both sets of diagonals. Note that the figure is periodic, with one period shown between the vertical lines. \label{aplaincpp}}
\end{figure}
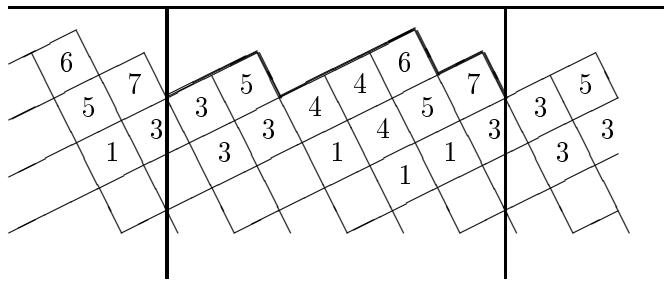

\begin{figure}

$$\mathfig{0.48}{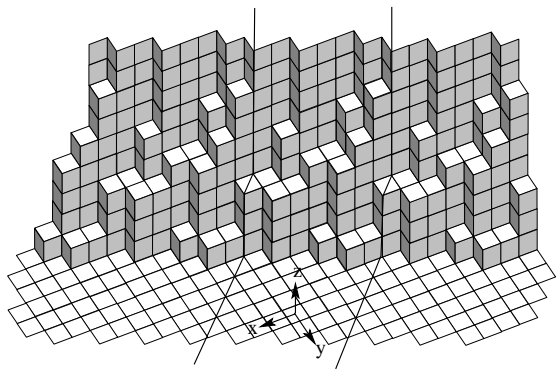}$$

\caption[The three dimensional representation of a cylindric plane partition]{The three dimensional representation of the cylindric plane partition shown in Figure \ref{aplaincpp}. Place a pile of boxes at each position, with the highest given by the entry in that position. The boundary is depicted by an infinitely high "wall". Note that the figure is periodic, with one period shown between the dark lines.    \label{3Dcpp_intro}} 
\end{figure}

Cylindric plane partition can be depicted by a semi-infinite array of integers, as in Figure \ref{aplaincpp}, or via a three dimensional interpretation, as in Figure \ref{3Dcpp_intro}.

\begin{Theorem} \label{abacus_cpp_intro} (see Section \ref{bpp})
There is a weight preserving bijection between descending abacus configurations with a given compactification $\psi_{(0)}$ and cylindric plane partitions with a certain boundary (depending on $\psi_{(0)}$).
\end{Theorem}

Thus the crystal structure on the set of descending abacus configurations transfers to a crystal structure on the set of cylindric plane partitions. We describe the actions of $e_i$ and $f_i$ on cylindric plane partitions explicitly in Section \ref{cpp_crystal}, although for the most part we find it easier to work directly with the descending abacus configurations. The following follows from Theorem \ref{abacus_cpp_intro}, once we have studied descending abacus configurations in more detail.

\begin{Theorem} \label{generating_intro} (see Section \ref{rcpp})
The generating function for cylindric plane partitions with a given boundary is given by the $q$-character of a certain irreducible representation of $\agl_n$.
\end{Theorem}

The highest weight of the irreducible representation involved is easily calculated from the boundary of the cylindric plane partitions, as we explain in Section \ref{rcpp}. The generating function for cylindric plane partitions has previously been calculated by Borodin \cite{Borodin:2006}. His answer looked very different, but we directly show they agree.

\subsection{Rank level duality}

A cylindric plane partition of type $(\ell, n)$ can be read either clockwise on counterclockwise (left to right or right to left in Figure \ref{aplaincpp}). Thus Theorem \ref{generating_intro} gives two ways to calculate the generating function. In one, it is the $q$-character (or principally graded character) of a certain irreducible level $\ell$ representation of $\agl_n$. In the other it is the $q$-character of a certain irreducible level $n$ representation of $\agl_\ell$. Thus we recover the following rank-level duality result, first observed by Frenkel in \cite{IFrenkel:1982}. A more precise statement is given in Section \ref{rld}.

\begin{Theorem} (Frenkel)
Let $\Lambda$ be a level $\ell$ dominant integral weight of $\agl_n$, and $W_\Lambda$ the corresponding irreducible representation of $\agl_n$. Then there is a corresponding level $n$ dominant integral weight $\Lambda'$ of $\agl_\ell$ such that
\begin{equation}
\dim_q(W_\Lambda)= \dim_q(W_{\Lambda'}).
\end{equation}
\end{Theorem}

\noindent As we explain in Section \ref{rld}, $\Lambda'$ is easily calculated from $\Lambda$. Although this result is not new, our proof is very different from that given by Frenkel.

\subsection{Notable connections to known models}

In Section \ref{like_kashiwara} we relate our realization of the irreducible integrable highest weight crystals of $\asl_n$ to the Kyoto path model developed by Kashiwara et. al. in \cite{KKMMNN1} and \cite{KKMMNN2} (see \cite{Hong&Kang:2000} for a more recent explanation). This is done by exhibiting an explicit crystal isomorphism between the highest irreducible component of the abacus model and the Kyoto path model for a particular perfect crystal and ground state path. The "descending strands" in our model (see Figure \ref{descending_figure}) correspond to the factors in the semi-infinite tensor product of perfect crystals used in the Kyoto path model. 

We would also like to mention a 1991 paper by Jimbo, Misra, Miwa and Okado \cite{JMMO:1991} which contains some results relevant to Chapter \ref{abacus_chapter}. In particular, they present a similar realization for the crystal of any irreducible integrable highest weight representation of $\asl_n$.

%\end{document}

\chapter{Background and notation} \label{background_chapter}

In this chapter we present an assortment results which we will use later on. We suggest the casual reader skip to Chapters \ref{commutor_chapter} and \ref{abacus_chapter}, and refer to the background material as needed. 

\section{Conventions} \label{notation}
We must first fix some notation. For the most part we follow conventions from \cite{CP}.

$\bullet$ $\g$ is a complex simple Lie algebra with Cartan algebra $ \mathfrak{h} $, and $A = (a_{ij})_{i,j \in I}$ is its Cartan matrix. 

$\bullet$ $ \langle \cdot , \cdot \rangle $ denotes the paring between $ \mathfrak{h} $ and $ \mathfrak{h}^\star $ and $ ( \cdot , \cdot) $ denotes the usual symmetric bilinear form on either $ \mathfrak{h}$ or $ \mathfrak{h}^\star $.  Fix the usual bases $ \alpha_i $ for $ \mathfrak{h}^\star $ and $ H_i $ for $\mathfrak{h}$, and recall that $ \langle H_i, \alpha_j \rangle = a_{ij} $.  

$\bullet$ $ d_i = (\alpha_i, \alpha_i)/2 $, so that $ (H_i, H_j) = d_j^{-1} a_{ij} $.  Let $ B $ denote the matrix $ (d_j^{-1} a_{ij}) $. 

$\bullet$ $ q_i = q^{d_i} $.  

$\bullet$ $\rho$ is the weight satisfying $(\alpha_i, \rho)=d_i$ for all $i$.

$\bullet$ $ H_\rho $ is the element of $ \mathfrak{h} $ such that $ \langle \alpha_i, H_\rho \rangle = d_i = (\alpha_i, \rho) $ for all $ i $.

$\bullet$ $\theta$ is the diagram automorphism such that $w_0 (\alpha_i) = - \alpha_{\theta(i)},$ where $w_0$ is the longest element in the Weyl group.

$\bullet$ $U_q(\g)$ is the quantized universal enveloping algebra associated to $\g$, generated over $\mathbb{C}(q)$ by $E_i$, $F_i$ for all  $i \in I$, and $K_w$ for $w$ in the co-weight lattice of $\g$. As usual, let $K_i= K_{H_i}.$ We use conventions as in \cite{CP}. For convenience, we recall the exact formula for the coproduct:
\begin{equation} \label{coproduct}
\begin{cases}
\Delta{E_i} & = E_i \otimes K_i + 1 \otimes E_i \\
\Delta{F_i} &= F_i \otimes 1 + K_i^{-1} \otimes F_i \\
\Delta{K_i} &= K_i \otimes K_i
\end{cases}
\end{equation}

$\bullet$ In Chapter \ref{commutor_chapter}, we sometimes need to adjoin a fixed $k^{th}$ root of unity to the base field $\bc(q)$, which we denote by $q^{1/k}$. The integer $k$ depends on $\g$, and always divides twice the dual Coxeter number. 

$\bullet$ $[n]= \frac{q^n - q^{-n}}{q-q^{-1}},$ and $X^{(n)} = \frac{X^n}{[n][n-1] \cdots [2]}.$

$\bullet$ $V_\lambda$ is the irreducible representation of $U_q(\g)$ with highest weight $\lambda$. 

$\bullet$ $B_\lambda$ is a fixed global basis for $V_\lambda$, in the sense of Kashiwara (see \cite{Kashiwara:1991}). $b_\lambda$ and $b_\lambda^{\low}$ are the highest weight and lowest weight elements of $B_\lambda$ respectively.

$\bullet$ We reserve $\Lambda$ to be an integrable highest weight of $\asl_n$.

\section{The quantum Weyl group} \label{quantum_weyl_group}

Following Lusztig \cite[Part VI]{Lusztig:1993} and \cite[Section 5]{L2}, we introduce an action of the braid group of type $ \g $ on any $V_\lambda$, which is know as the quantum Weyl group. In order to be precise, we first introduce a completion $\widetilde{U_q(\g)}$ of $U_q(\g)$. We then define the quantum Weyl group as a subgroup of the invertible elements of $\widetilde{U_q(\g)}$. We will be most interested in the action of the braid group element $T_{w_0}$ corresponding to the longest element in the Weyl group. 

\subsection{The completion $\widetilde{U_q(\g)}$} \label{comp_section}

We will be working in the completion  $\widetilde{U_q(\g)}$ of $U_q(\g)$ with respect to the weak topology generated by all matrix elements of finite dimensional representations. This section includes two equivalent explicit definitions of $\widetilde{U_q(\g)}$ (Definition \ref{hatdef} and Corollary \ref{hatdef2}), as well as some basic results about its structure. Most importantly, we show that $\widetilde{U_q(\g)}$ is isomorphic to the direct product of the endomorphism rings of all $V_\lambda$. Thus an element of $\widetilde{U_q(\g)}$ is equivalent to a choice of $x \in \End ({V_\lambda})$ for each $\lambda \in P_+$. 

\begin{Definition} \label{hatdef} Let $R$ be the ring consisting of series $\sum_{k=1}^\infty X_k$, where each $X_k \in U_q(\g)$ and, for any fixed $\lambda$, $X_k \cdot V_\lambda = 0$ for all but finitely many $k$. Notice that there is a well defined action of $R$ on any $V_\lambda$. Let $I$ be the two sided ideal in  $R$ consisting of elements which act as zero on all $V_\lambda$. Then  $\widetilde{U_q(\g)}$ is defined to be $R/I$.
\end{Definition}

\begin{Comment}
This is equivalent to the completion with respect to the topology mentioned above, since $U_q(\g)$ is semi simple, so the set of matrix elements of finite dimensional representations is point-separating for $U_q(\g)$. In particular the natural map of $U_q(\g)$ to $\widetilde{U_q(\g)}$ is an embedding. 
\end{Comment}

This completion has a simple description as follows:

\begin{Theorem} \label{comp}
$\widetilde{U_q(\g)}$ is isomorphic as an algebra to $\displaystyle \prod_{\lambda \in P_+} \End_{\bc(q)} (V_\lambda)$.
\end{Theorem}

Before proving Theorem \ref{comp} we will need two technical lemmas.

\begin{Lemma} \label{isp}
There is an element $p_\lambda \in U_q(\g)$ such that 
\begin{enumerate}
\item $p_\lambda (v_\lambda) = v_\lambda$

\item For any $\mu \neq \lambda$, $p_\lambda$ sends the $\mu$ weight space of $V_\lambda$ to 0.

\item $p_\lambda V_\mu = 0$ unless $\langle \mu-\lambda, \rho^\vee \rangle > 0$ or $\mu = \lambda$.
\end{enumerate}
\end{Lemma}

\begin{proof}
Fix a lowest weight vector $v_\lambda^{\text{low}} \in V_\lambda$.
$V_\lambda$ is a quotient of $U_q^-(\g) \cdot v_\lambda$, so we can choose some $F \in U_q^-(\g)$ such that $Fv_\lambda = v_\lambda^{\text{low}}$. Similarly, we can choose some $E \in U^+_q(\g)$ such that $E v_\lambda^{\text{low}} = v_\lambda$. Then $p':= EF$ clearly satisfies the first two conditions.

For each $i \in I$, let $R_i = E_i^{( \langle \lambda, \alpha_i^\vee \rangle )} F_i^{( \langle \lambda, \alpha_i^\vee \rangle )}$. Let
\begin{equation}
p_\lambda= \left(\prod_{i \in I} R_i\right) p',
\end{equation}
where the product is taken in any order. It is straightforward to see that this element satisfies the desired properties.
\end{proof}

\begin{Lemma} \label{isend}
Let $I_\lambda$ be the kernel of the action of $U_q(\g)$ on $V_\lambda$. Then
$U_q(\g) /I_\lambda $ is isomorphic to $\End_{\bc (q) } V_\lambda$.
\end{Lemma}

\begin{proof}
Let $d = \dim (V_\lambda)$.
Using the PBW basis in the $E$s, there is a $d$ dimensional subspace $\mathcal{F}$ of elements in $U_q^+(\g)$ that act non-trivially on $V_\lambda,$ and in fact such that  $p_\lambda \mathcal{F}$ is still $d$ dimensional, where  $p_\lambda$ is as in Lemma \ref{isp}. One can tensor this space with the PBW operators from $U_q^-(\g)$ to get a $d^2$ dimensional subspace of $U_q(\g)$ that acts non-trivially on $V_\lambda$. The result follows.
\end{proof}

\begin{proof}[{\bf Proof of Theorem \ref{comp}}]
Using Lemmas \ref{isp} and \ref{isend}, we can realize any endomorphism of $V_\lambda$ using an element of $U_q(\g)$ that kills $V_\mu$ unless $\langle \mu-\lambda, \rho^\vee \rangle >0$, or $\mu = \lambda$. The result follows.
\end{proof}

We include the following result to show how our definition of $\widetilde{U_q(\g)}$ relates to other completions that appear in the literature. This could also be taken as the definition of $\widetilde{U_q(\g)}$.

\begin{Corollary} \label{hatdef2}
Let each $\lambda \in P_+$, let $I_\lambda$ be the two sided ideal of $U_q(\g)$ generated by all $E_i^{\langle \lambda, \alpha_i^\vee \rangle}$ and $F_i^{\langle \lambda, \alpha_i^\vee \rangle}$. Let $\displaystyle U_q''(\g)= \lim_{\leftarrow} U_q(\g)/I_\lambda$, using the partial order on weights where $\mu \leq \lambda$ if and only if $\lambda - \mu \in P_+$. $U_q''(\g)$ acts in a well defined way on any finite dimensional module, so there is a map $U_q''(\g) \rightarrow \widetilde{U_q(\g)}$. This is an isomorphism.
\end{Corollary}

\begin{proof} \label{whdef2}
The same argument as we used to prove Theorem \ref{comp} shows that the image is $\displaystyle \prod_{\lambda \in P_+} \End(V_\lambda)$, which is all of $ \widetilde{U_q(\g)}$. The map is injective by the definition of $U''_q(\g)$.
\end{proof}

It will also be nessisary for us to consider the following completion of $U_q(\g) \otimes U_q(\g)$:

\begin{Definition} \label{bigbig_completion_def}
Let $\widetilde{U_q(\g) \otimes U_q(\g)}$ be the completion of  $U_q(\g) \otimes U_q(\g)$ in the weak topology generated all matrix elements it's action on $V \otimes W$ for all pairs of finite dimensional representations $V$ and $W$ of $U_q(\g)$.
\end{Definition}

\begin{Comment}
 The completion $\widetilde{U_q(\g)}$ is related to the algebra $\dot{U}$ from \cite[Chapter 23]{Lusztig:1993} as follows. $\dot{U}$ acts in a well defined way on each irreducible representation $V_\lambda$, and no non-zero element of $\dot{U}$ acts as zero on every $V_\lambda$. Hence $\dot{U}$ naturally embeds in $\widetilde{U_q(\g)}$. There is a canonical basis $\dot{B}$ for $\dot{U}$. All but finitely many elements of $\dot{B}$ act as zero on any given $V_\lambda$ (see \cite{Lusztig:1993} Remark 25.2.4 and Section 23.1.2), so the space of all formal (infinite) linear combinations of elements of $\dot{B}$ also maps to $\widetilde{U_q(\g)}$. This map is bijective, and so $\widetilde{U_q(\g)}$ is naturally identified with the space of all formal linear combinations of elements of $\dot{B}$.
 \end{Comment}

It is clear the  $\widetilde{U_q(\g)}$ has the structure of a ring, and that it acts in a well defined way on finite representations. It also has a well defined topological coalgebra structure, with the coproduct of $u$ defined by the action of an element $u$ on a tensor product $V \otimes W$. This is only a topological coproduct because it maps $ \widetilde{U_q(\g)}$ into $\displaystyle \prod_{\lambda, \mu} \End_{\bc(q)} V_\lambda \otimes \End_{\bc(q)} V_\mu$, which can be though of as a completion of  $\displaystyle \prod_{\lambda} \End_{\bc(q)} V_\lambda \otimes  \prod_{\mu}  \End_{\bc(q)} V_\mu$. The restriction of this coproduct to $U_q(\g)$ agrees with the normal coproduct so, since $U_q(\g)$ is a dense subalgebra of $\widetilde{U_q(\g)}$, we see that $\widetilde{U_q(\g)}$ is a topological Hopf algebra. 

We will need to consider the group of invertible elements of $\widetilde{U_q(\g)}$ acting on $\widetilde{U_q(\g)}$ by conjugation. This action preserves the algebra structure of $\widetilde{U_q(\g)}$, but does not preserve the coproduct.

\begin{Definition} \label{Aactdef}
Let $X$ be an invertible element in $\widetilde{U_q(\g)}$. Define $C_X$ (conjugation by $ X $) to be the algebra automorphism of $\widetilde{U_q(\g)}$ defined by $u \rightarrow X u X^{-1}$.
\end{Definition}

\begin{Comment}
We caution the reader that $C_X$ is not that Hopf theoretic adjoint action of $X$, as defined in, for example, \cite{CP}. 
\end{Comment}

\begin{Comment} \label{thediagram}
For any invertible $X \in \widetilde{U_q(\g)},$ the action of $X $  on representations is compatible with the automorphism $ C_X $ in the sense that, for any representation $V$, the following diagram commutes:
\begin{equation}
\xymatrix{
V  \ar@(dl,dr) \ar@/ /[rrr]^{X} &&& V \ar@(dl,dr) \\
\widetilde{U_q(\g)}  \ar@/ /[rrr]^{C_X} &&& \widetilde{U_q(\g)}. \\
}
\end{equation} 
In general, $C_X$ does not preserve the subalgebra $U_q(\g)$ of $\widetilde{U_q(\g)}$, although it does in all cases we consider here.   
\end{Comment}

\subsection{Definition of the quantum Weyl group}

We first define the action of the generators $T_i$. Our conventions are such that $T_i$ is $T_{i,-1}^{"} = T_{i,1}'^{-1}$ in the notation from \cite{Lusztig:1993}. 
 \begin{Definition} \label{defti} (see \cite[5.2.1]{Lusztig:1993})
 $T_i$ is the element of  $\widetilde{U_q(\g)}$ that acts on a weight vector $v$  by:
 \begin{equation}
 T_i(v)= \sum_{\begin{array}{c} a,b,c \geq 0 \\ a-b+c=(\wt(v), \alpha_i) \end{array}} (-1)^b q_i^{ac-b}E_i^{(a)} F_i^{(b)} E_i^{(c)} v.
 \end{equation}
\end{Definition}

By \cite[Theorem 39.4.3]{Lusztig:1993}, these $T_i$ generate an action of the braid group on each $V_\lambda$, and thus a map from the braid group to $\widetilde{U_q(\g)}$. This realization of the braid group is often referred to as the quantum Weyl group. It is related to the classical Weyl group by the fact that, for any weight vector $v \in V$, $\wt(T_i(v))= s_i(\wt(v))$.

\begin{Theorem} [{see \cite[Theorem 8.1.2]{CP}  or \cite[Section 37.1.3]{Lusztig:1993}}]
The conjugation action of the braid group on $\widetilde{U_q(g)}$  (see Definition \ref{Aactdef}) preserves the subalgebra $U_q(\g)$, and is defined on generators by:
\begin{equation} \label{defT}
\begin{cases}
C_{ T_i} (E_i)= -F_iK_i \\
C_{T_i} (F_i)= -K_i^{-1} E_i \\
C_{T_i} (K_H)= K_{s_i(H)} \\
C_{T_i} (E_j) = \sum_{r=0}^{-a_{ij}} (-1)^{r-a_{ij}} K_i^{-r} E_i^{(-a_{ij}-r)} E_j E_i^{(r)} \mbox{  if  } i \neq j \\
C_{T_i} ( F_j) = \sum_{r=0}^{-a_{ij}} (-1)^{r-a_{ij}} K_i^{r} F_i^{(r)} F_j F_i^{(-a_{ij}-r)} \mbox{  if  } i \neq j.
\end{cases}
\end{equation}
\end{Theorem}

Fix some $w$ in the Weyl group $W$, and a reduced decomposition of $w$ into simple reflection $w= s_{i_1} \cdots s_{i_k}$. By \cite[Section 2.1.2]{Lusztig:1993}, the element $T_w \in \widetilde{U_q(\g)} $ defined by
\begin{equation} \label{braidw}
T_w:= T_{i_1} \cdots T_{i_k}
\end{equation}
is independent of the reduced decomposition. Furthermore, the following holds.

\begin{Lemma} [see \cite{CP} Proposition 8.1.6]
Let $w \in W$ be such that $w(\alpha_i)= \alpha_j$. Then $C_{T_w} (E_i)= E_j$.
\end{Lemma}

\subsection{The action of $T_{w_0}$}

Let $w_0$ be the longest element of the Weyl group, and $T_{w_0}$ the corresponding element of the braid group given by Equation (\ref{braidw}).

\begin{Lemma} \label{stillE}
The action of $C_{T_{w_0}}$ on $U_q(\g)$ is given by
\begin{equation}
\begin{cases}
C_{T_{w_0}} (E_i) = -F_{\theta(i)} K_{\theta(i)} \\
C_{T_{w_0}} (F_i) = -K_{\theta(i)}^{-1} E_{\theta(i)} \\
C_{T_{w_0}} (K_H) = K_{w_0(H)}, \mbox{ so that } C_{T_{w_0}} (K_i)  = K_{\theta(i)}^{-1}
\end{cases}
\end{equation}
\end{Lemma}

\begin{proof}
Fix $i$. Then $T_{w_0}$ can be written as $T_{w_0} = T_{\theta(i)} T_w$ for some $w$ in the Weyl group. By the definition of $\theta$,  $C_{T_{w_0}} (E_i)$ is in the weight space $-\alpha_{\theta(i)}$. It follows that $C_{T_w}(E_i)$ is in the weight space $\alpha_{\theta(i)}$. Hence by Lemma \ref{stillE}, $C_{T_w}(E_i)= E_{\theta(i)}$. Therefore, by (\ref{defT}), $C_{T_{w_0}} (E_i)=  -F_{\theta(i)} K_{\theta(i)} $, as required. A similar proof works for $F_i$.  The action on $ K_H $ is straightforward.
\end{proof}

\begin{Comment}
Note that $ C_{T_{w_0}} $ is not a coalgebra anti-automorphism, so we cannot use $ T_{w_0} $ to construct a commutativity constraint in the manner of Section \ref{autint_intro}.  We will first need to correct $ T_{w_0} $.  There are essentially two natural ways of doing this --- one leads to the standard braiding and the other to Drinfeld's coboundary structure.
\end{Comment}

We now understand the action of $C_{T_{w_0}}$ on $U_q(\g)$. We also need to understand how $T_{w_0}$ acts on any finite dimensional representation. By Comment \ref{ximake} it is sufficient to know $T_{w_0} (v_\lambda)$ for each $\lambda \in P_+$. 

\begin{Lemma} \label{th:Ti}
Let $ V $ be any representation, and $ v \in V $a weight vector such that $ E_i \cdot v = 0 $.  Then $ T_i(v) =  (-1)^n q^{d_i n} F_i^{(n)} v $, where $ n = \langle \wt(v), \alpha_i^\vee \rangle $. 
\end{Lemma}

\begin{proof}
Fix $v \in V$ with $E_i(v)=0$, and let $n= \langle \wt(v), \alpha_i^\vee \rangle $. It follows from $U_q(sl_2)$ representation theory that $F_i^{n+1} (v)=0$. The lemma then follows directly from the definition of $ T_i $ (Definition \ref{defti}).
\end{proof}

The following can be found in  \cite[Lemma 39.1.2]{Lusztig:1993}  recalling that our $T_i$ is equal to $T_{i,1}'^{-1}$ in the notation from that book, although we find it convenient to include a proof. 

\begin{Proposition}  \label{jjthings}
Let $ w= s_{i_1} \cdots s_{i_\ell} $ be a reduced word. For each $1 \leq k \leq \ell$, the following statements hold.
\begin{enumerate}

\item \label{jj1} $E_{i_{k+1}} T_{i_k} \cdots T_{i_1} (v_\lambda) = 0$.

\item \label{jj2} $T_{i_k} \cdots T_{i_1} (v_\lambda) = (-1)^{n_1+ \cdots +n_k} q^{d_{i_1} n_1 + \dots + d_{i_k} n_k} F_{i_k}^{(n_k)} \cdots F_{i_1}^{(n_1)} v_\lambda$,\\ where $n_j = \langle  s_{i_1} \cdots s_{i_{j-1}}\alpha_{i_j}^\vee, \lambda \rangle $. 
\end{enumerate}
\end{Proposition}

\begin{proof}
Note that $\wt(E_{i_{k+1}} T_{i_k} \cdots T_{i_1} (v_\lambda)) = s_{i_{k}} \cdots s_{i_1} \lambda + \alpha_{i+1}$, so it suffices to show that the dimension of the $ s_{i_{k}} \cdots s_{i_1} \lambda + \alpha_{i_{k+1}}$ weight space in $V_\lambda$ is zero. The dimensions of weight spaces are invariant under the Weyl group, so we may act by $s_{i_1} \cdots s_{i_k}$ and instead show that the $\lambda +s_{i_1} \cdots s_{i_k} \alpha_{i_{k+1}} $ weight space of $V_\lambda$ is zero. But $s_{i_1} \cdots s_{i_k} s_{i_{k+1}}$ is a reduced word in the Weyl group, which implies that $s_{i_1} \cdots s_{i_k} \alpha_{i_{k+1}} $ is a positive root. Since $\lambda$ is the highest weight of $V_\lambda$, part  (\ref{jj1}) follows.
  
Part (\ref{jj2}) Follows by repeated use of (\ref{jj1}) and Lemma \ref{th:Ti}.
\end{proof}

\begin{Definition} \label{lowdef}
Fix a highest weight vector $v_\lambda \in V_\lambda$. Define the corresponding lowest weight vector $ v_{\lambda}^{\text{low}} \in V_\lambda $ by 
\begin{equation} \label{twl} T_{w_0} v_\lambda = (-1)^{\langle 2 \lambda, \rho^\vee \rangle} q^{(2\lambda, \rho)} v_{\lambda}^{\text{low}}. \end{equation}
By Proposition \ref{jjthings} part (\ref{jj2}), this is equivalent to defining 
\begin{equation} \label{vlambdalow}
v_{\lambda}^{\text{low}} = F_{i_m}^{(n_m)} \cdots F_{i_1}^{(n_1)} v_\lambda,
\end{equation}
where $ w_0 = s_{i_1} \cdots s_{i_m} $ is any reduced expression for the longest element in the Weyl group, and $n_j = \langle  s_{i_1} \cdots s_{i_{j-1}}\alpha_{i_j}^\vee, \lambda \rangle. $ \end{Definition}

\begin{Comment}
It follows from Proposition \ref{bottotop} below that $v_\lambda$ and $v_\lambda^{\low}$ are also related by
\begin{equation}
T_{w_0} v_\lambda^{\low} =v_\lambda.
\end{equation}
This is a simpler statement that Equation (\ref{twl}), although is somewhat more difficult to prove directly.
\end{Comment}

\section{Braided categories, coboundary categories and $R$-matrices} \label{coboundary_section}

We briefly review the definitions of braided and coboundary categories, and how $Rep(U_q(\g))$ can be given such structures. 

\begin{Definition} \label{cat_ax}
Let $\mathcal{C}$ be a monoidal category with  a natural system of isomorphisms $\sigma= \{ \sigma_{A,B} : A \otimes B \rightarrow B \otimes A \}$ for all $A,B \in \mathcal{C}$.

\begin{enumerate}
\item \label{square_ax} We say $\sigma$ {\bf squares to one} if for all $A$ and $B$, $\sigma_{A,B} \circ \sigma_{B,A}$ is the identity.

\item \label{braided_ax} We say $\sigma$ satisfies the {\bf hexagon axiom} if for all $A,B,C \in \mathcal{C}$, the following diagrams commute
 \begin{equation*}
\xymatrix{
A \otimes B \otimes C  \ar@/ /[drr]_{\Id \otimes \sigma_{B,C}} \ar@/ /[rrrr]^{\sigma_{A \otimes B, C}} &&&& C \otimes A \otimes B  \\
&& A \otimes C \otimes B  \ar@/ /[urr]_{\sigma_{A,C} \otimes \Id}  \\
}
\end{equation*} 

 \begin{equation*}
\xymatrix{
A \otimes B \otimes C  \ar@/ /[drr]_{\sigma_{A,B} \otimes \Id} \ar@/ /[rrrr]^{\sigma_{A, B \otimes C}} &&&& B \otimes C \otimes A  \\
&& B \otimes A \otimes C  \ar@/ /[urr]_{\Id \otimes \sigma_{A,C}}  \\
}
\end{equation*}

Note: These diagrams would be hexagons if we included the associator maps.

\item \label{leaf_ax} We say $\sigma$ satisfies the {\bf leaf axiom} if for all $A,B,C \in \mathcal{C}$ the following diagram commutes:
 \begin{equation*}
\xymatrix{
A \otimes B \otimes C  \ar@/ /[d]^{\sigma_{A, B \otimes C}} \ar@/ /[rrr]^{\small{\Id} \otimes \sigma_{B,C}} &&& A \otimes C \otimes B \ar@/ /[d]^{\sigma_{A, C \otimes B}} \\
B \otimes C \otimes A  \ar@/ /[rrr]^{\sigma_{B,C} \otimes \Id} &&&   C \otimes B \otimes A. \\
}
\end{equation*} 
\end{enumerate}
\end{Definition}

\begin{Definition} \label{sym_cat_def}
A {\bf symmetric monoidal category} is a pair $(\mathcal{C}, \sigma)$, where $\mathcal{C}$ is a monoidal category and $\sigma = \{ \sigma_{A,B} : A, B \in \mathcal{C} \}$ satisfies Definition \ref{cat_ax} parts (\ref{square_ax}) and (\ref{braided_ax}).
\end{Definition}

\begin{Definition} \label{braided_cat_def}
A {\bf braided monoidal category} is a pair $(\mathcal{C}, \sigma)$, where $\mathcal{C}$ is a monoidal category and $\sigma = \{ \sigma_{A,B} : A, B \in \mathcal{C} \}$ satisfies Definition \ref{cat_ax} part (\ref{braided_ax}).
\end{Definition}

\begin{Definition} \label{coboundary_cat_def}
A {\bf coboundary monoidal category} is a pair  $(\mathcal{C}, \sigma)$, where $\mathcal{C}$ is a monoidal category and $\sigma = \{ \sigma_{A,B} : A, B \in \mathcal{C} \}$ satisfies Definition \ref{cat_ax} parts (\ref{square_ax}) and (\ref{leaf_ax}).
\end{Definition}

\begin{Comment}
Definition \ref{cat_ax} Part (\ref{leaf_ax}) follows from parts \ref{square_ax} and \ref{braided_ax}, so the definitions of braided monoidal categories and coboundary categories are both relaxations of the conditions for a symmetric monoidal category.
\end{Comment}

As we have mentioned, $Rep(U_q(\g))$ can be given the structure of a braided monoidal category. This is usually done by constructing on element $R$ in a certain completion of $U_q(\g) \otimes U_q(\g)$, known as the universal R matrix (see for example \cite{CP}). 
\begin{Definition}
A universal $R$-matrix for $U_q(\g)$ is an element $R$ in a completion of $U_q(\g) \otimes U_q(\g)$ which acts in a well defined way on  $V \otimes W$ for any $V$ and $W$ in $Rep(U_q(\g))$, and such that the system of maps $\sigma_{V,W} := \Flip \circ R: V \otimes W \rightarrow W \otimes V$ gives $Rep(U_q(\g))$ the structure of a braided monoidal category.
\end{Definition}
The universal $R$ matrix is not truly unique. However, it exists, and there is a well studied standard choice.

In \cite[Section 3]{Drinfeld:1990}, Drinfeld presents a way to modify $\sigma^{br}$ to give $Rep(U_q(\g))$ the structure of a coboundary category. We now review the construction of this "unitarized" system of isomorphisms, which we denote $\sigma^D$.

Consider the "ribbon element" $Q$. This belongs to a certain completion of $U_q(\g)$, and is uniquely determined by the fact that it acts on each irreducible representation $ V_\lambda$ as multiplication by $ q^{(\lambda, \lambda + 2 \rho)} $ (see for example \cite{BK}).  
The following is well known.

\begin{Proposition} \label{ror} For and $V,W \in Rep(U_q(\g))$, 
\begin{equation} \sigma^{br}_{W,V} \circ \sigma^{br}_{V,W} = (Q^{-1} \otimes Q^{-1}) \Delta(Q), \end{equation}
where $\Delta(Q)$ denotes the action of $Q$ on $V \otimes W$.
\end{Proposition}

The element $ Q $ is central, and $Q$ admits a central square root, denoted by  $ Q^{1/2} $, which  acts on $V_\lambda$ as multiplication by the constant $ q^{(\lambda, \lambda)/2 + (\lambda, \rho) }$.  

\begin{Definition}
$\sigma^D= \{ \sigma_{V,W}^D \}$ where
\begin{equation}
\sigma^D_{V,W}= \sigma^{br} (Q^{1/2} \otimes Q^{1/2}) \Delta(Q^{-1/2}).
\end{equation}
\end{Definition}
It was shown by Drinfeld in \cite[Section 3]{Drinfeld:1990} that $(Rep(U_q(\g)), \sigma^D)$ is a coboundary category. Actually, Drinfeld defines $\sigma^D$ using the unitarized $R$-matrix
\begin{equation} \label{QR}
\bar{R} = R (Q^{1/2} \otimes Q^{1/2}) \Delta(Q^{-1/2}) =  (Q^{1/2} \otimes Q^{1/2})R \Delta(Q^{-1/2}) ,
\end{equation}
so that $\sigma^D = swap \circ \bar{R}$.

\section{Crystals and crystal bases} \label{crystal_section}
In this section, we introduce crystal bases and abstract crystals, and discuss their structure.  We refer the reader to \cite{Kashiwara:1991} for a more detailed overview of these topics.  We will also use some results from \cite{CP}. Unfortunately, the conventions in \cite{CP} and \cite{Kashiwara:1991} do not agree. Thus, in order that our conventions remain consistent, we have been forced to modify certain statements from \cite{Kashiwara:1991}. In particular, we use crystal basis at $q=\infty$ rather then at $q=0$.

\subsection{Crystal bases}

\begin{Definition}
Let $\Aa_\infty= \mathbb{C}[q]_\infty$ be the algebra of rational functions in one variable $q$ over $\bc$ whose denominators are not divisible by $q$.
\end{Definition}

\begin{Definition} \label{Kop} Fix a finite dimensional representation $V$ of $\g$, and $i \in I$.  Define the Kashiwara operators $ \tilde{F}_i, \tilde{E}_i : V \rightarrow V $ by linearly extending the following definition
\begin{equation}
\begin{cases}
\tilde{F}_i(F_i^{(n)}(v)) = F_i^{(n+1)} (v) \\
\tilde{E}_i (F_i^{(n)}(v))= F_i^{(n-1)} (v).
\end{cases}
\end{equation}
for all $ v \in V $ such that  $E_i (v)=0$. 
\end{Definition}

\begin{Comment} \label{symt}
It follows from the representation theory of $\text{sl}_2$ that  $ \tilde{F}_i$ and $ \tilde{E}_i $ can also be defined by linearly extending 
\begin{equation}
\begin{cases}
\tilde{E}_i(E_i^{(n)}(v)) = E_i^{(n+1)} (v) \\
\tilde{F}_i (E_i^{(n)}(v))= E_i^{(n-1)} (v).
\end{cases}
\end{equation}
for all $ v \in V $ such that  $F_i (v)=0$. Thus the operators are symmetric under interchanging the roles of $E_i$ and $F_i$, even if the definition does not appear to be. 
\end{Comment}

\begin{Definition}
A crystal basis of a representation $ V$ is a pair $(\LL, B) $, where $ \LL $ is an $\Aa_\infty$-lattice  of $ V$ and $B$ is a basis for $ \LL/q^{-1} \LL$, such that
\begin{enumerate}
\item $\LL $ and $ B $ are compatible with the weight decomposition of $ V $.
\item $\LL $ is invariant under the Kashiwara operators and $ B \cup 0 $ is invariant under their residues $ e_i := \tilde{E}_i^\modq, f_i := \tilde{F}_i^\modq : \LL/q^{-1}\LL \rightarrow \LL/q^{-1} \LL$.
\item For any $ b, b' \in B $, we have $ e_i b = b' $ if and only if $ f_i b' = b $.
\end{enumerate}
\end{Definition}

The following three theorems of Kashiwara are crucial to us.

\begin{Theorem}[\cite{Kashiwara:1991}, Theorem 1]
Let $ V, W $ be representations with crystal bases $ (\LL, A)$ and $(\mathcal{M}, B) $ respectively.  Then $( \LL \otimes \mathcal{M}, A \otimes B) $ is  a
crystal basis of $ V \otimes W $.
\end{Theorem}

\begin{Theorem}[\cite{Kashiwara:1991}, Theorem 2] \label{K2}  Let $ \LL_\lambda $ be the $ \Aa_0 $ module generated by the $ \tilde{F}_i $ acting on $v_\lambda $ and let $ B_\lambda $ be the set of non-zero vectors in $ \LL_\lambda/q^{-1} \LL_\lambda $ obtained by acting on $v_\lambda$ with any sequence of $ \tilde{F}_i $. Then
$(\LL_\lambda, B_\lambda) $ is a crystal basis for $ V_\lambda $.
\end{Theorem}

Theorem \ref{K2} gives a choice of crystal basis for any $V_\lambda$, unique up to an overall scalar. The following result shows that these are all the crystal basis, and furthermore that any crystal basis of a reducible representation $V$ is a direct sum of such bases.

\begin{Theorem}[\cite{Kashiwara:1991}, Theorem 3] \label{uniquecb}
Let $ V $ be a representation of $U_q(\g) $ and let $ (\LL, B) $ be a crystal basis for $ V $.  Then there exists an isomorphism of $ U_q(\g) $ representations $ V \cong \oplus_{j} V_{\lambda_j} $ which takes $ (\LL, B) $ to $ (\oplus_j \LL_{\lambda_j}, \cup_j B_{\lambda_j}) $.
\end{Theorem}

\subsection{Abstract crystals and tensor products}
So far we have been considering crystal bases.  It will also be important for us to consider (abstract) crystals. 

\begin{Definition}
An (abstract) crystal is a finite set $ B $ along with operators $ e_i, f_i : B \rightarrow B \cup \{0\} $ and $ \wt : B \rightarrow P $ which obey certain axioms (see \cite{cactus}).
\end{Definition}

Every crystal basis $ (\LL, B) $ gives an abstract crystal.  Namely, we choose $ B $ to be the underlying set and define $ e_i := \tilde{E}_i^\modq, f_i := \tilde{F}_i^\modq : B \rightarrow B $.  The weight map is defined using the decomposition of the crystal basis into weight spaces. A crystal can be represented by a colored directed graph, where the vertices consist of the elements in $B$, and there is an $i$-colored edge from $b$ to $b'$ if and only if $f_i(b)=b'$. 

\begin{Comment}
From now on we will reserve the word crystal to mean an abstract crystal which arises from a crystal basis as described above.
\end{Comment}

The tensor product rule for $\g$ modules leads to a tensor product rule for crystals, which we will now review. We then present an equivalent definition of the tensor product rule using strings of brackets. This fits more closely with our later definitions, and also helps explain why many realizations of crystals (for example, the realization of $\mbox{sl}_n$ crystals using Young tableaux) make use of brackets.

We start by defining three elements in the root lattice of $\g$ associated to each element $b \in B$.
Let $\g$ be a symmetrizable Kac-Moody algebra, with simple roots indexed by $I$. Let $B$ the crystal of an integrable representation $V$ of $\g$. For each $b \in B$ and $i \in I$, define:
\begin{align}
\varepsilon_i (b) & := \max \{ m : e_i^m (b) \neq 0 \} \\
\varphi_i (b) & := \max \{ m : f_i^m (b) \neq 0 \}.
\end{align}
These are always finite because $V$ is integrable. 
\begin{Definition}
Let $\Lambda_i$ be the fundamental weight associated to $i \in I$.  For each $b \in B$, define three elements in the weight lattice of $\g$ by:
\begin{enumerate}

\item $\displaystyle \varphi(b):= \sum_{i \in I} \varphi_i (b) \Lambda_i$

\item $\displaystyle \varepsilon(b) := \sum_{i \in I} \varepsilon_i (b) \Lambda_i$

\item $\displaystyle \wt (b):=\varphi(b)-\varepsilon(b)$.

\end{enumerate}
\end{Definition}

\begin{Comment} \label{wtcomment}
It turns out that $\wt(b)$ will always be equal to the weight of the corresponding canonical basis element $v(b)$ (see for example \cite{Hong&Kang:2000}).
\end{Comment}

We now give the tensor product rule for crystals, using conventions from \cite{Hong&Kang:2000}.
If $A$ and $B$ are two crystals, the tensor product $A \otimes B$ is the crystal whose underlying set is $\{ a \otimes b: a \in A, b \in B \}$, with operators $e_i$ and $f_i$ defined by:
\begin{equation} \label{edef} e_i (a \otimes b)=
\begin{cases}
e_i  (a) \otimes b, \quad \text{if}\quad  \varphi_i(a) \geq \varepsilon_i(b)\\
a \otimes e_i  (b),\quad \text{otherwise}
\end{cases} \end{equation}
\begin{equation} \label{fdef} f_i (a \otimes b)=
\begin{cases}
f_i  (a) \otimes b, \quad \text{if} \quad  \varphi_i(a) > \varepsilon_i(b)\\
a \otimes f_i (b),\quad \text{otherwise}.
\end{cases} \end{equation}

\noindent This can be reworded as follows. In this form it is known as the signature rule: 

\begin{Lemma} \label{bdef}
For $b \in B$, let $S_i(b)$ be the string of brackets $) \cdots )( \cdots ($, where the number of  $``)"$ is $\varepsilon_i(b)$ and the number of $``("$ is $\varphi_i(b)$. Then the actions of $e_i$ and $f_i$ on $A \otimes B$ can be calculated as follows: 
\begin{equation}
e_i (a \otimes b) \hspace{-2pt} =  \hspace{-2pt}
\begin{cases}
a \otimes e_i(b) \neq 0 \quad \mbox{if the first uncanceled } ``)"  \mbox{ from the right in } S_i(a) S_i(b) \mbox{ is in } S_i(b) \\
e_i(a) \otimes b \neq 0 \quad \mbox{if the first uncanceled } ``)"  \mbox{ from the right in } S_i(a) S_i(b) \mbox{ is in } S_i(b) \\
0 \hspace{0.92in} \mbox{if there is no uncanceled } ``)" \mbox{ in } S_i(a)S_i(b)
\end{cases}
\end{equation}
\begin{equation}
f_i (a \otimes b)  \hspace{-2pt} =  \hspace{-2pt}
\begin{cases}
f_i(a) \otimes b \neq 0 \quad \mbox{if the first uncanceled } ``("  \mbox{ from the left in } S_i(a) S_i(b) \mbox{ is in } S_i(a) \\
a \otimes f_i (b) \neq 0 \quad \mbox{if the first uncanceled } ``("  \mbox{ from the left in } S_i(a) S_i(b) \mbox{ is in } S_i(b) \\
0 \hspace{0.92in} \mbox{if there is no uncanceled } ``(" \mbox{ in } S_i(a)S_i(b)
\end{cases}
\end{equation}
\end{Lemma}

\begin{proof} This formula for calculating $e_i$ and $f_i$ follows immediately from Equations (\ref{edef}) and (\ref{fdef}). To see that $e_i(a \otimes b) \neq 0$ if there are any uncanceled $``)"$, notice that in this case you always act on a factor that contributes at least one $``)"$, and hence has $\epsilon_i > 0$. By the definition of $\epsilon_i$, $e_i$ does not send this element to $0$. The proof for $f_i$ is similar. \end{proof}
The advantage of Lemma \ref{bdef} over equations (\ref{edef}) and (\ref{fdef}) is that we can easily understand the actions of $e_i$ and $f_i$ on the tensor product of several crystals:

\begin{Corollary} \label{mycrystaldef}
 Let $B_1, \ldots, B_k$ be crystals of integrable representations of $\g$. Let $b_1 \otimes \cdots \otimes b_k \in B_1 \otimes \cdots \otimes B_k$. For each $1 \leq j \leq k$, let $S_i(b_j)$ be the string of brackets $) \cdots )( \cdots ($, where the number of  $``)"$ is $\varepsilon_i(b_j)$ and the number of $``("$ is $\varphi_i(b_j)$. Then:
 \begin{enumerate} 
\item \label{cc1}  
$\displaystyle e_i (b_1 \otimes \cdots \otimes  b_k) =
\begin{cases}
b_1 \otimes  \cdots e_i(b_j)  \cdots \otimes b_k \neq 0 \quad \begin{array}{l}  \mbox{if the first uncanceled } ``)"  \mbox{ from the } \\ \mbox{right in } S_i(b_1) \cdots S_i(b_k) \mbox{ is in } S_i(b_j) \end{array} \\
0 \hspace{1.8in} \mbox{if there is no uncanceled } ``)".
\end{cases}$

\item \label{cc2}
$\displaystyle f_i (b_1 \otimes \cdots \otimes  b_k)= 
\begin{cases}
b_1 \otimes  \cdots f_i(b_j)  \cdots \otimes b_k \neq 0 \quad \begin{array}{l}  \mbox{if the first uncanceled } ``("  \mbox{ from the} \\ \mbox{left in } S_i(b_1) \cdots S_i(b_k)  \mbox{ is in } S_i(b_j) \end{array} \\
0 \hspace{1.8in} \mbox{if there is no uncanceled } ``(".
\end{cases}$

\item \label{cc3}
$\varepsilon_i(b_1 \otimes \cdots \otimes b_k)$ is the number of uncanceled $``)"$ in $S_i(b_1) \cdots S_i(b_k)$.

\item  \label{cc4}
$\varphi_i (b_1 \otimes \cdots \otimes b_k)$ is the number of uncanceled $``("$ in $S_i(b_1) \cdots S_i(b_k)$.

\end{enumerate}
 \end{Corollary}
 
 \begin{proof}
 Parts (\ref{cc1}) and (\ref{cc2}) follow by iterating Lemma \ref{bdef}. To see part (\ref{cc3}), notice that, if the first uncanceled $``)"$ is in $S_i (b_j)$, then $\varepsilon_i (b_j) \geq 1$ and $e_i$ acts on $b_j$. Hence $e_i$ changes $S_i(b_j)$ by reducing the number of $``)"$ by one, and increasing the number of $``("$ by one. The only affect on $S_i(b_1) \cdots S_i (b_k)$ is that the first uncanceled $``)"$ is changed to $``("$. This reduces the number of uncanceled $``)"$ by one. $e_i$ will send the element to $0$ exactly when there are no uncanceled $``)"$ left. Hence part (\ref{cc3}) follows by the definition of $\varepsilon_i$. Part (\ref{cc4}) is similar. \end{proof}
 
 \section{Commutors}
 
 We review the crystal commutor, as defined by Henriques and Kamnitzer. We then review a related coboundary structure on $U_q(\g)$ representations, and discuss how the crystal commutor arises as a combinatorial limit of this construction.
 
 \subsection{The crystal commutor} \label{commutor_section}

In \cite{cactus}, Henriques and Kamnitzer established the existence (and uniqueness) of a Sch\"utzenberger involution $ \xi_B : B \rightarrow B $ on crystals associated to representations of $U_q(\g)$. This satisfies the properties
\begin{equation}
\xi_B(e_i \cdot b) = f_{\theta(i)} \cdot \xi_B(b), \quad \xi_B(f_i \cdot b) = e_{\theta(i)} \cdot \xi_B(b), \quad wt(\xi_B(b)) = w_0 \cdot wt(b).
\end{equation}

Following a suggestion of A. Berenstein, the Sch\"utzenberger involution was used in \cite{cactus} to define the commutor for crystals by the formula
\begin{equation}
\begin{aligned} \label{commuter_definition}
  \sigma^{HK}_{A,B} : A\otimes B &\rightarrow B \otimes A \\
 a \otimes b &\mapsto \xi_{A \otimes B} (\xi_B(b) \otimes \xi_A(a)) = \mathrm{Flip} \circ (\xi_A \otimes \xi_B) (\xi_{A \otimes B} (a \otimes b)).
 \end{aligned}
 \end{equation}
The second expression here is just the inverse of the first expression, and the equality is proved in \cite{cactus}.

\begin{Theorem}[\cite{cactus}, Theorem 6] $\g$-Crystals, with the above tensor product rule and commutor, forms a coboundary category.
\end{Theorem}

\subsection{Lifting the crystal commutor to representations} \label{xionreps}

In the previous section we described Sch\"utzenberger involution $\xi$ as an involution on the crystal associated to a representation $V$ of $U_q(\g)$, and how $\xi$ is used to define the crystal commutor $\sigma^{HK}$. We now describe (following \cite{cactus}) how to modify this construction to obtain an involution of the actual representation $V$, and hence a commutor for $U_q(\g)$ representations. 

There is a one dimensional family of maps $ V_\lambda \rightarrow V_\lambda $ which exchange the actions of $ E_i, F_i $ with $ F_{\theta(i)}, E_{\theta(i)} $. Define $ \xi_{V_\lambda} $ to be the unique such map which takes the highest weight basis vector $v_\lambda $ to the lowest weight vector $v_\lambda^{\text{low}}$ (see Definition \ref{lowdef}).

Hence we have constructed an element $ \xi \in \widetilde{U_q(\g)} $ such that conjugation by $ \xi $ is given by
\begin{equation}
\begin{cases}
C_{\xi}(E_i)  = F_{\theta(i)}\\
C_{\xi}(F_i) = E_{\theta(i)}\\
C_{\xi}(K_H) = K_{w_0 \cdot H}
\end{cases}
\end{equation}

We can now define a commutor for the category of $U_q(\g)$ representations. For any representations $V $ and $W$ of $U_q(\g)$, define a map from $V \otimes W$ to $W \otimes V$ by
\begin{equation}
\sigma^{\xi}_{V,W} := \xi_{W \otimes V} \circ (\xi_W \otimes \xi_V) \circ \mbox{Flip}  = \mathrm{Flip} \circ \xi_V \otimes \xi_W \circ \xi_{V \otimes W}.
\end{equation}
Note that $ \xi $ is a coalgebra anti-automorphism, so, as in Section \ref{autint_intro}, $ \sigma^{\xi} $ is an isomorphism of $ U_q(\g) $ modules. Henriques and Kamnitzer \cite{cactus} show $ \sigma^{\xi} $ gives $Rep(U_q(\g)$ the structure of a coboundary monoidal category (see Definition \ref{coboundary_cat_def}). They also note that this coboundary structure is not unique. In Section \ref{Rcommutor_section}, we will show that Drinfeld's coboundary structure $\sigma^D$ can be obtained as a slight modification of $\sigma^{\xi}$, thus relating the two ways of putting coboundary structures on the category of $U_q(\g)$ representations.

\subsection{The crystal commutor as a combinatorial limit} 
We now show that Sh\"utzenberger involution on representations, as defined in Section \ref{xionreps}, induces Sch\"utzenberger involution on crystal bases, as defined in Section \ref{commutor_section}. This implies that the coboundary structure $\sigma^{\xi}$ on $U_q(\g)$ representation induces the crystal commutor $\sigma^{HK}$ on a tensor product of crystal bases. 
We begin with the following lemmas.

\begin{Lemma} \label{lowlemma}
Fix an irreducible representation $V_\lambda$ with a chosen highest weight vector $v_\lambda$. Let $v_\lambda^\text{low}$ be as in Definition \ref{lowdef}. Then, for any reduced word $ w_0 = s_{i_1} \cdots s_{i_m} $, we have
\begin{equation} \label{adf1}
v_{\lambda}^{\text{low}} = \tilde{F}_{i_m}^{n_m} \cdots \tilde{F}_{i_1}^{n_1} v_\lambda \end{equation}
\begin{equation} \label{adf2}
v_{\lambda}= \tilde{E}_{\theta(i_m)}^{n_m} \cdots \tilde{E}_{\theta(i_1)}^{n_1} v_\lambda^\text{low} \end{equation}
where, as in Proposition \ref{jjthings}, $n_j = \langle  s_{i_1} \cdots s_{i_{j-1}}\alpha_{i_j}^\vee, \lambda \rangle $.
\end{Lemma}
\begin{proof}
By Proposition \ref{jjthings}, for each $0  \leq k <m$, we have $E_{i_{k+1}}F_{i_k}^{(n_k)} \cdots F_{i_1}^{(n_1)} v_\lambda=0$.  So Equation (\ref{adf1}) follows from the definition of the Kashiwara operators (Definition \ref{Kop}).

Now $w_0= s_{i_1} \cdots s_{i_m}$ is a reduced word in $W$, which implies that $s_{\theta(i_m)} \cdots s_{\theta(i_1)}$ is as well. Thus by Equation (\ref{adf1}),
\begin{equation} \label{anewe}
v_{\lambda}^{\text{low}} = \tilde{F}_{\theta(i_1)}^{\ell_m} \cdots \tilde{F}_{\theta(i_m)}^{\ell_1} v_\lambda,\end{equation}
where $\ell_j = \langle  s_{\theta(i_m)} \cdots s_{\theta(i_{m-j+2})}\alpha_{\theta(i_{m-j+1})}^\vee, \lambda \rangle $. For all $j$, 
\begin{align}
\label{si1} s_{i_1} \cdots s_{i_{j-1}}\alpha_{i_j}^\vee  &= w_0 w_0  s_{i_1} \cdots s_{i_{j-1}}\alpha_{i_j}^\vee \\
 &\label{si2}= -w_0 s_{\theta(i_1)} \cdots s_{\theta(i_{j-1})} \alpha_{\theta(i_j)}^\vee \\
 &\label{si3}=-s_{\theta(i_m)} \cdots s_{\theta(i_1)}  s_{\theta(i_1)} \cdots s_{\theta(i_{j-1})} \alpha_{\theta(i_j)}^\vee \\
&\label{si4}=-s_{\theta(i_m)} \cdots s_{\theta(i_j)} \alpha_{\theta(i_j)}^\vee \\
&\label{si5}=s_{\theta(i_m)} \cdots s_{\theta(i_{j+1})} \alpha_{\theta(i_j)}^\vee.
\end{align}
Here (\ref{si1}) follows because $w_0$ is an involution, and (\ref{si2}) follows because for all $i \in I$, $w_0 (\alpha_i)= - \alpha_{\theta(i)}$ and $w_0 s_i w_0^{-1} = s_{\theta(i)}$. Thus $n_j = \ell_{m-j+1}$ so, by Equation (\ref{anewe}), 
\begin{equation}
v_{\lambda}^{\text{low}} = \tilde{F}_{\theta(i_1)}^{n_1} \cdots \tilde{F}_{\theta(i_m)}^{n_m} v_\lambda.\end{equation}
By Definition \ref{Kop}, this is equivalent to Equation (\ref{adf2}).
\end{proof}

The following follows from \cite[Theorem 3.3 (b)]{L2}, although we provide a proof.

\begin{Proposition}\label{bottotop} The action of the element $ \xi_{V_\lambda} $ defined in Section \ref{xionreps} is an involution.  In particular
\begin{equation} \xi_{V_\lambda}(v_\lambda^{\low}) = v_\lambda. \end{equation}
\end{Proposition}

\begin{proof}
Recall that $\xi_{V_\lambda}$ interchanges the action of $E_i$ and $F_i$. By Comment \ref{symt}, $\xi_{V_\lambda}$ also interchanges the action of $\tilde{E}_i$ and $\tilde{F}_i$. Thus, applying $\xi_{V_\lambda}$ to both sides of Equation (\ref{adf1}), 
\begin{equation}
\xi(v_\lambda^\text{low})= \tilde{E}_{\theta(i_m)}^{n_m} \cdots \tilde{E}_{\theta(i_1)}^{n_1} v_\lambda^\text{low}.
\end{equation}
The result then follows by Equation (\ref{adf2}).
\end{proof}

We can now describe how $ \xi_V $ acts on a crystal basis:  
\begin{Theorem} \label{xionbasis}
Let $ (\LL, B) $ be a crystal basis for a representation $ V $. Then the following holds.
\begin{enumerate}
\item  $ \xi_{V}(\LL) = \LL$. 
\item By (i), $ \xi_V $ gives rise to a map between $ \xi_V^\modq : \LL/q^{-1} \LL \rightarrow \LL/q^{-1}\LL $.  For each $ b \in B $, we have 
\begin{equation} \xi^\modq_V(b) =  \xi_B(b). \end{equation}
\end{enumerate}
\end{Theorem}

\begin{proof}
First, we note that it is sufficient to prove the theorem in the case that $ (\LL,B) =(\LL_\lambda, B_\lambda) $.  The general case of the theorem then follows from an application of Theorem \ref{uniquecb}.

So  assume that $ V = V_\lambda, \LL = \LL_\lambda, B = B_\lambda $.  Note that $ \xi_{V_\lambda} $ exchanges the action of $ E_i $ and $ F_{\theta(i)} $ and hence, by Comment \ref{symt}, interchanges the actions of $\tilde{E}_i $ and $ \tilde{F}_{\theta(i)} $.  Since $ \LL $ is generated by $ \tilde{F}_i $ acting on $ v_\lambda$, we see that $ \xi_{V_\lambda}(\LL) $ is generated by $ \tilde{E}_i $ acting on $ \xi_{V_\lambda}(v_\lambda) = v_\lambda^{\text{low}} $.  Lemma \ref{lowlemma} shows that $ v_\lambda^{\text{low}} \in \LL $ so, since $ \LL $ is invariant under the action of the $ \tilde{E}_i $, we conclude that $ \xi_{V_\lambda}(\LL) \subset \LL $. By Proposition \ref{bottotop}, $\xi_{V_\lambda}^2$ is the identity, so in fact we must have $\xi_{V_\lambda}(\LL)=\LL$. 

For part (ii), note that $ v_\lambda^{\text{low}} $ is obtained by acting on $ v_\lambda $ with the $ \tilde{F}_i$.  Hence its reduction $ \text{mod} q^{-1}\LL $ must lie in $ B$, and in fact must by the lowest weight element in $B$.  The result follows because $C_\xi$ acts on  the set of Kashiwara operators according to $F_i \leftrightarrow E_{\theta(i)}$.
\end{proof}

The following establishes the connection between $\sigma^{\xi}$ and $\sigma^{HK}$:

\begin{Corollary}
Let $ (\LL_1, B_1) $ and $(\LL_1, B_2)$ be crystal bases for two representations $ V $ and $W$. Then the following holds.
\begin{enumerate}
\item  $ \sigma^{\xi}_{V,W}(\LL_1 \otimes \LL_2) = \LL_2 \otimes \LL_1$. 
\item By (i), $ \sigma^\xi_{V,W} $ gives rise to a map between $ \sigma_{V,W}^{\xi \quad \modq} : \LL/q^{-1} \LL \rightarrow \LL/q^{-1}\LL $.  For each $ b_1 \otimes b_2 \in B_1 \otimes B_2 $, we have 
\begin{equation} \sigma_{V,W}^{\xi \quad \modq} (b_1 \otimes b_2) =  \xi^{HK}_{B_1, B_2} (b_1 \otimes b_2).\end{equation}
\end{enumerate}
\end{Corollary}

\begin{proof}
Follows from Theorem \ref{xionbasis} and definitions.
\end{proof}

\begin{Comment}
There is an even stronger connection between $ \xi_{V_\lambda} $ and the canonical (or global) basis $B_\lambda^c $ for $ V_\lambda $: It follows from \cite[Chapter 21]{Lusztig:1993} that $ \xi_{V_\lambda} $ is the linear extension of the set map $ \xi_\lambda : B_\lambda^c \rightarrow B_\lambda^c $. However, because this fact does not hold for tensor products of canonical bases, it will not be useful for us. That is why we state the weaker fact above which holds for all crystal bases.
\end{Comment}

\section{The abacus}  \label{lotsastuff}

Here we explain the abacus used by James and Kerber in \cite{JK:1981}. We start by defining a bijection between partitions and rows of beads. This is essentially the correspondence between partitions and semi-infinite wedge products (see for example \cite[Chapter 14]{Kac:1990}); the positions of the beads correspond to the factors in the wedge product. As in \cite{JK:1981}, this one row is transformed into several parallel rows of beads on an abacus.

We use the "Russian" diagram of a partition, shown in Figure \ref{partition}. Place a bead on the horizontal axis under each down-sloping segment on the edge of the diagram. The corresponding row of beads uniquely defines the partition. Label the horizontal axis so the corners of all boxes are integers, with the vertex at $0$. For a partition $\lambda= (\lambda_1, \lambda_2, \ldots)$, where we define $\lambda_i=0$ for large $i$, the positions of the beads will be $\lambda_i - i + 1/2$. The empty partition corresponds to having beads at all negative positions in $\bz + 1/2$, and none of the positive positions. Adding a box to the partition corresponds to moving a bead one step to the right.

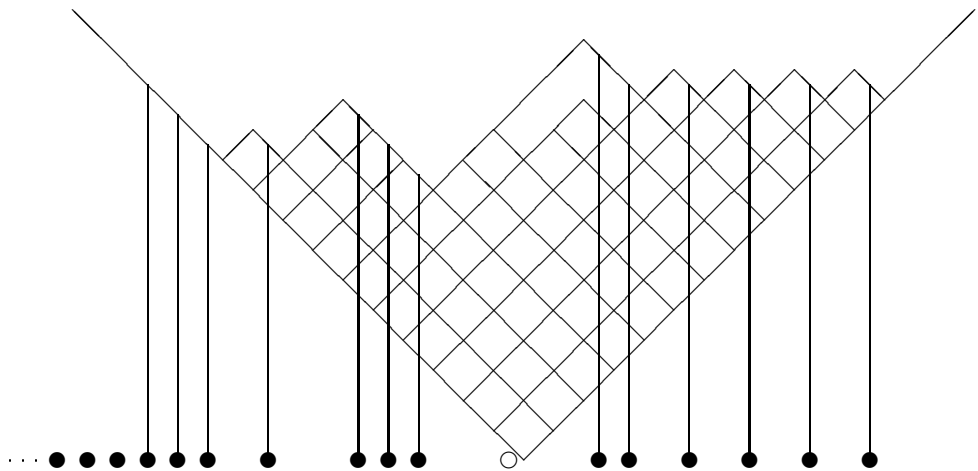
\begin{figure}
\begin{center}
\input{addribbon}
\end{center}
\caption[Adding a ribbon to a partition]{Adding a 4-ribbon to a partition corresponds to moving one bead forward 4 positions, possibly jumping over some beads in between. In this example, the bead that was at position $-0.5$ is moved forward four places to position $3.5$. The diagram has changed above positions $0.5, 1.5$ and $2.5$, but adding the ribbon won't change the slope of the diagram there, so the beads stay in the same place.
\label{addribbon}}
\end{figure}

Recall that a {\it ribbon} (or rim-hook) is a skew partition $\lambda \backslash  \mu$ whose diagram is connected and has at most one box above each position on the horizontal axis (i.e. it is a strip one box wide). Adding an $\ell$-ribbon (that is, a ribbon with $\ell$ boxes) to a partition $\lambda$ moves one bead exactly $\ell$ steps to the right in the corresponding row of beads, possibly jumping over other beads. See Figure \ref{addribbon}. 
\setlength{\unitlength}{0.4cm}

\begin{figure}
\begin{center}
\begin{picture}(30,15)
\input{partition}
\end{picture}
\end{center}
\caption[Partitions and bead strings]{Partitions and bead strings. This figure shows the bead string corresponding to the partition  (12,11,10,9,7,5,3,3,3,1). We fix the origin to be the vertex of the partition, so the beads are in positions
$  \ldots, -12.5, -11.5, -10.5, -8.5, -5.5, -4.5, -3.5, -0.5, 2.5, 5.5, 7.5, 9.5, 11.5.$  \label{partition}}
\end{figure}

\setlength{\unitlength}{0.5cm}

It is often convenient to work with just the "bead" picture. In order to avoid confusion, we denote the empty spaces by white beads, and indicate the position of the origin by a line. The example shown in Figure \ref{partition} becomes:

\vspace{-0.22in}
\begin{center}
\begin{picture}(30,2)
\label{bead_strand}
\input{bead_strand}
\end{picture}
\end{center}
\vspace{0.2in}

\vspace{-0.08in}

\noindent We then put the beads into groups of $\ell$ starting at the origin, as shown below for $\ell=4$:

\vspace{-0.18in}
\begin{center}
\begin{picture}(30,2)
\label{grouped_beads}
\input{grouped_beads}
\end{picture}
\end{center}

\vspace{0.13in}

\noindent Rotating each group 90 degrees counterclockwise, and compressing, we get $\ell$ rows:

\begin{center}
\input{first_abacus}
\end{center}

\noindent We will call this the level $\ell$ abacus (in this example $\ell=4$).
Adding an $\ell$-ribbon now corresponds to moving one bead forward one position, staying on the same row. For instance, adding the ribbon as in Figure \ref{addribbon} corresponds to moving the third black bead from the right on the top row, which gives:

\begin{center}
\begin{picture}(30,4)
\input{second_abacus}
\end{picture}
\vspace{0.1in}
\end{center}

As explained in (\cite{JK:1981} Chapter 2.7), this model immediately gives some interesting information about the partition: The $\ell$ rows can be interpreted as $\ell$ partitions, using the correspondence between partitions and rows of beads (shifting the origin if necessary). This is known as the $\ell$-quotient of $\lambda$.  We can also consider the partition obtained by pushing the beads on each row as far to the left as they will go, but not changing rows (and only doing finitely many moves). This is the $\ell$-core.

\chapter{The crystal commutor and related constructions} \label{commutor_chapter}

This Chapter contains several results related to the crystal commutor. Sections \ref{Rcommutor_section} and \ref{commutorKashiwara_section} are based on joint work with Joel Kamnitzer (\cite{Rcommutor} and \cite{commutorKashiwara}). Section \ref{RTheta_section} is partly based on \cite{RTheta}, although the construction given here is somewhat different and works in greater generality.

\section{The crystal commutor and Drinfeld's unitarized $R$-matrix} \label{Rcommutor_section}

Let $A$ and $B$ be the crystals of two representations of a simple complex Lie algebra $\g$. Using the Sch\"utzenberger involution, Henriques and Kamnitzer \cite{cactus} defined a natural isomorphism $A \otimes B \rightarrow B \otimes A$, which they call the crystal commutor. This gives $\g$-crystals the structure of a coboundary category.

By an analogous construction, Henriques and Kamnitzer also defined a commutor $\sigma^\xi_{V, W} : V \otimes W \rightarrow W \otimes V$, where $V$ and $W$ are finite dimensional representations of $U_q(\g)$. There is some choice in lifting the Sch\"utzenberger involution to representations, so the commutor here is not unique. 

There is a more standard isomorphism from $V \otimes W$ to $W \otimes V$, called the braiding. This is defined by $v \otimes w \mapsto \Flip \circ  R ( v \otimes w )$, where $R$ is the universal $R$ matrix. In \cite[Section 3]{Drinfeld:1990}, Drinfeld introduced a "unitarized" $R$ matrix $\bar{R}$. He showed that the map $V \otimes W \rightarrow W \otimes V$ given by $v \otimes w \mapsto \Flip \circ  \bar{R} (v \otimes w)$ is a coboundary structure on the category of $U_q(\g)$ representations.

The first purpose of this section is to relate these two ways of putting a coboundary structure on the category of $U_q(\g)$ representations, thus answering a question from \cite{cactus}. We then show that, for any two crystal bases, Drinfeld's commutor preserves the associated lattices and acts by the crystal commutor on the bases (up to some negative signs).  Thus the crystal commutor is essentially a combinatorial limit of Drinfeld's commutor for representations. We feel this explains why the crystal commutor is a coboundary structure, and not a braiding, as one might naively expect.

\subsection{Realizing $\bar{R}$ in the form $(Y^{-1} \otimes Y^{-1} ) \Delta{Y}$.} \label{RY}

We will construct the unitarized $R$ matrix in the desired form by modifying a similar result for the standard $R$ matrix, due to Kirillov-Reshetikhin and Levendorskii-Soibelman. Their result is stated as an expression in the $h$-adic completion of $U_h(\g) \otimes U_h(\g)$, although it in fact does give a well defined action on $V \otimes W$ for any representations $V$ and $W$ of $U_q(\g)$, so is well defined in $\widetilde{U_q(\g) \otimes U_q(\g)}$ (see Definition \ref{bigbig_completion_def}). In order to use Theorem \ref{sR}, for this section only, we will write some expressions in the $h$-adic completions of $U_h(\g)$ and  $U_h(\g) \otimes U_h(\g)$, and simply note that all the ones we use are well defined in $\widetilde{U_q(\g)}$ and $\widetilde{U_q(\g) \otimes U_q(\g)}$ as well. Translating the conventions in \cite{KR:1990} and \cite{LS:1991} into ours, we obtain the following result.

\begin{Theorem}[{\cite[Theorem 3]{KR:1990}, \cite[Theorem 1]{LS:1991}}] \label{sR}
With notation as in Section \ref{notation}, the standard $R$-matrix for $U_h(\g)$ can be realized as
\begin{equation}
R=  \exp \left( h \sum_{i, j \in I} (B^{-1})_{ij} H_i \otimes H_j \right) (T_{w_0}^{-1} \otimes T_{w_0}^{-1}) \Delta(T_{w_0}).
\end{equation}
\end{Theorem}

\begin{Definition} \label{defJ}
Let $J$ be the operator which acts on a finite dimensional representation $V$ of $U_q(\g)$ by multiplying each vector of weight $ \mu $ by $q^{( \mu, \mu) /2 + ( \mu, \rho )}$. It is a straightforward calculation to see that $J$ can be realized in a completion of  $U_h(\g)$ by
\begin{equation} \label{Jdefeq}
J:= \exp  \left[ h \left( \frac{1}{2} \sum_{i,j} \left( (B^{-1})_{ij} H_i H_j \right) + H_\rho \right) \right].
\end{equation}
\end{Definition}

Actually, $ (\mu, \mu)/2 + (\mu, \rho)$ can in some cases be a fraction. Once $\g$ is fixed, the possible denominators are bounded. To be precise we should adjoin a fixed $k^{th}$ root of $q$ to our base field, for some $k$ depending on $\g$. This causes no difficulty. 

\begin{Comment}
It follows from Lemma \ref{Dj} below that Theorem \ref{sR} is equivalent to saying $R= (X^{-1}\otimes X^{-1} ) \Delta(X)$, where $X= J T_{w_0}$.
\end{Comment}

\begin{Definition} \label{defY} $Y$ is the element in the completion of $ U_q(\g)$ defined by $Y:= Q^{-1/2} J T_{w_0} $.
\end{Definition} 
 
 We are now ready to state the main result of this section. 
\begin{Theorem} \label{answer1}
The unitarized $R$ matrix can be realized as
\begin{equation}
\bar{R}= (Y^{-1} \otimes Y^{-1} ) \Delta(Y).
\end{equation}
\end{Theorem}

\begin{Comment}
In fact, $Y$ is a well defined operator on $U_q(\g)$ over $\bc(q)$. That is, unlike for the standard $R$ matrix, we do not actually need to adjoin a $k^{th}$ root of $q$.
\end{Comment}

We prove Theorem \ref{answer1} by a direct calculation, using Theorem \ref{sR}. 
We will need the following technical lemma:

\begin{Lemma} \label{Dj}
\begin{equation}
\Delta(J) = ( J \otimes J) \exp \Big( h\sum_{i,j \in I} (B^{-1})_{ij} H_i \otimes H_j \Big) 
\end{equation}
\end{Lemma}

\begin{proof}
\begin{align}
\Delta(J) &= \Delta \bigg(\exp  \bigg[ h \Big( \frac{1}{2} \sum_{i,j}  (B^{-1})_{ij} H_i H_j  + H_\rho \Big) \bigg] \bigg) \\
&= \exp \bigg[  h \bigg( \frac{1}{2}  \sum_{i,j} (B^{-1})_{ij} (H_i \otimes 1 + 1 \otimes H_i)(H_j \otimes 1 + 1 \otimes H_j)  + H_\rho \otimes 1 + 1 \otimes H_\rho \bigg) \bigg]\\
& = \nonumber \exp \Bigg[ h \bigg( \frac{1}{2}  \sum_{i,j} (B^{-1})_{ij} H_i H_j \otimes 1 + H_\rho \otimes 1 \bigg) \Bigg] \times \\
& 
\hspace{0.3in}  \times \exp \bigg[ h \bigg(  \frac{1}{2}  \sum_{i,j}  (B^{-1})_{ij} 1 \otimes H_i H_j  + 1 \otimes H_\rho \bigg) \bigg]  \exp \bigg[ h \sum_{i,j \in I} ( B^{-1})_{ij} H_i \otimes H_j  \bigg] \\
&=  (J \otimes J) \exp \Big[ h\sum_{i,j \in I} ( B^{-1})_{ij} H_i \otimes H_j \Big].
\end{align}
\end{proof}

%%%%%%%%%%%%We no longer seem to use this Lemma%%%%%%%%%%%
%\begin{Lemma} \label{jt}
%$J T_{w_0} = T_{w_0} J K_\rho^{-2} = K_\rho^2 T_{w_0} J$
%\end{Lemma}

%\begin{proof}
%Choose some weight vector $v$ of weight $\mu$ in some irreducible representation $V_\lambda$ of $U_q(\g)$. It is sufficient to show that
%\begin{equation}
%J T_{w_0} \cdot v= T_{w_0} J K_\rho^{-2} \cdot v= K_\rho^2 T_{w_0} J \cdot v.
%\end{equation}
%This is a simple calculation using Definition \ref{defJ}, and recalling that $T_{w_0}$ acts as $w_0$ on weights.
%\end{proof}

\begin{proof}[Proof of Theorem \ref{answer1}]
From the definition, we have that
\begin{align}
(Y^{-1} \otimes Y^{-1} ) \Delta(Y)  \hspace{-0.02in}&= \hspace{-0.02in}
( T_{w_0}^{-1} J^{-1} Q^{1/2} \otimes  T_{w_0}^{-1} J^{-1} Q^{1/2} ) \Delta( Q^{-1/2}J T_{w_0}) \\
 &= \hspace{-0.02in} ( T_{w_0}^{-1} \otimes T_{w_0}^{-1}) (J^{-1} \otimes J^{-1}) (Q^{1/2} \otimes Q^{1/2})  \Delta(Q^{-1/2})\Delta(J) \Delta(T_{w_0}) \\
 &= \hspace{-0.02in} (Q^{1/2} \otimes Q^{1/2})  ( T_{w_0}^{-1} \otimes T_{w_0}^{-1}) (J^{-1} \otimes J^{-1}) \Delta(J) \Delta(T_{w_0}) \Delta(Q^{-1/2}),
 \end{align}
where the last equality follows because $Q^{1/2}$ is central.
Then by Lemma \ref{Dj}:
\begin{align}
 \nonumber & (Y^{-1} \otimes Y^{-1} ) \Delta(Y)\\
 &= (Q^{1/2} \otimes Q^{1/2} )( T_{w_0}^{-1} \otimes T_{w_0}^{-1} ) 
\exp \Big( h \sum_{i, j \in I} (B^{-1})_{ij} H_i \otimes H_j \Big)
\Delta(T_{w_0}) \Delta(Q^{-1/2}) \\
\label{t13}&= (Q^{1/2} \otimes Q^{1/2}  ) 
\exp \Big( h \sum_{i, j \in I} (B^{-1})_{ij} H_i \otimes H_j \Big)
( T_{w_0}^{-1} \otimes T_{w_0}^{-1}) \Delta(T_{w_0}) \Delta(Q^{-1/2}) \\
\label{t14} &= (Q^{1/2} \otimes Q^{1/2}) R \Delta(Q^{-1/2}).
\end{align}
Here Equation (\ref{t13}) follows because $T_{w_0}$ permutes weight spaces as $w_0$, so $T_{w_0} H_i T_{w_0}^{-1} = -H_{\theta(i)}$. Equation (\ref{t14}) follows by Theorem \ref{sR}. The theorem follows by Equation (\ref{QR}).
\end{proof}

\subsection{Realizing $\bar{R}$ using Sch\"utzenberger involution} \label{RS}

This section contains our first main result of Chapter \ref{commutor_chapter} (Corollary \ref{answer2}), which realizes the unitarized $R$-matrix using a slight modification of Sch\"utzenberger involution. 
Fix two representations $V$ and $W$ of $U_q(\g)$. As discussed in Section \ref{xionreps}, there is a natural isomorphism $\sigma^\xi_{V,W} : V \otimes W \rightarrow W \otimes V$ defined by
\begin{equation} \label{xx1}
\sigma^{\xi}_{V,W} = \Flip \circ (\xi_V^{-1} \otimes \xi_W^{-1}) \circ \xi_{V \otimes W},
\end{equation}
where $\xi$ is Sch\"utzenberger involution. We have added inverses to the expression to make it more like Theorem \ref{answer1}. At the moment this has no effect, since $\xi$ is an involution, but it will be important later on. 

The commutor $\sigma^{\xi}$ endows the Category of $U_q(\g)$ representations with a coboundary structure, as defined by Drinfeld \cite[Section 3]{Drinfeld:1990}. In \cite{cactus}, Henriques and Kamnitzer note that one can multiply the action of $\xi$ on each irreducible representation by $\pm 1$, with the signs chosen independently for each $V_\lambda$, and Equation (\ref{xx1}) still defines a coboundary structure. They ask if there is a choice of signs such that the resulting commutor coincides with $\Flip \circ \bar{R}$, where $\bar{R}$ is the unitarized $R$ matrix. It turns out that we need  little bit more freedom. At the end of this section, we realize $\Flip \circ \bar{R}$ in terms of Sch\"utzenberger involution, where the action of $\xi$ on each irreducible representation is rescaled by certain $4^{th}$ roots of unity. It is convenient to first work with a different modification of $\xi$.

\begin{Definition} \label{xi11}
$ \xi'' $ is the element of $\widetilde{U_q(\g)}$ which acts on a weight vector $v \in V_\lambda$ by $ \xi''(v) = (-1)^{\langle\mu - w_0 (\lambda), \rho^\vee\rangle } \xi(v)$, where $ \mu $ is the weight of $ v$. Notice, that $\mu - w_0 ( \lambda)$ is always in the root lattice, so $\langle\mu - w_0 ( \lambda), \rho^\vee\rangle$ is always an integer.
\end{Definition}

\begin{Proposition} \label{adanddiag}
$T_{w_0}, \xi''$ and $J$ are all invertible in $\widetilde{U_q(\g)}$, the actions of $C_{T_{w_0}}$, $C_{\xi''}$ and $C_{Q^{-1/2} J} = C_{J}$ all preserve the subalgebra $U_q(\g)$, and:
\begin{enumerate}

\item \label{adt} $
\begin{cases}
C_{T_{w_0}}(E_i) = -F_{\theta(i)} K_{\theta(i)} \\
C_{T_{w_0}}(F_i) = -K_{\theta(i)}^{-1} E_{\theta(i)} \\
C_{T_{w_0}}(K_H) = K_{w_0(H)}, \text{ so that } C_{T_{w_0}}(K_i) = K_{\theta(i)}^{-1}
\end{cases}
$

\item \label{adj} $
\begin{cases}
C_{J}(E_i) = K_i E_i  \\
C_{J}(F_i) = F_i K_i^{-1}\\
C_{J}(K_H) =  K_H
\end{cases}$

\item \label{adx}
$
\begin{cases}
C_{\xi''}(E_i)  = -F_{\theta(i)}\\
C_{\xi''}(F_i) = -E_{\theta(i)}\\
C_{\xi''}(K_H) = K_{w_0 \cdot H}
\end{cases}
$
\end{enumerate}
Furthermore, $ Y = \xi'' $ where, as in Section \ref{RY}, $Y= Q^{-1/2} J T_{w_0}$.
\end{Proposition}

\begin{proof}
It is clear from the definitions that these elements act as invertible endomorphisms on each $V_\lambda$, and hence by Theorem \ref{comp} they are invertible. Their conjugation actions preserve $U_q(\g)$ by (i), (ii) and (iii), which we prove below.

(\ref{adt}) This is Lemma \ref{stillE}.

(\ref{adj}) Let $v$ be a vector of weight $\mu$ in some finite dimensional representation. It is a straightforward calculation to see that $J(E_i v) = K_i E_i (J(v))$, $J(F_i(v))= F_i K_i^{-1} (J(v))$ and $J(K_H (v)) = K_{w_0 \cdot H} J(v)$.

(\ref{adx}) For $C_{\xi''} (E_i)$ and $C_{\xi''}( F_i)$ this follows immediately from Definition \ref{xi11} and the definition of Sch\"utzenberger involution (see Section \ref{xionreps}). It is a straightforward calculation to show that for any weight vector $v$, $\xi''(K_H (v)) = K_{w_0 (H)} \xi'' (v)$. It follows that $C_{\xi''} (K_H) = K_{w_0(H)}$.

It remains to show that $Y=\xi''$.
A direct calculation using (i), (ii) and (iii) shows that $C_{J} C_{T_{w_0}}= C_{\xi''}$. Since $Q^{-1/2}$ is central, this implies that $C_Y= C_{\xi''}$. It is immediate from Definitions that, for  all $\lambda$, $Y(v_\lambda)= \xi''(v_\lambda) = (-1)^{\langle 2\lambda, \rho^\vee \rangle} v_\lambda^\text{low}$. Thus the result follows since, as in Comment \ref{ximake}, an element $\xi \in \widetilde{U_q(\g)}$ is uniquely determined by $C_\xi$ and the action of $\xi$ on all $v_\lambda$. \end{proof}

The following corollaries give us the desired realization of the unitarized $R$ matrix.

\begin{Corollary} \label{Rxi11}
The unitarized R matrix acts on a tensor product $V \otimes W$ by
\begin{equation}
\bar{R} (v \otimes w )= (\xi_V''^{-1} \otimes \xi_W''^{-1}) \circ \xi_{V \otimes W}'' (v \otimes w)
\end{equation}
\end{Corollary}

\begin{proof}
This follows from Theorem \ref{answer1} since, by Proposition \ref{adanddiag}, $Y=\xi''$.
 \end{proof}

This is not quite what we were looking for since, for $v \in V_\lambda$, the relationship between $\xi(v)$ and $ \xi''(v) $ depends on the weight of $ v $, not just on $\lambda$.  To fix this problem, define $\xi' : V_\lambda \rightarrow V_\lambda $ by $ \xi'(v) = i^{2\langle\lambda, \rho^\vee\rangle} \xi(v)$. Notice that $\langle\lambda, \rho^\vee\rangle$ is in general only a half integer, so multiples of $i$ do appear.
We immediately deduce the following.
\begin{Corollary} \label{answer2}
The unitarized R matrix acts on a tensor product $V \otimes W$ by
\begin{equation}
\bar{R} (v \otimes w) = (\xi_V'^{-1} \otimes \xi_W'^{-1}) \circ \xi_{V \otimes W}' (v \otimes w)
\end{equation}
\end{Corollary}

\begin{proof}
Follows from Corollary \ref{Rxi11} by a straightforward calculation.
\end{proof}

\begin{Comment}
Notice that $\xi'_{V_\lambda} \circ \xi'_{V_\lambda} = (-1)^{\langle 2 \lambda, \rho^\vee \rangle} \Id$, and in particular $\xi'$ does not in general square to the identity, as required by Henriques and Kamnitzer. However, their argument can be modified slightly to show directly that $\Flip \circ  (\xi_V'^{-1} \otimes \xi_W'^{-1}) \circ\xi'_{V \otimes W}$ still defines a commutor which satisfies the axioms of a coboundary category. We do not include this, since it follows from the corresponding fact for  $\Flip \circ \bar{R}$.
\end{Comment}

\subsection{The action of $\bar{R}$ on a tensor product of crystal bases} \label{ocb} \label{RC}
This section contains our second main result, namely an explicit relationship between Drinfeld's commutor and the crystal commutor. Roughly, Theorem \ref{main2} shows that the crystal commutor is the crystal limit of Drinfeld's commutor (modulo some signs).  

Recall that if $ (\LL, B) $ is  a crystal basis for a representation $ V$, we have both a linear map $ \xi'_V : V \rightarrow V $ and a map of sets $ \xi_B : B \rightarrow B $ (coming from regarding $ B $ as an abstract crystal).  
The following proposition follows immediately from Theorem \ref{xionbasis}.

\begin{Proposition} \label{xi'onbasis}
Let $ (\LL, B) $ be a crystal basis for a representation $ V $.  Then:
\begin{enumerate}
\item  $ \xi'_V(\LL) = \LL$. 
\item By (i), $ \xi'_V $ gives rise to a map between $ {\xi'}_V^\modq : \LL/q^{-1} \LL \rightarrow \LL/q^{-1}\LL $.  For each $ b \in B $, we have 
\begin{equation} {\xi'}^\modq_V(b) = i^{\langle \lambda, 2 \rho^\vee \rangle} \xi_B(b) \end{equation} 
where $ \lambda $ is the highest weight of the crystal component containing $ b $. \qed
\end{enumerate}
\end{Proposition}

Now let $(\LL, A)$ and $(\mathcal{M}, B)$ be crystal bases for two finite dimensional representations $V$ and $W$ of $U_q(\g)$. The crystal commutor defines a map $\sigma^{HK}_{A,B}: A \otimes B \rightarrow B \otimes A$.  This map comes from Drinfeld's commutor $ \Flip \circ \bar{R} $ in the following sense.

\begin{Theorem} \label{main2}
With the above setup:
\begin{enumerate}
\item $ \Flip \circ \bar{R}(\LL \otimes \mathcal{M}) = \mathcal{M} \otimes \LL $
\item By (i), $ \Flip \circ \bar{R}$ gives rise to a map 
\begin{equation} \Flip \circ \bar{R}^{(mod \, q^{-1}(\LL \otimes \mathcal{M}))} : (\LL \otimes \mathcal{M}) / q^{-1} ( \LL \otimes \mathcal{M}) \rightarrow (\mathcal{M} \otimes \LL) / q^{-1} (\mathcal{M} \otimes \LL). 
\end{equation} 
For all $ a \in A, b \in B$, 
\begin{equation}
 \Flip \circ \bar{R}^{(mod \, q^{-1} (\LL \otimes \mathcal{M}))} (a \otimes b) = (-1)^{\langle \lambda + \mu - \nu, \rho^\vee \rangle}  \sigma^{HK}_{A, B}(a \otimes b)
\end{equation}
where $ \lambda, \mu$ and  $\nu $ are the highest weights of the components of $A, B$ and $A \otimes B$ containing $ a, b$ and $ a \otimes b $ respectively.
\end{enumerate}

\end{Theorem}

\begin{proof}
By Corollary \ref{answer2} and Proposition \ref{xi'onbasis} applied to the crystal bases $ (\LL, A), (\mathcal{M}, B)$ and $(\LL \otimes \mathcal{M}, A \otimes B) $:
\begin{equation}
 \Flip \circ \bar{R}(\LL \otimes \mathcal{M}) = \text{Flip} \circ ({\xi'}_V^{-1} \otimes {\xi'}_W^{-1}) \circ \xi'_{V \otimes W} (\LL \otimes \mathcal{M}) = \text{Flip}(\LL \otimes \mathcal{M}) = \mathcal{M} \otimes \LL.
\end{equation}
This establishes (i).

(ii) follows directly from (i) and Proposition \ref{xi'onbasis}.
\end{proof}
 
Note that consistently working modulo the lattices, one can see that the coboundary properties of Drinfeld's commutor $ \sigma^{D} $ are transferred to the crystal commutor.  Of course it is very easy to prove the coboundary properties of the crystal commutor directly, but we feel this gives some explanation as to why these properties arise.

\section{The crystal commutor and Kashiwara's involution} \label{commutorKashiwara_section}

In this section we propose a new definition for the crystal commutor, and show that it agrees with Henriques and Kamnitzer's definition in the finite type case. This has appeared in \cite{commutorKashiwara}. The advantage of this new definition is that it is well defined for any symmetrizable Kac-Moody algebra. The disadvantage is that we do not know how to prove directly that it gives a coboundary structure on the category of $U_q(\g)$-crystals, so this has only been established for finite type $\g$. Our new definition makes use of Kashiwara's involution on the "Verma crystal" $B_\infty$, so we first review this construction.

\subsection{Kashiwara's involution on $B_\infty$} \label{kinv}
Here we give a purely combinatorial description of $B_\infty$ and of Kashiwara's involution. As we mentioned in Section \ref{kcom_intro}, it is perhaps more natural to think of Kashiwara's involution as coming from a related involution on $U_q^-(\g)$. For more on this point of view we refer the reader to \cite{Kashiwara:1995}.
For any dominant weights $\lambda$ and $\gamma$, there is an inclusion of crystals $ B_{\lambda+\gamma} \rightarrow B_\lambda \otimes B_\gamma$ which sends $b_{\lambda+\gamma}$ to $b_\lambda \otimes b_\gamma$. The following is immediate from the tensor product rule:

\begin{Lemma} \label{in_gamma}
The image of the inclusion $ B_{\lambda+\gamma} \rightarrow B_\lambda \otimes B_\gamma $ contains all elements of the form $ b \otimes b_\gamma $ for $ b \in B_\lambda $.  \qed
\end{Lemma}

Lemma \ref{in_gamma} defines a map $ \iota_{\lambda}^{\lambda + \gamma} : B_\lambda \rightarrow B_{\lambda + \gamma} $ which is $ e_i $ equivariant and takes $ b_\lambda $ to $ b_{\lambda + \gamma} $.  
These maps make $\{ B_\lambda \}$ into a directed system, and the limit of this system is $ B_\infty $.  There are $ e_i $ equivariant maps $ \iota_\lambda^\infty : B_\lambda \rightarrow B_\infty $. When there is no danger of confusion we denote $\iota_\lambda^\infty$ simply by $\iota$.  

The crystal $ B_\infty $ has additional combinatorial structure. In particular, we will need:
\begin{enumerate}
\item The map $\tau : B_\infty \rightarrow \Lambda_+ $ defined by  $\tau(b) = \min \{ \lambda :  b \in \iota (B_\lambda) \}.$
\item The map  $ \varepsilon : B_\infty \rightarrow \Lambda_+ $ given by, for any $b \in B_\infty$ and any $\lambda$ such that $b \in \iota_\lambda^\infty(B_\lambda)$, $\epsilon(b)= \epsilon(\iota_\lambda^{\infty -1}(b))$, where $\epsilon$ is defined on $B_\lambda$ as in Section \ref{notation}. 
\item Kashiwara's involution $*$ (for the construction of this involution see \cite[Theorem 2.1.1]{Kashiwara:1993}).
\end{enumerate}

These maps are related by the following result of Kashiwara.

\begin{Proposition}[\cite{Kashiwara:1995}, Prop. 8.2] \label{th:prop}
Kashiwara's involution preserves weights and satisfies 
\begin{equation}
\tau(*b) = \varepsilon(b), \quad \varepsilon(*b) = \tau(b). 
\end{equation} 
\end{Proposition}

\begin{Comment}
All of this combinatorial structure can be seen easily using the MV polytope model \cite{Kamnitzer:2006}.  The inclusions $ \iota $ correspond to translating polytopes.  The maps $ \tau $ and $\varepsilon $ are given by counting the lengths of edges coming out of the top and bottom vertices.  The involution $ * $ corresponds to negating a polytope.  From this description, the proof of the above proposition is immediate.
\end{Comment}

\subsection{A definition of the crystal commutor using Kashiwara's involution} \label{defK}

Recall that the crystal commutor is a natural system of crystal isomorphisms $\sigma^{HK}_{B, C}: B \otimes C \rightarrow  C \otimes B$, for all integrable highest weight $U_q(\g)$ crystals $B$ and $C$. It is uniquely defined once we know what it  does to highest weight elements of every tensor product 
$B_\lambda \otimes B_\mu$ of irreducible crystals. These are all of the form $b_\lambda \otimes c$, and must be sent to highest
weight elements of $B_\mu \otimes B_\lambda$, which in turn are of the form $b_\mu \otimes b$. 

Let $b_\lambda \otimes c$ be highest weight in $B_\lambda \otimes B_\mu$.
By the tensor product rule for crystals, $ \varepsilon(c) \le \lambda$, so, by Proposition \ref{th:prop}, $\tau(*c) \leq \lambda$, or equivalently $*c \in \iota (B_\lambda)$. For this reason $*c$ can be considered an element of $B_\lambda$. Using Proposition \ref{th:prop} once more, we see that $b_\mu \otimes *c$ must be highest weight. Hence there is a unique isomorphism of crystals $B_\lambda \otimes B_\mu \rightarrow B_\mu \otimes B_\lambda$ which takes each highest weight element $b_\lambda \otimes c$ to $b_\mu \otimes *c$. The following shows that this isomorphism is equal to the crystal commutor.

\begin{Theorem}
\label{maintheorem_commutor}
If $ b_\lambda \otimes c$ is a highest weight element in $B_\lambda \otimes B_\mu$, then
$\sigma^{HK}_{B_\lambda,B_\mu} (b_\lambda \otimes c)= b_\mu \otimes *c$.

\end{Theorem}

 \subsection{Proof of Theorem \ref{maintheorem_commutor}}
 \label{proof}

One of the main tools we will need is the notion of Kashiwara data (also called string data), first studied by Kashiwara (see for example \cite{Kashiwara:1995} section 8.2).
Fix a reduced word ${\bf i}$ for $w_0$, by which we mean $ \wi = (i_1, \dots, i_m) $, where each $ i_k $ is a node of the Dynkin diagram, and $ w_0 = s_{i_1} \cdots s_{i_m} $. The downward Kashiwara data for $b \in B_\lambda$ with respect to $ \wi $ is the sequence of non-negative integers $(p_1, \ldots p_m)$ defined by
\begin{gather}
p_1 := \varphi_{i_1}(b), \quad p_2 := \varphi_{i_2}(f_{i_1}^{p_1} b), \quad \dots, \quad p_m := \varphi_{i_m}(f_{i_{m-1}}^{p-1} \ldots f_{i_1}^{p_1} b).
\end{gather}
That is, we apply the lowering operators in the direction of $i_1$ as far as we can, then in the direction $i_2$, and so on. The following result is due to Littelmann \cite[section 1]{Littelmann:1998}.

\begin{Lemma} \label{th:bottom} After we apply these steps, we reach the lowest element of the crystal. That is:
\begin{equation} f_{i_m}^{p_m} \ldots f_{i_1}^{p_1} b = b_\lambda^{low}. \end{equation}
Moreover, the map $B_\lambda \rightarrow \mathbb{N}^m$ taking $b \rightarrow (p_1, \ldots p_m)$
is injective. \qed
\end{Lemma}

Similarly, the upwards Kashiwara data for $b \in B_\lambda$ with respect to ${\bf i}$ is the sequence $(q_1, \ldots q_m)$ defined by
\begin{gather}
q_1 := \varepsilon_{i_1}(b), \quad q_2 := \varepsilon_{i_2}(e_{i_1}^{q_1} b), \quad \dots, \quad q_m := \varepsilon_{i_m}(e_{i_{m-1}}^{q-1} \ldots e_{i_1}^{q_1} b).
\end{gather}

We introduce the notation $ w_k^\wi := s_{i_1} \cdots s_{i_k} $.

\begin{Lemma} \label{KDL}  In the crystal $ B_\lambda $, we have the following:
\begin{enumerate}
\item The downward Kashiwara data for $b_\lambda$ is given by $p_k =\langle w_{k-1}^\wi \cdot \alpha_{i_k}^\vee , \lambda \rangle.$ 
\item For each $k$, $\varepsilon_{i_k}  ( f_{i_{k-1}}^{p_{k-1}}  \ldots f_{i_1}^{p_1} b_\lambda )= 0.$
\end{enumerate}
\end{Lemma}

\begin{proof} Let $(p_1, \ldots p_m)$ be the downwards Kashiwara data for $b_\lambda$, and
let $\mu_k$ be the weight of $ f_{i_k}^{p_k} \ldots f_{i_1}^{p_1} b_\lambda$. Since $f_{i_k}^{p_k} \ldots f_{i_1}^{p_1} b_\lambda$ is the end of an $\alpha_{i_k}$ root string, we see that
\begin{equation} \label{reflect} \mu_k= s_{i_k} \cdot \mu_{k-1} -  a_k\alpha_{i_k}, \end{equation} 
 where $a_k = \varepsilon_{i_k}  ( f_{i_{k-1}}^{p_{k-1}}  \ldots f_{i_1}^{p_1} b_\lambda )$.  Using this fact at each step,
\begin{equation} \mu_m = w_0 \cdot \lambda - \sum_{k=1}^m a_k s_{i_m} \ldots s_{i_{k+1}} \cdot \alpha_{i_k}. \end{equation}
By Lemma \ref{th:bottom}, we know that $f_{i_m}^{p_m} \ldots f_{i_1}^{p_1} b_\lambda = b_\lambda^{low}$, so that $\mu_m =  w_0 \cdot \lambda$. 
Hence 
\begin{equation}  \sum_{k=1}^m a_k s_{i_m} \ldots s_{i_{k+1}} \cdot \alpha_{i_k} = 0. \end{equation}
Now, $s_{i_m} \cdots s_{i_k} $ is a reduced word for each $k$, which implies that
$s_{i_m} \ldots s_{i_{k+1}} \alpha_{i_k}$ is a positive root. Thus
each $a_k$ is zero, proving part (ii).

Equation (\ref{reflect}) now shows that $\mu_k = s_{i_k} \ldots s_{i_1} \cdot \lambda$, for all $k$. In particular that  $f_{i_k}^{p_k}$ must perform the reflection $s_{i_k}$ on the weight $s_{i_{k-1}} \ldots s_{i_1} \cdot \lambda .$ Therefore, 
\begin{equation} p_k = \langle \alpha_{i_k}^\vee  , s_{i_{k-1}} \ldots s_{i_1} \cdot \lambda \rangle
= \langle w_{k-1}^{\bf i} \alpha_{i_k}^\vee, \lambda \rangle. \end{equation}
\end{proof}

\begin{Lemma} Let $b_\lambda \otimes c$ be a highest weight element of $B_\lambda \otimes B_\mu$. 
Let $b \otimes b_\mu^{\text{low}}$ be the lowest weight element of the component containing 
$b_\lambda \otimes c$. Let $ (p_1, \ldots p_m)$ be the downward Kashiwara data for $c$ with respect to ${\bf i}$, and $(q_1, \ldots q_m)$ the upward Kashiwara data for $b$ with respect to ${\bf i}^\rev := (i_m, \dots, i_1) $. Then, for all k, 
$p_k + q_{m-k+1} = \langle w_{k-1}^\wi \cdot \alpha_{i_k}^\vee ,  \nu\rangle,$ where $\nu= \wt (b_\lambda \otimes c)$.
\label{updowndata}
\end{Lemma}

\begin{proof}
Let $r_k= \langle w_{k-1}^\wi \cdot \alpha_{i_k}^\vee , \nu\rangle.$ By part (i) of Lemma \ref{KDL}, $(r_1, \ldots r_m)$ is the downward Kashiwara data for $b_\lambda \otimes c$. Define $b_k \in B_\lambda$ and $c_k \in B_\mu$ by 
$b_k \otimes c_k = f_{i_k}^{r_k} \ldots f_{i_1}^{r_1} ( b_\lambda \otimes c )$. Part (ii) of Lemma \ref{KDL}, along with the definition of Kashiwara data, shows that, for each $1 \leq k \leq m$, 
\begin{equation} e_{i_k} (b_{k-1} \otimes c_{k-1})= 0 \hspace{0.2in} \mbox{and} \hspace{0.2in} f_{i_k} (b_k \otimes c_k) = 0. \end{equation}
In particular, the tensor product rule for crystals implies
\begin{equation} e_{i_k} b_{k-1} = 0 \hspace{0.3in} \mbox{and} \hspace{0.2in} 
f_{i_k} c_k  = 0. \end{equation}
Define $p_k$ to be the number of times $f_{i_k}$ acts on $c_{k-1}$ to go from $b_{k-1} \otimes c_{k-1}$ to
$b_k \otimes c_k$, and $q_{m-k+1}$ to be the number of times $f_{i_k}$ acts on $b_{k-1}$. Since $f_{i_k} c_k=0$, we see that $\varphi_{i_k} (c_{k-1}) = p_k$. Hence, by definition $(p_1, \ldots p_m)$ is the downward Kashiwara data for $c$ with respect to $ \wi $. Similarly,
$e_{i_k} b_{k-1} = 0$, so $\varepsilon_{i_k} (c_k) = q_{m-k+1}$. By Lemma \ref{th:bottom}, $b_m= b$, so this implies that $(q_1, \ldots q_m)$ is the upward Kashiwara
data for $b$ with respect to $ \wi^\rev $. 
Since $p_k + q_{m-k+1}$ is the number of times that $ f_{i_k} $ acts on $ b_{k-1} \otimes c_{k-1} $ to reach $b_k \otimes c_k$, we see that $ p_k + q_{m-k+1} = r_k $.
\end{proof}

Let $b_\lambda \otimes c$ be a highest weight element in $B_\lambda \otimes B_\mu$. As discussed in section \ref{defK}, $*c$ can be considered as an element of $B_\lambda$.

\begin{Lemma} Define $\nu = wt (b_\lambda \otimes c)$. Let $(p_1, \ldots p_m)$ be the downward Kashiwara data for $c \in B_\mu$ with respect to ${\bf i}.$ Let  $(q_1, \ldots q_m)$ be the downward Kashiwara data for $*c \in B_\lambda$ with respect to the decomposition 
$\theta ({\bf i}^\rev) := (\theta(i_m), \dots, \theta(i_1))$ of $w_0$. Then, for all $k$, 
$p_k +q_{m-k+1} = \langle w_{k-1}^\wi \cdot \alpha_{i_k}^\vee, \nu \rangle.$
\label{*data}
\end{Lemma}

\begin{proof}
The proof will depend on results from \cite{Kamnitzer:2006} on the MV polytope model for crystals.  In particular, within this model it is easy to express Kashiwara data and the Kashiwara involution.

Let $ P = P(M_\bigdot) $ be the MV polytope of weight $ (\nu - \lambda, \mu) $ corresponding to $ c $.  Then by \cite [Theorem 6.6]{Kamnitzer:2006}, 
\begin{equation} \label{eq:pM} 
p_k = M_{w_{k-1}^\wi \cdot \Lambda_{i_k}} - M_{w_k^\wi \cdot \Lambda_{i_k}}.
\end{equation}

Now, consider $ P $ as a stable MV polytope (recall that this means that we only consider it up to translation).  Then by \cite[Theorem 6.2]{Kamnitzer:2006}, we see that $ *(P) = -P $.

The element $ \iota_\mu^{-1} *\iota_\lambda(c) \in B_\lambda $ corresponds to the MV polytope $ \nu - P $ and hence has BZ datum $ N_\bigdot $,  where $M_\bigdot$ and $N_\bigdot$ are related by
\begin{equation} \label{eq:N}
M_\gamma = \langle \gamma, \nu \rangle + N_{-\gamma}.
\end{equation}

Let $\wi' = \theta(\wi^{\rev}) $.  
Then,
\begin{equation}
 -w_{k-1}^{\wi} \cdot \Lambda_{i_k} = w_{m-k +1}^{\wi'} \cdot \Lambda_{i'_{m-k+1}} \text{ and } 
 -w_k^{\wi} \cdot \Lambda_{i_k} = w_{m-k}^{\wi'} \cdot \Lambda_{i'_{m-k+1}}.
 \end{equation}
Combining the last 3 equations, we see that
\begin{equation} \label{eqpk}
\begin{aligned}
p_k &= \langle w_{k-1}^{\wi} \cdot \Lambda_{i_k}, \nu \rangle + N_{-w_{k-1}^{\wi} \cdot \Lambda_{i_k}} - \langle w_k^{\wi} \cdot \Lambda_{i_k}, \nu \rangle - N_{-w_k^{\wi} \cdot \Lambda_{i_k}} \\
&= N_{w_{m-k+1}^{\wi'} \cdot \Lambda_{i'_{m-k+1}}} - N_{w_{m-k}^{\wi'} \cdot \Lambda_{i'_{m-k+1}}} + \langle w_{k-1}^{\wi} \cdot \Lambda_{i_k} - w_k^{\wi} \cdot \Lambda_{i_k}, \nu \rangle.
\end{aligned}
\end{equation}
Applying \cite[Theorem 6.6] {Kamnitzer:2006} again,
\begin{equation} \label{eqqk}
q_k = N_{w_{k-1}^{\wi'} \cdot \Lambda_{i'_k}} - N_{w_k^{\wi'} \cdot \Lambda_{i'_k}}.
\end{equation}
We now add equation (\ref{eqpk}) and (\ref{eqqk}), substituting $m-k+1$ for $k$ in the second equation, to get
\begin{equation} \label{eq:p+q} p_k+q_{m-k+1}= \langle w_{k-1}^{\wi} \cdot \Lambda_{i_k} - w_k^{\wi} \cdot \Lambda_{i_k}, \nu \rangle = \langle w_{k-1}^{\wi} \cdot (\Lambda_{i_k} - s_{i_k} \cdot \Lambda_{i_k}), \nu \rangle = \langle w_{k-1}^{\wi} \cdot \alpha_{i_k}^\vee, \nu \rangle.
\end{equation}
\end{proof}

\begin{proof}[Proof of Theorem \ref{maintheorem_commutor}]
We know that $\sigma^{HK}_{B_\lambda, B_\mu} (b_\lambda \otimes c) = b_\mu \otimes b$ for some $b \in B_\lambda$. By the definition of $\sigma^{HK}_{B_\lambda, B_\mu}$ (see  Section \ref{commutor_section}),  we see that 
\begin{equation}
\xi (b_\lambda \otimes c)=  (\xi \circ \xi) (\mathrm{Flip}(\sigma^{HK}(b_\lambda \otimes c))) =  \xi(b) \otimes b_\mu^{low}.
\end{equation}
In particular, $\xi(b)  \otimes b^{low}_\mu$ is the lowest weight element of the component of $B_\lambda \otimes B_\mu$ containing $b_\lambda \otimes c$. 

Fix a reduced expression ${\bf i}= (i_1, \ldots, i_m)$ for $w_0$. Let $(p_1, \ldots, p_m)$ be the downward Kashiwara data for $c$ with respect to ${\bf i}$, and let $(q_1, \ldots q_m)$ be the downward Kashiwara data for $b$ with respect to $\theta({\bf i}^\rev):= (\theta(i_m), \cdots, \theta(i_1))$. 
Notice that $(q_1, \ldots q_m)$ is also the upward Kashiwara data for $\xi (b)$ with respect to ${\bf i}^{\rev}:= (i_m, \ldots i_1)$, since $\xi$ interchanges the action of $f_i$ and $e_{\theta(i)}$. Hence by Lemma \ref{updowndata}, 
\begin{equation} \label{sequ}  p_k+ q_{m-k+1} = \langle w_{k-1}^\wi \cdot \alpha_{i_k}^\vee, \nu \rangle \end{equation}
for all $k$, where $\nu$ is the weight of $b_\lambda \otimes c$.

As discussed in Section \ref{defK}, $*c \in \iota(B_\mu)$, and so can be considered as an element of $B_\mu$. Let $(q_1', \ldots q_m')$ by the downward Kashiwara data for $*c \in B_\mu$ with respect to $\theta ( {\bf i}^{\rev} )$. 
By Lemma \ref{*data} we have 
\begin{equation} \label{*equ} p_k + q_{m-k+1}' = \langle w_{k-1}^\wi \cdot \alpha_{i_k}^\vee, \nu \rangle . \end{equation}
Comparing equations (\ref{sequ}) and (\ref{*equ}), $q_k = q_k'$ for each $1 \leq k \leq m$. That is, the downward Kashiwara data for $b$ and $*c$ with respect to $\theta({\bf i}^{\rev})$ are identical.
Hence by Lemma \ref{th:bottom}, $ b = *c $. \end{proof}

\section{A new formula for the $R$-matrix of a symmetrizable Kac-Moody algebra} \label{RTheta_section}

Let $\g$ be a finite type complex simple Lie algebra, and let $U_q(\g)$ be the corresponding quantized universal enveloping algebra. As we discussed in Section \ref{autint_intro}, Kirillov-Reshetikhin \cite{KR:1990} and Levendorskii-Soibelman \cite{LS:1991} developed a formula for the universal R-matrix 
\begin{equation} \label{KReq} R= (X^{-1} \otimes X^{-1}) \Delta(X), \end{equation}
where $X$ belongs to a completion of $U_q(\g)$. The element $X$ is constructed using the braid group element $T_{w_0}$ corresponding to the longest word of the Weyl group, and as such only makes sense when $\g$ is of finite type. 

The element $X$ in Equation (\ref{KReq}) defines a vector space endomorphism $X_V$ on each representation $V$, and in fact $X$ is defined by this system $\{ X_V \}$ of endomorphisms. With these point of view, Equation (\ref{KReq}) is equivalent to the claim that, for any finite dimensional representations $V$ and $W$ of $U_q(\g)$, and any $u \in V \otimes W$,
\begin{equation} \label{KReq2}
R (u) =(X_V^{-1} \otimes X_W^{-1}) X_{V \otimes W} (u).
\end{equation}
The endomorphisms $X_V$ permute weight spaces according to the longest element of the Weyl group, and as such only makes sense when $\g$ is of finite type.

In the present work we replace $X_V$ with an endomorphism $\Theta_V$ which preserves weight spaces. We show that, for any symmetrizable Kac-Moody algebra $\g$, and any integrable highest weight representations $V$ and $W$ of $U_q(\g)$, the action of the universal $R$-matrix on $u \in V \otimes W$ is given by
\begin{equation} \label{teq1}
R (u)= (\Theta_V^{-1} \otimes \Theta_W^{-1}) \Theta_{V \otimes W} (u). 
\end{equation}
There is a technical difficulty because $C_{\Theta}$ is not linear over the base field $\bc(q)$, but instead is compatible with the $``\barr"$ automorphism of $\bc(q)$ which inverts $q$. Consequently, the endomorhsim $\Theta_V$ is $\barr$-linear instead of linear, and in fact depends on a choice of a $\barr$-involution on $V$. In the case of an irreducible representation $V_\lambda$, the endomorphism $\Theta_{V_\lambda}$ depends only on a choice of highest weight vector $v_\lambda \in V_\lambda$. To make Equation (\ref{teq1}) precise we must define the action of $\Theta$ on $V_{\lambda} \otimes V_\mu$ in terms of chosen highest weight vectors $v_\lambda \in V_\lambda$ and $v_\mu \in V_\mu$, and then show that the composition $(\Theta_{V_\lambda}^{-1} \otimes \Theta_{V_\mu}^{-1}) \Theta_{V_\lambda \otimes V_\mu}$ does not depend on $v_\lambda$ and $v_\mu$. 

The system of endomorphisms $\Theta$ was previously studied in \cite{RTheta}, where it was used to construct the universal $R$-matrix when $\g$ is of finite type. Essentially we have extended this previous work to include all symmetrizable Kac-Moody algebras. However, the action of $\Theta$ on a tensor product in defined differently here than in  \cite{RTheta}, so the constructions of $R$ are a-priori not identical, and we have not in fact proven that the construction in \cite{RTheta} gives the universal $R$-matrix in all cases.

\subsection{The endomorphism $\Theta$ and the system of isomorphisms $\sigma^\Theta$.} \label{maketheta}

We now introduce a $\bc$-algebra automorphism $C_\Theta$ of $U_q(\g)$: 
\begin{equation}
\begin{cases}
C_\Theta (E_i) =   E_{i} K_i^{-1} \\
C_\Theta (F_i) =   K_i F_{i} \\
C_\Theta (K_i) = K_{i}^{-1} \\
C_\Theta (q) = q^{-1}.
\end{cases}
\end{equation}
Notice that $C_\Theta$ inverts $q$, so it is not a $\bc(q)$ algebra automorphism, but is instead $\barr$ linear. That is, it is compatible which the $\barr$-involution of $\bc(q)$ which inverts $q$.

One can check that $C_\Theta$ is an algebra involution and a coalgebra anti-involution. It is perhaps useful to work over ${\Bbb Q}(q)$ instead of $\bc(q)$, and imagine that $q$ is specialized to a complex number on the unit circle (although not a root of unity). Thus bar-involution should be thought of as complex conjugation. The following defines an involution of $U_q(\g)$ compatible with bar-involution on $\bc(q)$, and will prove very useful:

\begin{Definition} \label{bardef}
$C_\barr: U_q(\g) \rightarrow U_q(\g)$ is the $\bc$-algebra involution defined by
\begin{equation}
\begin{cases}
C_\barr{q} = q^{-1}\\
C_\barr{K}_i = K_i^{-1}  \\
C_\barr{E}_i = E_i \\
C_\barr{F}_i = F_i.
\end{cases}
\end{equation}
\end{Definition}

We must define a $\bc$-vector space automorphism $\Theta_{V_\lambda}$ of each $V_\lambda$ which is compatible with $C_\Theta$. Since $C_\Theta$ is bar-linear instead of linear, it will be necessary to first fix an automorphism of $V_\lambda$ compatible with $\barr$-involution on $\bc(q)$. If we fix a highest weight vector $v_\lambda \in V_\lambda$, there is a unique involution  $\barr_{(V_\lambda, v_\lambda)}$ that fixes $v_\lambda$, and such that the following diagram commutes:

\begin{equation*}
\xymatrix{
V_\lambda  \ar@(dl,dr) \ar@/ /[rrr]^{\barr_{(V_\lambda, v_\lambda)}} &&& V_\lambda \ar@(dl,dr) \\
U_q(\g)  \ar@/ /[rrr]^{\C_{\barr}} &&&   U_q(\g). \\
}
\end{equation*}
When it does not cause confusion we will denote $\barr_{(V_\lambda, v_\lambda)} (v)$ by $\bar{v}$. In general we will use the notation $(V, bar_V)$ to denote a representation $V$ together with a chosen involution on that representation compatible with the involution $C_\barr$ of $U_q(\g)$.

\begin{Definition} \label{Thetadef}
The action of $\Theta_{(V_\lambda, v_\lambda)}$ on $V_\lambda$ is defined to send $v_\lambda$ to $q^{-( \lambda, \lambda)/2 + ( \lambda, \rho )} v_\lambda$, and to be compatible with $C_\Theta$. \end{Definition}

 It is straightforward to check that for any weight vector $b \in V_\lambda$ of weight $\mu$,  $\Theta_{(V_\lambda, B_\lambda)} (b) = q^{-( \mu, \mu)/2 + ( \mu, \rho )} \bar{b}$. Thus we can define:
 
 \begin{Definition} \label{good_theta_def}  For any representation $V$ with a chosen involution $\barr_V$ compatible with $C_\barr$, define $\Theta_{(V, \barr_V)}$ to act on an weight vector $b$ by
 \begin{equation}
 \Theta_{(V, \barr_V)} (b)=  q^{-( \mu, \mu)/2 + ( \mu, \rho )} \barr_{V}(b).
 \end{equation}
Note that by the above discussion $\Theta_{(V, \barr_V)}$ is compatible with $C_\Theta$.
\end{Definition}

The construction described in Section \ref{autint_intro} uses the action of $\xi$ on $V \otimes W$. Thus we will need to define how $\Theta$ acts on a tensor product. In particular, we need a well defined notion of tensor product in the category whose objects consist of a pair $(V, \barr_V)$ of a representation and an involution compatible with $C_\barr$. It is enough to define an involution on $V_\lambda \otimes V_\mu$ compatible with $C_\barr$ given highest weight vectors $v_\lambda \in V_\lambda$ and $v_\mu \in V_\mu$. The following naive guess is not compatible with $C_\barr$, but will be useful in defining an involution which is.

\begin{Definition} \label{bbd}
Fix $(V, \barr_V)$ and $(W, \barr_W)$, where $\barr_V$ and $\barr_W$ are involutions of $V$ and $W$ compatible with $C_\barr$. Let $(\barr_V \otimes \barr_W)$ be the vector space involution on $V \otimes W$ defined by $f(q) v \otimes w \rightarrow f(q^{-1}) \bar{v} \otimes \bar{w}$ for all $f(q) \in \bc(q)$ and $v \in V, w \in W$.
\end{Definition}

\begin{Comment}
It is straightforward to check that the action of $(\barr_V \otimes \barr_W)$ on a vector in $V \otimes W$ does not depend on its expression as a sum of element of the form $f(q) v \otimes w$. The resulting map is a $\bc$-linear involution, and in particular is invertible.
\end{Comment}

\begin{Lemma} \label{tensor_involution}
Let $v_\nu$ be a singular weight vector in $V_\lambda \otimes V_\mu$, and write
\begin{equation}
v_\nu= \sum_{j=1}^N b_j \otimes c_j,
\end{equation}
where each $b_j$ is a weight vector of $V_\lambda$, and each $c_j$ is a weight vector of $V_\mu$. Then
\begin{equation} \label{tbar}
\barr(v_\nu):= \sum_{j=0}^N q^{(\mu,\mu)-(\wt(c_j),\wt(c_j))+ 2(\mu-\wt(c_j), \rho)} \bar{b}_{j} \otimes  \bar{c}_{j}
\end{equation}
is also singular.
\end{Lemma}

\begin{proof}
For any vector $v \in B_\lambda \otimes B_\mu$, introduce that notation $v^\beta$ to be indicate the part of $v$ in $V_\lambda(\nu-\beta) \otimes V_\mu(\beta)$ (here $V_\mu(\beta)$ denotes the $\beta$ weight space in $V_\mu$).

Fix $i \in I$. The vector $v_\nu$ is singular, so $E_i v_\nu =0$ and hence $(E_i v_\nu)^\beta=0$ for all $\beta$. One can easily calculate that:
\begin{equation}
0= (E_i v_\nu)^\beta = \sum_{\wt(c_j) = \beta} q^{(\beta, \alpha_i)} E_i b_j \otimes c_j + \sum_{\wt(c_j)=\beta-\alpha_i} b_j \otimes E_i c_j.
\end{equation}
Using Equation (\ref{tbar}): 
\begin{align}
  (E_i \barr(v_\nu))^\beta & = \sum_{\wt(c_j)  = \beta} q^{(\mu,\mu)-(\beta,\beta)+2(\mu-\beta, \rho)} q^{(\beta, \alpha_i)} E_i \bar{b}_j \otimes \bar{c}_j \\
\nonumber & \hspace{0.2in} + \sum_{\wt(c_j)=\beta-\alpha_i} q^{(\mu,\mu)-(\beta-\alpha_i,\beta-\alpha_i)+2(\mu-\beta+\alpha_i, \rho)} \bar{b}_j \otimes E_i \bar{c}_j \\
&=  q^{(\mu,\mu)-(\beta-\alpha_i,\beta-\alpha_i)+2(\mu-\beta+\alpha_i, \rho)}  \times \\
\nonumber & \hspace{0.5in} \times \left(  \sum_{\wt(c_j) = \beta} q^{-(\beta, \alpha_i)} E_i \bar{b}_j \otimes \bar{c}_j + \sum_{\wt(c_j)=\beta-\alpha_i} \bar{b}_j \otimes E_i \bar{c}_j \right) \\
&=  q^{(\mu,\mu)-(\beta-\alpha_i,\beta-\alpha_i)+2(\mu-\beta+\alpha_i, \rho)}  (\barr_{(V_\lambda, v_\lambda)} \otimes \barr_{(V_\mu, v_\mu)}) (E_i v_\nu)^\beta,
\end{align}
where $ (\barr_{(V_\lambda, v_\lambda)} \otimes \barr_{(V_\mu, v_\mu)}) $ is the involution from Definition \ref{bbd}.
But $E_i(v_\nu)^\beta=0$, and $ (\barr_{(V_\lambda, v_\lambda)} \otimes \barr_{(V_\mu, v_\mu)})$ is an invertible linear map on $V_\lambda \otimes V_\mu$, so this implies that $(\barr_{(V_\lambda, v_\lambda)} \otimes \barr_{V_\mu, v_\mu}) E_i(v_\nu)^\beta=0$. Since this holds for all $\beta$, we see that $E_i\barr(v_\nu)=0$. Hence $\barr(v_\nu)$ is a singular vector.
\end{proof}

We can now define an involution $\barr_{(V_\lambda, v_\lambda) \otimes (V_\mu, v_\mu)}$ on $V_\lambda \otimes V_\mu$ which is compatible with $C_\barr$, and which depends only on $v_\lambda$ and $v_\mu$.

\begin{Definition} \label{tb2}
Let $\barr_{(V_\lambda, v_\lambda) \otimes (V_\mu, v_\mu)}$ be the unique involution on $V_\lambda \otimes V_\mu$ which agrees with the involution $\barr$ from Lemma \ref{tensor_involution} on singular vectors, and is compatible with $C_\barr$.
\end{Definition}

This can be extended by naturallity to give a $\barr$-involution on $(V, \barr_V) \otimes (W, \barr_W)$. Thus Definition \ref{good_theta_def} gives a well defined action of $\Theta$ on the tensor product, and we can proceed to define $\sigma^\Theta$ following the methods described in Section \ref{autint_intro}.

\begin{Definition} \label{sigma_theta_def}
$\sigma^\Theta_{V,W} := \Flip \circ (\Theta_V^{-1} \otimes \Theta_W^{-1}) (\Theta_{V \otimes W}): V \otimes W \rightarrow W \otimes V,$

\noindent where the action of $\Theta_{V \otimes W}$ is defined using $\barr_{(V_\lambda, v_\lambda) \otimes (V_\mu, v_\mu)}$ from Definition \ref{tb2}.
\end{Definition}
As discussed in Section \ref{autint_intro}, these $\sigma^\Theta_{V,W}$ combine to give a natural family of isomorphisms, which we will denote by $\sigma^\Theta$. Note that $\sigma^\Theta$ is defined using a choice of involution on $V$ and $W$, but in fact, by Proposition \ref{ind_prop} below, it is independent of this choice.

\begin{Comment}
One can check that $\Theta_V$ is an involution, so the inverses in Definition \ref{sigma_theta_def} are in some sense unnecessary. We include them because $\Theta_{V}$ should really be thought of as an isomorphism between $V$ and the module which is $V$ as a $\bc$ vector space, but with the action of $U_q(\g)$ twisted by $C_\Theta$. We have not specified the action of $\Theta$ on this new module. The way the formula is written, $\Theta$ is always acting on $V, W$ or $V \otimes W$ with the usual action, where it has been defined.  \end{Comment}

\begin{Proposition} \label{ind_prop}
$\sigma^\Theta_{V,W}: V \otimes W \rightarrow W \otimes V$ is independent of the choice of involutions $\barr_V$ and $\barr_W$.
\end{Proposition}

\begin{proof}
It is sufficient to consider the case of $\sigma_{V_\lambda, V_\mu}^\Theta$. The involutions on $V_\lambda$ and $V_\mu$ only depend on a choice of highest weight vector in each of $V_\lambda$ and $V_\mu$. It is straightforward to check that re-scaling these highest weight vectors has no effect on $\sigma^\Theta_{V_\lambda, V_\mu}$.
\end{proof}

\subsection{Constructing the standard $R$-matrix} \label{theproof}

In this section we prove that the action of $(\Theta_V^{-1}  \otimes \Theta_W^{-1}) \Theta_{V \otimes W}$ on $V \otimes W$ is equal to the action of the standard $R$-matrix. Equivalently,  $\sigma^{\Theta}$ is the standard braiding. 

Let $(V_\lambda, v_\lambda)$ and $(V_\mu, v_\mu)$ be irreducible representations with chosen highest weight vectors. Every vector $w \in V_\lambda \otimes V_\mu$ can be writen as
\begin{equation}
 w= v_\lambda \otimes c_0 + b_{k-1} \otimes c_1 + \ldots + b_1 \otimes c_{k-1} + b_0 \otimes v_\mu,
\end{equation}
where, for $0 \leq j \leq k-1$, $b_j$ is a weight vector of $V_\lambda$ of weight strictly less then $\lambda$, and $c_j$ a weight vector of $V_\mu$ of weight strictly less then $\mu$. Furthermore, the vectors $b_0 \in V_\lambda$ and $c_0 \in V_\mu$ are uniquely determined by $w$. Thus we can define projections from $V_\lambda \otimes V_\mu$ to $V_\lambda$ and $V_\mu$ as follows:

\begin{Definition} \label{p_def} The projections $p_{\lambda, \mu}^1: V_\lambda \otimes V_\mu \rightarrow V_\lambda$ and $p_{\lambda, \mu}^2: V_\lambda \otimes V_\mu \rightarrow V_\mu$ are given by, for all $w \in  V_\lambda \otimes V_\mu $, 
\begin{align} p_{\lambda,\mu}^1(w):= b_0 \\
p_{\lambda, \mu}^2(w):= c_0.
\end{align}
\end{Definition}

\begin{Lemma} \label{find_highest}
Let $S_{\lambda, \mu}$ be the space of singular vectors in $V_\lambda \otimes V_\mu$. 
The restrictions of the maps $p_{\lambda, \mu}^1$ and $p_{\lambda, \mu}^2$ from Definition \ref{p_def} to $S_{\lambda, \mu}$ are injective. \end{Lemma}

\begin{proof}
It follows from the tensor product rule for crystal bases (see Section \ref{crystal_section}) that the span of $p_{\lambda, \mu}^2(s)$ over singular vectors $s$ of $V_\lambda \otimes V_\mu$ has the same dimension as $S_{\lambda, \mu}$ itself. Hence $p_{\lambda, \mu}^2$ is injective.

To see that $p_{\lambda,\mu}^1$ is also injective, note that, since $C_\Theta$ is a coalgebra anti-automorphism,  $\Flip \circ  (\Theta_{V_\lambda} \otimes \Theta_{V_\mu})$ is an isomorphism of representations. Thus the map 
\begin{equation} s \rightarrow \Flip \circ  (\Theta_{V_\lambda} \otimes \Theta_{V_\mu}) s
\end{equation}
 is a vector space isomorphisms between the space of singular vectors of $V_\lambda \otimes V_\mu$ and the space of singular vectors of $V_\mu \otimes V_\lambda$. Composing this with $p_{\mu, \lambda}^2$, we see that the map
\begin{equation}
s \rightarrow q^{(\mu, \mu)/2 -(\mu, \rho)+(\wt(s)-\mu, \wt(s)-\mu/2)-(\wt(s)-\mu, \rho)} \barr_{V_\lambda} (p_{\lambda, \mu}^1(s))
\end{equation}
is injective. But then $p_{\lambda,\mu}^1$ itself must also be injective.
\end{proof}

The following fact is well known (see for instance \cite{Lusztig:1993}). As in Section \ref{RY}, we will write expressions using $h$, where $q=e^{h}$, and simply note that the action on any (tensor product of) representations of $U_q(\g)$ is well defined (although can sometimes introduce fractional powers of $q$). 

\begin{Lemma} The standard universal $R$ matrix can be expressed as
\begin{equation} \label{R12}
R_h =  \exp \left( h \sum_{i,j} (B^{-1} )_{ij} H_i \otimes H_j \right) \left( 1+ \sum_{\begin{array}{c} \txt{positive integral} \\ \text{weights } \beta \end{array} } X_\beta \otimes Y_\beta \right)
\end{equation}
where each $X_\beta$ is an element of $U_h(\g)$ of weight $\beta$, and $Y_\beta$ is an element of $U_h(\g)$ of weight $-\beta$.  \qed
\end{Lemma}

\begin{Lemma} \label{JLemma} Let $v \in V$ and $w \in W$ be weight vectors. Then 
\begin{equation} \label{Jmain}
\exp \left( h \sum_{i,j} (B^{-1} )_{ij} H_i \otimes H_j \right) (v \otimes w) =
q^{(\wt(v), \wt(w))} (v \otimes w),
\end{equation}
where $q=e^h$.
\end{Lemma}

\begin{proof}
Define 
\begin{equation} \label{J1}
J:=   \exp  \left[ h \left( \frac{1}{2} \sum_{i,j} \left( (B^{-1})_{ij} H_i H_j \right) + H_\rho \right) \right].
\end{equation}
It is straightforward to check that $J$ acts on a weight vector $b$ in any representation $V$ by 
\begin{equation} \label{J2}
J(b)= q^{(\wt(b), \wt(b))/2 + (\wt(b), \rho)} b.
\end{equation}
Note that $J$ can act by fractional powers of $q$, which is why we must adjoint a fixed $k^{th}$ root of $q$ to $\bc(q)$ (see Section \ref{notation}). 
Calculating
$(J^{-1} \otimes J^{-1}) \Delta(J) (v \otimes w)$ using Equation (\ref{J1}) gives the left side of equation (\ref{Jmain}), and calculating  $(J^{-1} \otimes J^{-1}) \Delta(J) (v \otimes w)$ using Equation (\ref{J2}) gives the right side of Equation (\ref{Jmain}). Thus the two sides are equal. Details of these calculations can be found in \cite{Rcommutor}.
\end{proof}

We are now ready to prove our final main result of Chapter \ref{commutor_chapter}:

\begin{Theorem} \label{RTheta_main}
$(\Theta_V^{-1}  \otimes \Theta_W^{-1}) \Theta_{V \otimes W}$ acts on $V \otimes W$ as the standard $R$-matrix.
\end{Theorem}

\begin{proof}
By naturallity it is enough to consider the case when $V$ and $W$ are irreducible.
We will actually show that $\sigma^\Theta$ acts on $V_\lambda \otimes V_\mu$ as the standard braiding, since this is equivalent. Since $C_\Theta$ is a coalgebra anti-autamorphism, using the discussion in Section \ref{autint_intro} we see that $\sigma^\Theta$ is an isomorphism of representations. Since $\Flip \circ R$ is also an isomorphism, it suffices to show that they agree on singular vectors in $V_\lambda \otimes V_\mu$. So, let $v_\nu$ be a singular vector of weight $\nu$, and write  
\begin{equation}
v_\nu= b_\lambda \otimes c_0 + b_{k-1} \otimes c_1 + \ldots + b_1 \otimes c_{k-1} + b_0 \otimes b_\mu
\end{equation}
where, for $0 \leq j \leq k-1$, $b_j$ is a weight vector of $V_\mu$ of weight strictly less then $\mu$. By Definitions \ref{good_theta_def} and \ref{sigma_theta_def}, 
\begin{align}  \sigma^\Theta(v_\nu) 
 &= q^{-(\wt(b), \wt(b))/2  -(\mu, \mu)/2+(\mu+\wt(b), \mu+\wt(b))/2 } b_\mu \otimes b_0 + \ldots \\
& =q^{(\wt(b), \mu)}  b_\mu \otimes b_0 + \ldots,
\end{align}
where $\ldots$ representations terms of the form $c \otimes b$ with $\wt(c) < \mu$.

It follows from Definition \ref{R12} and Lemma \ref{JLemma} that 
\begin{equation}
\Flip \circ R (v_\nu)= q^{(\wt(b), \mu)} b_\mu \otimes b_0 + \ldots,
\end{equation}
where again $\ldots$ representations terms of the form $c \otimes b$ where $\wt(c) < \mu$. Both $\sigma^\Theta(v_\nu)$ and $\Flip \circ R (v_\nu)$ are singular vectors in $V_\mu \otimes V_\lambda$. Hence by Lemma \ref{find_highest} they are equal.
\end{proof}

It is an immediate consequence of Theorem \ref{RTheta_main} and Equation (\ref{QR}) that we have a new expression for Drinfeld's unitarized $R$-matrix:

\begin{Corollary}
Let $\Gamma_{V_\lambda}$ be the $\bc$-linear endomorphism of $(V, \barr_V)$ defined by composing $ \Gamma_{(V, \barr_V)} \circ Q^{-1/2}$, where $Q$ is as in Section \ref{coboundary_section}. Then
$(\Gamma_V^{-1}  \otimes \Gamma_W^{-1}) \Gamma_{V \otimes W}$ acts on $V \otimes W$ as Drinfeld's unitarized $R$-matrix. \qed
\end{Corollary}

\section{Future directions}

The motivating question behind this work was to find a coboundary structure on the category of $\g$-crystals, when $\g$ is any symmetrizable Kac-Moody algebra. This has still not been done, although the results in section \ref{commutorKashiwara_section} lead to a natural conjecture:

\begin{Conjecture} \label{CrystConjecture}
The definition of the crystal commutor from Section \ref{defK} is a coboundary structure.\footnote{This Conjecture has recently been proven by Savage \cite[Theorem 6.4]{Savage:2008}}
\end{Conjecture}

The results from Section \ref{RTheta_section} are partly aimed at proving this conjecture. We have given a construction of the standard braiding of the form $(\Theta^{-1} \otimes \Theta^{-1}) \Delta(\Theta)$, and from there a construction of Drinfeld's coboundary structure of the same form. We hope to use methods similar to Section \ref{Rcommutor_section} to show that Drinfeld's coboundary structure acts in a nice way on crystal bases, and that this will lead to a proof of Conjecture \ref{CrystConjecture}. However, $\Theta$ itself does not preserve crystal lattices, so the methods cannot be applied without significant modification.  

We now mention a few other questions arising from the work in this chapter.

\begin{Question}
For each connected subgraph $\Gamma$ of the Dynkin diagram, let $\g_\Gamma$ be the corresponding Levi subalgebra of $\g$, and $U_h(\g_\Gamma)$ be the corresponding Levi subalgebra of $U_h(\g)$. One can define $\xi''_\Gamma= Q_\Gamma^{-1/2} J_\Gamma T_{w_\Gamma}$, where $Q_\Gamma$ and $J_\Gamma$ act on a representation $V$ of $U_h(\g)$ via the obvious functor to $U_h(\g_\Gamma)$ representations, and $T_{w_\Gamma}$ is the braid group element corresponding to the longest word in the Weyl group of $\g_\Gamma$. These $\xi''_\Gamma$ are invertible, so they generate a group $\mathcal{W}$ acting on $U_h(\g)$ (which preserves only the algebra structure), and on representations of $U_h(\g)$. What group is this?

\end{Question}

We believe that answering this question would be an important step in understanding the relationship between the braid group and the cactus group. In the case $\g = \mbox{sl}_n$ we hope that $\mathcal{W}$ is closely related to the cactus group, where, if $\Gamma$ consists of nodes $s$ through $t-1$ of the Dynkin diagram, $\xi''_\Gamma$ corresponds to the generator $r_{[s,t]}$ of the cactus group (see \cite[Section 3]{cactus}). It cannot agree exactly since ${\xi''_\Gamma}^2 \neq \Id$ (only ${\xi''_\Gamma}^4 = \Id$). However, one may be able to use the methods of Section \ref{ocb} to modify the $\xi''$ and obtain exactly the cactus group. 

\begin{Question}
The commutor gives an action of the $n$-fruit cactus group $J_n$, and hence its group algebra, on tensor products $V \otimes \cdots \otimes V$ of $U_q(\g)$ modules.  Does this action factor through any quotient algebra in some special cases?  For example what about the case when $ \g = \mathrm{sl}_n $ and each $ V $ is the standard representation.  What about the corresponding action on tensor products of crystals?
\end{Question}

For the case of the braiding, the corresponding question has a nice answer in the above special case.  The action of the braid group factors through the Hecke algebra and we have the quantum Schur-Weyl duality.

The next question concerns an alternate approach to proving Conjecture \ref{CrystConjecture}. Kashiwara's $ * $-involution gives $ B_\infty $ an additional crystal structure, defined by $ f_i^* \cdot b := * \circ f_i \circ * (b) $.  Let $ B_\infty^i $ denote the crystal with vertex set $\mathbb{Z}_{\geq 0}$, where $ e_j, f_j $ act trivially for $ j \ne i $ and $ e_i, f_i $ act as they do on the usual $B_\infty $ for $ \mathfrak{sl}_2 $. Kashiwara \cite[Theorem 2.2.1]{Kashiwara:1993} showed that the map
\begin{equation}
\begin{aligned}
B_\infty &\rightarrow B_\infty \otimes B_\infty^i \\
b &\mapsto (e_i^*)^{\varepsilon_i(*b)}(b) \otimes \varepsilon_i(*b)
\end{aligned}
\end{equation}
is a morphism of crystals with respect to the usual crystal structures on each side.  We can think of this fact as an additional property of $ * $.

On the other hand, the crystal commutor $ \sigma^{HK} $ also has an additional property, which is called the cactus relation.  This relation states that if $ A, B, C $ are crystals, then 
\begin{equation} \sigma^{HK}_{A, C \otimes B} \circ (1 \otimes \sigma^{HK}_{B,C}) = \sigma^{HK}_{B \otimes A, C} \circ (\sigma^{HK}_{A, B} \otimes 1). \end{equation} (See \cite[Theorem 3]{cactus}).

\begin{Question}
Is there a relation between this additional property of Kashiwara's involution $ * $ and the cactus relation for the commutor $ \sigma^{HK}$?
\end{Question}

\chapter{Combinatorial models for $\asl_n$ crystals and applications} \label{abacus_chapter}

Much of this chapter has previously appeared in \cite{abacus}. 

\section{Crystal structures}  \label{crystal_def}

We start by defining two families of crystal structures for $\asl_n$. In the first, the vertices are partitions. In the second, the vertices are configurations of beads on an abacus. In each case the family is indexed by a positive integer $\ell$. We will refer to $\ell$ as the level of the crystal, since, using the results of this section, it is straightforward to see our level $\ell$ crystal structure decomposes into a union of crystals corresponding to level $\ell$ irreducible representations of $\asl_n$. The crystal graph of every irreducible $\asl_n$ representation appears as an easily identifiable subcrystal of the abacus model. Some irreducible representation do not appear using the partition model as we define it, although they do after a simple change of conventions (just shift the coloring in Figure \ref{ribboncrystal}). There is a slight subtlety in the case of $\asl_2$, since the proof of Theorem \ref{crystal_proof_th} fails. Theorem \ref{hipart} and the results from Section \ref{cylindric_partitions_and_abacus} do hold in this case, as we show in Section \ref{like_kashiwara} using the Kyoto path model. However, we cannot prove the full strength of Theorem \ref{crystal_proof_th} for $n =2$. 

\subsection{Crystal structure on partitions} \label{crystal_for_ribbons}
We now define level $\ell$ crystal operators for $\asl_n$, acting on the set of partitions (at least when $n \geq 3$). Figure \ref{ribboncrystal} illustrates the definition for $n=3$ and $\ell=4$. Color the boxes of a partition with $n$ colors $c_0, c_1, \ldots c_{n-1}$, where all boxes above position $k$ on the horizontal axis are colored $c_s$ for $\displaystyle s \equiv  \huge{ \lfloor} k /l \huge{\rfloor}$ modulo $n$. To act by $f_i$, place a $``("$ above the horizontal position $k$ if boxes in that position are colored $c_i$, and  you  can add an $\ell$-ribbon whose rightmost box is above $k$. Similarly, we put a $``)"$ above each position $k$ where boxes are colored $c_i$ and you can remove an $\ell$-ribbon whose rightmost box is above $k$. $f_i$ acts by adding an $\ell$-ribbon whose rightmost box is below the first uncanceled $``("$ from the left, if possible, and sending that partition to $0$ if there is no uncanceled $``("$. Similarly, $e_i$ acts by removing an $\ell$-ribbon whose rightmost box is below the first uncanceled $``)"$ from the right, if possible, and sending the partition to $0$ otherwise.

Theorem \ref{crystal_proof_th} below will imply that, for $n \geq 3$, these operators do in fact endow the set of partitions with an $\asl_n$ crystal structure. 

\subsection{Crystal structure on abacus configurations}
\label{abcryst}

As discussed in Section \ref{lotsastuff}, each partition corresponds to a certain configuration of beads on an abacus (abacus configuration). Thus the operators $e_i$ and $f_i$ also define operators on the set of abacus configurations coming from partitions. We know from Section \ref{lotsastuff} that adding (removing) an $\ell$-ribbon to a partition corresponds to moving one bead forward (backwards) one position in the corresponding $\ell$-strand abacus. This allows us to describe the operators $e_i$ and $f_i$ directly using the abacus. This new definition of operators $e_i$ and $f_i$  makes sense on all abacus configurations, regardless of whether or not they actually come from partitions.

The resulting operators $e_i$ and $f_i$ on abacus configurations can be defined as follows (see Figure \ref{abacuscrystal} for an example).
Color the gaps between the columns of beads with $n$ colors, putting  $c_0$ at the origin, and $c_{[ i \mod n ]}$ in the $i^{th}$ gap, counting left to right. We will include the colors in the diagrams by writing in the $c_i$'s below the corresponding gap. The operators $e_i$ and $f_i$ are then calculated as follows: Put a $``("$ every time a bead could move to the right across color $c_i$, and a $``)"$ every time a bead could move to the left across $c_i$. The brackets are ordered moving up each $c_i$ colored gap in turn from left to right. We group all the brackets corresponding to the same gap above that gap. $f_i$ moves the bead corresponding to the first uncanceled $``("$ from the left one place forward, if possible, and sends that element to 0 otherwise. Similarly, $e_i$ moves the bead corresponding the the first uncanceled $``)"$ from the right one space backwards, if possible, and sends the element to $0$ otherwise. 

We will now show that the operators $e_i$ and $f_i$ give the set of abacus configurations the structure of an $\asl_n$ crystal, when $n \geq 3$.  It follows immediately that our level $\ell$ operators on the set of partitions also gives us a crystal structure, since the map sending a partition to the corresponding abacus configuration preserves the operators $e_i$ and $f_i$. 

\begin{Theorem} \label{crystal_proof_th}
Fix $n \geq 3$ and $\ell \geq 1$. Define a colored directed graph $G$ as follows:  The vertices of $G$ are all configurations of beads on an $\ell$-strand abacus, which have finitely many empty positions to the left of the origin, and finitely many full positions to the right of the origin. There is a $c_i$-colored edge from $\psi$ to $\phi$ if and only if $f_i (\psi) = \phi$. Then each connected component of $G$ is the crystal graph of some integrable highest weight representation of $\asl_n$. 
\end{Theorem}

\begin{figure}
\begin{center}
\input{crystal_proof_diagram}
\end{center}
\caption[Showing that the abacus forms a crystal]{Consider an abacus configuration $\psi$ for $\asl_3$. Cut the abacus along the gaps colored $c_0$ and stack the pieces into a column as shown (this is the step that fails for $\asl_2$). Notice that we can calculate $f_1, f_2, e_1$ and $e_2$ on the column abacus using the same rule as for the row abacus, since we will look at the $c_1$ or $c_2$ gaps in exactly the same order. Each row of the column abacus, considered by itself, generates either a trivial crystal, or the crystal of one of the two fundamental representations of $\slt$. By Corollary \ref{mycrystaldef}, the rule for calculating $e_i$ and $f_i$ on the column abacus is exactly the same as the rule for taking the tensor product of all the non-trivial $\slt$ crystals. So, $\psi$ generates a component of the tensor product of a bunch of crystals of $\slt$ representations, which is itself the crystal of an $\slt$ representation.  \label{crystal_proof_diagram}}
\end{figure}
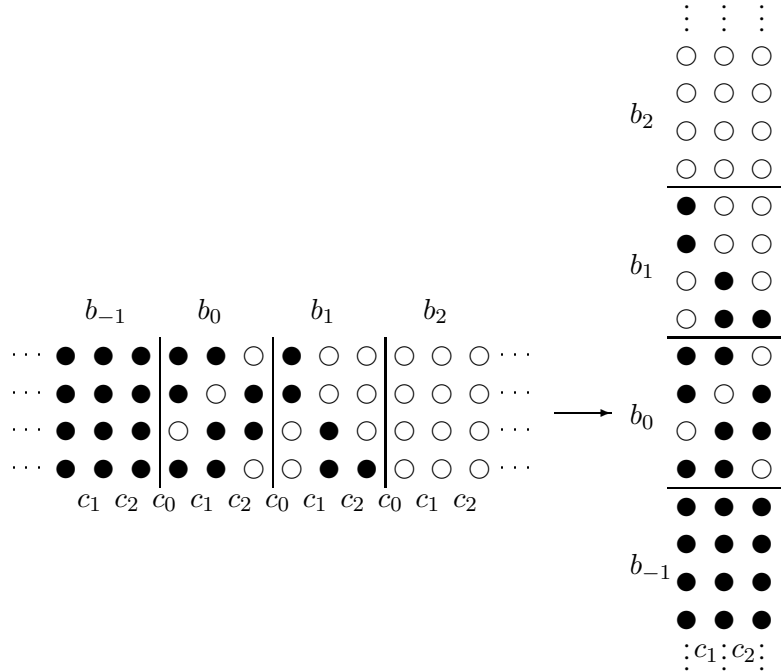

\begin{proof}
By \cite[Proposition 2.4.4]{KKMMNN1}, (see also \cite{Stembridge}), it is sufficient to show that, for each pair $0 \leq i < j<n$, each connected component of the graph obtained by only considering edges of color $c_i$ and $c_j$ is 
\begin{equation}
\begin{cases} \mbox{An} \hspace{2pt} \mbox{sl}_3 \hspace{2pt} \mbox{crystal if} \hspace{2pt}  |i-j|=1 \hspace{2pt} \mbox{mod}(n) \\  
\mbox{An sl}_2 \times \mbox{sl}_2 \hspace{2pt} \mbox{crystal otherwise.}
\end{cases}
\end{equation}
Choose some abacus configuration $\psi$, and some $0 \leq i < j \leq n-1$, and consider the connected component containing $\psi$ of the subgraph of $G$ obtained by only considering edges colored $c_i$ and $c_j$. If $|i-j| \neq 1$ $ \mbox{mod} (n)$, then $f_i$ and $f_j$ clearly commute (as do $e_i$ and $e_j$). 
Also, it is clear that if we only consider one color $c_i$, then $G$ is a disjoint union on finite directed lines. This is sufficient to show that the component containing $\psi$ is an $\mbox{sl}_2 \times \mbox{sl}_2$ crystal, as required.

Now consider the case when $|i-j|=1$ $ \mbox{mod} (n)$. The abacus model is clearly symmetric under shifting the colors, so, without loss of generality, we may assume $i=1$ and $j=2$. Also, we need only consider the abacus for $\asl_3$, since if there are columns of beads not bordering a $c_1$ or $c_2$ gap, they can be ignored without affecting $f_1$ or $f_2$. As shown in Figure \ref{crystal_proof_diagram}, each connected component of this crystal is the crystal of an integrable $\slt$ representation, as required. 
\end{proof}

\subsection{More on the abacus model, including the highest irreducible part} \label{highpart_section}
This section is a little technical. We introduce some definitions, then present some results about the structure of the crystals defined in Section \ref{abcryst}. Lemma \ref{structure_lemma} is the main result we need for the applications in Section \ref{cylindric_partitions_and_abacus}. We state and prove this lemma without using the fact that this is an $\asl_n$ crystal, so it holds even for $\asl_2$, when the proof of Theorem \ref{crystal_proof_th} fails. We finish the section by identifying a subcrystal which realizes the crystal graph for any irreducible integrable representation of $\asl_n$ of level $\ell >0$ (Theorem \ref{hipart}). We suggest the casual reader look at the statements of Lemma  \ref{structure_lemma} and Theorem \ref{hipart}, and refer to the definitions as needed; the details of the proofs can safely be skipped.

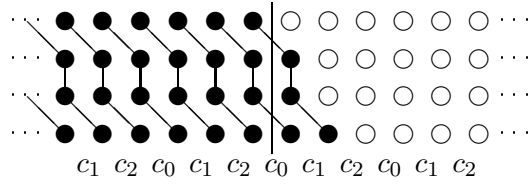
\begin{figure}
\begin{center}

\begin{picture}(20,3.5)
\input{compact1}
\end{picture}
\vspace{0.2in}
\caption[The highest weight element]{The highest weight element for an irreducible $\asl_3$ crystal of highest weight $\Lambda_0 + 2 \Lambda_1 + \Lambda_2$. As explained in Section \ref{like_kashiwara}, this corresponds to the ground-state path in the Kyoto path model. Each $4$-tuple of beads connected by a line corresponds to an element of $B_4$.  \label{compact1} }
\end{center}
\end{figure}

\begin{figure} \vspace{-0.4in}
\begin{center}
\begin{picture}(20,5)
\input{good_abacus}
\end{picture}
\vspace{0.2in}

\end{center}

\caption[A tight descending abacus]{A tight descending abacus. This configuration represents an element of the irreducible crystal generated by the compact abacus configuration shown in Figure \ref{compact1}. \label{good_abacus}}

\end{figure}
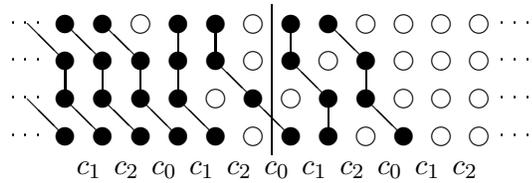

\begin{Definition}
Let $\psi$ be an abacus configuration. The {\bf compactification} of $\psi$, denoted $\psi_{(0)}$, is the configuration obtained by pushing all the (black) beads to the left, using only finitely many moves, and not changing the row of any bead. For example, Figure \ref{compact1} is the compactification of the configuration in Figure  \ref{good_abacus}.
\end{Definition}

\begin{Definition}
Let $r$ be a row of beads with finitely many negative positions empty and finitely many positive positions full. For each $i \in \bz_{>0}$, let $r^i$ denote the positions of the $i^{th}$ bead in $r$, counting from the right. 
\end{Definition}

\begin{Definition}
Let $r$ and $s$ be two rows of beads. We say $r \leq s$ if $r^i \leq s^i$ for all $i >0$. 
\end{Definition}

\begin{Definition} \label{moveup}
Let $\psi$ be an abacus configuration. Define a row of beads $\psi_i$ for any $i \in \bz$, by letting $\psi_i$ be the $i^{th}$ row of the abacus if $0 \leq i \leq \ell-1$ (counting up and starting with 0), and extending to the rest of $\bz$ using $\psi_{i + \ell}^j := \psi_i^j-n$. That is, $\psi_{i+\ell}$ is $\psi_i$, but shifted $n$ steps to the left.
\end{Definition}

\begin{Definition} \label{descendingdef}
We say a configuration of beads $\psi$ is descending if 
$\psi_0 \geq \psi_1 \geq \ldots \geq \psi_{\ell -1} \geq \psi_{\ell}$. Equivalently, $\psi$ is descending it if $\psi_i \geq \psi_{i+1}$ for all $i \in \bz$.
\end{Definition}

\begin{Definition} Let $\psi$ be a descending abacus configuration. Then the {\bf descending stands} are the sets  $\psi_\bigdot^k$ consisting of the $k^{th}$ black bead from the right on each row.
\end{Definition}

\begin{figure}
\begin{center}
\input{nottight}
\end{center}
\caption[Tightening the abacus]{A descending abacus that is not tight. We have joined the descending strands $\psi_\cdot^k$ by lines.  $\psi_\cdot^2$ can be shifted as shown. The positions modulo $n$ (here $3$) do not change, the configuration is still decreasing, but they are closer to $\psi_\bigdot^3$. Hence, this configuration is not tight. The operation shown here is $T_2$ as in Definition \ref{dtk}. \label{nottight}}
\end{figure}
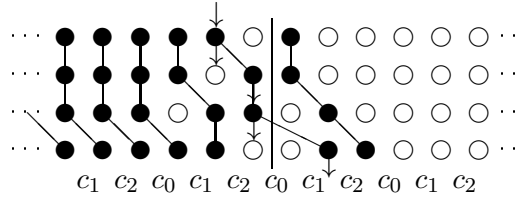

\begin{Comment}
Notice that a descending abacus configuration satisfies $\psi_0^k- n= \psi_\ell^k \leq \psi_{\ell-1}^k$ for all $k$. Thus, if we only draw $\ell$ rows of the abacus, each $\psi_\bigdot^k$ is at most $n+1$ columns wide (see Figure \ref{nottight}).
\end{Comment}

\begin{Definition} \label{dtk}
The tightening operator $T_k$ is the operator on abacus configuration $\psi$ which shifts the $k^{th}$ bead on each row down one row, if possible. Explicitly:
\begin{equation}
\begin{cases} 
\begin{cases}
T_k(\psi)_i^j = \psi_i^j \quad \mbox{if} \hspace{3pt} j \neq k \\
T_k(\psi)_i^k= \psi_{i+1}^k
\end{cases}
& \mbox{if} \hspace{3pt} \psi_{i+1}^k > \psi_i^{k+1} \hspace{2pt} \mbox{for all} \hspace{3pt} i \\
T_k (\psi) = 0 & \mbox{otherwise}
\end{cases}
\end{equation}
See Figure \ref{nottight} for an example. Similarly $T_k^*$ shifts the $k^{th}$ bead on each row up one row:  
\begin{equation}
\begin{cases} 
\begin{cases}
T_k^*(\psi)_i^j = \psi_i^j \quad \mbox{if} \hspace{3pt} j \neq k \\
T_k^*(\psi)_i^k= \psi_{i-1}^k
\end{cases}
& \mbox{if} \hspace{3pt} \psi_{i-1}^k < \psi_i^{k-1} \hspace{2pt} \mbox{for all} \hspace{3pt} i \\
T_k^* (\psi) = 0 & \mbox{otherwise}
\end{cases}
\end{equation}
\end{Definition}

\begin{Definition} \label{tightdef}
A descending configuration $\psi$ of beads is called {\bf tight} if, for each $k$, the beads $\psi^k_\bigdot$ are positioned as  far to the left as possible, subject to the set
$ \{ \psi^k_\bigdot \mod n \}$ being held fixed, and, for each $i$, $\psi_i^{k+1} < \psi^k_i$.  See Figure \ref{nottight}. This is equivalent to saying $T_k (\psi) =0$ for all $k$, since we only deal with descending configurations, so $T_k$ is the only way to move $\psi_\bigdot^k$ closer to $\psi_\bigdot^{k+1}$.
\end{Definition}

\begin{Lemma} \label{gglemma}
The set of descending abacus configurations (along with 0) is closed under the operators  $e_i$ and $f_i$ defined in Section \ref{abcryst}. Furthermore, the restriction of $e_i$ and $f_i$ to the set of descending abacus configurations can be calculated using the following rule:

For each $k$, let $S_i^k$ be the string of brackets $) \cdots ) ( \cdots ($ where the number of $``("$ is the number of beads in $\psi^k_\bigdot$  in position $i-1/2$ $\mbox{mod}(n)$, and the number of $``)"$ is the number of beads in  $\psi^k_\bigdot$  in position $i+1/2$ $\mbox{mod}(n)$. Let $S_i =  \cdots S_i^3 S_i^2 S_i^1$.

$\bullet$ If the first uncanceled $``)"$ from the right in $S_i$ comes from $S_i^k$, then $e_i$ moves a bead in $\psi_\bigdot^k$ one step to the left. The bead that moves is always the last bead of $\psi_\bigdot^k$ in position  $i+1/2$ $\mbox{mod}(n)$ that you encounter moving up the columns in turn from left to right. If there is no uncanceled $``)"$ in $S_i$, then $e_i$ sends that element to zero.

$\bullet$ If the first uncanceled $``("$ from the left in $S_i$ comes from $S_i^k$, then $f_i$ moves a bead in $\psi_\bigdot^k$ one step to the right. The bead that moves is always the first bead of $\psi_\bigdot^k$ in position  $i-1/2$ $\mbox{mod}(n)$ that you encounter moving up the columns in tun from left to right. If there is no uncanceled $``("$ in $S_i$, then $f_i$ sends that element to zero.
\end{Lemma}

\begin{proof}
Let $R_i$ be the string of brackets used to calculate  $e_i$ and $f_i$ in Section \ref{abcryst}. The descending condition implies:
\begin{enumerate}
\item \label{gg1} All the brackets in $R_i$ coming from beads in $\psi_\bigdot^k$ always come before all the brackets coming from beads in $\psi_\bigdot^{k-1}$.
\item \label{gg2} All the $``)"$ in $R_i$ coming from $\psi_\bigdot^k$ always come before all the $``("$ in $R_i$ coming from $\psi_\bigdot^k$
\end{enumerate}

These facts are both clear if you first do moves as in Figure \ref{cyclic_move} until the first bead of  $\psi_\bigdot^k$ in position $i+1/2$ $\mbox{mod}(n)$ is on the bottom row of the abacus. As argued in the caption, these moves commute with the actions of $e_i$ and $f_i$, and preserves the set of descending abacus configurations.

Let $R_i^k$ be the substring of $R_i$ consisting of brackets coming from  $\psi_\bigdot^k$. 
Then $S_i$ is obtained from $R_i$ by simply adding a string of canceling brackets $( \cdots ( ) \cdots )$ between each $R_i^k$ and $R_i^{k-1}$, where the number of $``("$ and $``)"$ is the number of pairs of touching beads, one in  $\psi_\bigdot^k$ and one in  $\psi_\bigdot^{k-1}$, that are in positions $1-1/2$ and $i+1/2$ $\mbox{mod}(n)$ respectively. Inserting canceling brackets does not change the first uncanceled $``("$.  Hence the first uncanceled $``("$ in $S_i$ will come from a bead in  $\psi_\bigdot^k$ if and only if the first uncanceled $``("$ in $R_i$ comes from a bead in  $\psi_\bigdot^k$. Therefore our new calculation of $f_i$ moves a bead in the right  $\psi_\bigdot^k$. It remains to show that the calculation of $f_i$ using $R_i$ always moves the first bead of  $\psi_\bigdot^k$ is position $i-1/2$ $\mbox{mod} (n)$. But this follows immediately from property (\ref{gg2}) above. The proof for $e_i$ is similar. Hence the new rule agrees with our definition of $e_i$ and $f_i$.

This new rule clearly preserves the set of descending abacus configurations.
\end{proof}

We are now ready to state and prove our main lemma concerning the structure of the operators $e_i$ and $f_i$ acting on abacus configurations:

\begin{figure}
\begin{picture}(40,4)
\input{cyclic_move}
\end{picture}
\vspace{0.02in}
\caption[Shifting the abacus]{Shifting the abacus. One can move the bottom row of the abacus to the top row, but shifted back by $n$. This will commute with the operations $f_i$, since when we create the string of brackets to calculate $e_i$ or $f_i$, we will still look at the $c_i$ gaps in exactly the same order. Also, it is clear from Definition \ref{descendingdef} that this preserves the set of descending abacus configurations. \label{cyclic_move}}
\end{figure}
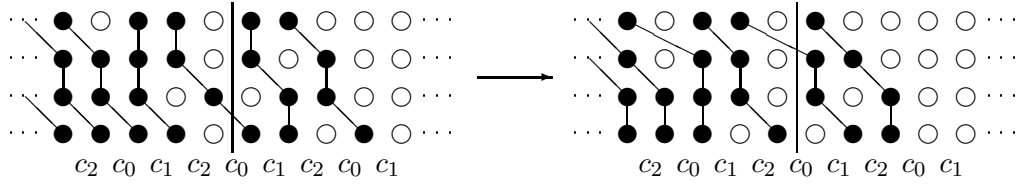

\begin{Lemma} \label{structure_lemma}
Fix $n \geq 2$, and $\ell \geq 1$. Consider the set of $\ell$-strand abacus configurations, colored with $c_0, \ldots c_{n-1}$ as shown in Figure \ref{abacuscrystal}. Let $G$ be the colored directed graph whose vertices are all descending abacus configurations $\psi$ with a given compactification $\psi_{(0)}$, and where there is an edge of color $c_i$ from $\psi$ to $\phi$ if and only if $f_i (\psi)=\phi$. Then:

\begin{enumerate}

\item \label{p2} The operators $e_i$ and $f_i$, restricted to the set of descending abacus configurations, commute with $ T_k$ and $T_k^*$.

\item \label{p3} The sources of $G$ (i.e. vertices that are not the end of any edge) are exactly those configurations that can be obtained from $\psi_{(0)}$ by a series of moves $T_k^*$ for various $k$.
 
\item \label{p3b} The set of tight descending abacus configurations is a connected component of $G$.

\item \label{p4} If we add an edge to $G$ connecting $\psi$ and $\phi$ whenever $T_k (\psi)=\phi$ for some $k$, then $G$ is connected. \end{enumerate}
\end{Lemma}

\begin{proof}

(\ref{p2}): Calculate $f_i$ as in Lemma \ref{gglemma}. Then $T_k$ clearly does not change the string of brackets $S_i$, and hence commute with $f_i$, as long as $T_k (\psi) \neq 0$ and $T_k \circ f_i (\psi) \neq 0$. The only potential problem is if $f_i \circ T_k (\psi) = 0$ but $T_k \circ f_i (\psi) \neq 0$ (or visa versa). 

Assume $f_i \circ T_k (\psi) = 0$ but $T_k \circ f_i (\psi) \neq 0$. Then $f_i$ must move the only bead $b$ of $\psi_\bigdot^k$ that hits $\psi_\bigdot^{k+1}$ when you apply $T_k$. $b$ must be the first bead of $\psi_\bigdot^k$ in position $i-1/2$ $\mbox{mod}(n)$ (moving up the columns and left to right). The last bead of $\psi_\bigdot^k$ in position $i+1/2$ $\mbox{mod}(n)$ must be on the the row below $b$. Since $b$ hits $\psi_\bigdot^{k+1}$ when you apply $T_k$, this must also be the last row of $\psi_\bigdot^{k+1}$ containing a bead in position $i-1/2$ $\mbox{mod}(n)$. By Lemma \ref{gglemma}, for $b$ to move, there must be at least as many beads of $\psi_\bigdot^k$ in position $i+1/2$ $\mbox{mod}(n)$ as beads in $\psi_\bigdot^{k+1}$ in position $i-1/2$ $\mbox{mod}(n)$. The descending condition then implies that the number of beads of $\psi_\bigdot^k$ in position $i+1/2$ $\mbox{mod}(n)$ must be equal to the number of beads of $\psi_\bigdot^{k+1}$ in position $i-1/2$ $\mbox{mod}(n)$. This is illustrated below:
\begin{center}
\begin{picture}(6,5)
\put(2.5,4){\vector(1,0){0.9}}
\put(2.5,3){\line(0,-1){0.7}}
\put(3.5,3){\line(0,-1){0.7}}
\put(2.5,1){\line(0,1){0.7}}
\put(3.5,1){\line(0,1){0.7}}
\put(2.5,4){\circle*{0.5}}
\put(2.5,3){\circle*{0.5}}
\put(2.5,1){\circle*{0.5}}
\put(3.5,3){\circle*{0.5}}
\put(3.5,1){\circle*{0.5}}
\put(2.5,5){\line(0,-1){1}}
\put(1.5,4){\line(1,-1){1}}
\put(2.5,4){\line(1,-1){1}}
\put(2.5,1){\line(1,-1){1}}
\put(3.5,1){\line(1,-1){1}}
\put(2.83,-0.3){$c_i$}
\end{picture}
\end{center}
But then $T_k \circ f_i (\psi)$ is in fact zero, since even after applying $f_i$ the bottom bead of $\psi_\bigdot^k$ in position $i+1/2$ $\mbox{mod}(n)$ cannot be shifted down without hitting $\psi_\bigdot^{k+1}.$ So the problem cannot in fact occur. The other case is similar, as are the cases involving $e_i$ or $T_k^*$.

(\ref{p3}): It is clear from part (\ref{p2}) that  applying operators $T_k^*$ to $\psi_{(0)}$ will always give a source (as long as the configuration is not sent to 0).  So,  let $\psi$ be a source, and we will show that $\psi$ can be tightened to $\psi_{(0)}$ by applying a series of operators $T_k$ for various $k$.  First, for each $k \geq 1$, define:
\begin{align}
 \varphi(\psi_\bigdot^k) := \sum_{j=0}^{\ell-1} \Lambda_i, \mbox{ where } i \equiv \psi_j^k + 1/2 \mbox{ mod} (n).  \\
  \varepsilon(\psi_\bigdot^k) := \sum_{j=0}^{\ell-1} \Lambda_i, \mbox{ where } i \equiv \psi_j^k - 1/2
\mbox{ mod} (n). 
\end{align}
We will use the notation $\varphi_i (\psi)$ (respectively $\varepsilon_i(\psi)$) to mean the coefficient of $\Lambda_i$ in $\varphi (\psi)$ (respectively $\varepsilon(\psi)$). 
For each $k \geq 1$, define $\psi |_k$ to be the abacus configuration obtained by removing the first $k-1$ black beads on each row, counting from the right. When we calculate $f_i$ as in Lemma \ref{gglemma}, all the brackets from $\psi|_k$ always come before the brackets from $\psi_\bigdot^m$, for any $m < k$. Hence, if $\psi$ is a source, then so is $\psi|_k$ for all $k$. We will prove the following statement for each $k \geq 1$:
\begin{align} \label{induct}
 \varphi(\psi|_{k+1}) = \varphi (\psi_\bigdot^{k+1}) = \varepsilon (\psi_\bigdot^{k}). 
\end{align}
$\psi$ can differ in only finitely many places from a compact configuration, so $\psi|_k$ is compact for large enough $k$, and (\ref{induct}) clearly holds for compact configurations. We proceed by induction, assuming (\ref{induct}) holds for some $k \geq 2$ and proving it still holds for $k-1$. 
Since $\varphi(\psi|_{k+1})= \varepsilon (\psi_\bigdot^{k}),$ the rule for calculating $f_i$ implies that $\varphi(\psi|_k)= \varphi(\psi_\bigdot^{k})$. Since $\psi|_{k-1}$ is still a source, we must have $\varepsilon_i(\psi_\bigdot^{k-1}) \leq \varphi_i(\psi_\bigdot^k)$ for each $i$. But $\sum_{i=0}^{n-1} \varepsilon_i(\psi_\bigdot^{k-1}) = \sum_{i=0}^{n-1} \varphi_i(\psi_\bigdot^{k}) = \ell$. Hence we see that  
$ \varphi (\psi_\bigdot^{k}) = \varepsilon (\psi_\bigdot^{k-1})$. So, (\ref{induct}) holds for $k-1$. 

By Definition \ref{tightdef} and the definitions of $\varphi(\psi)$ and $\varepsilon(\psi)$, (\ref{induct}) implies that we can tighten each $\psi_\bigdot^k$ until it is right next to $\psi_\bigdot^{k+1}$. Therefore, $\psi$ can be tightened to a compact configuration, which must be $\psi_{(0)}$ by the definition of $G$. Part (\ref{p3}) follows since $T_k(\psi)=\phi$ if and only if $T_k^*(\phi)=\psi$.

(\ref{p3b}): The graph $G$ is graded by $\bz_{\geq 0}$, where the degree of $\psi$ is the number of times you need to move one bead one step to the left to reach $\psi_{(0)}$. Every connected component to $G$ has at least one vertex $\psi_{\min}$ of minimal degree. $\psi_{\min}$ must be a source, since each $f_i$ is clearly degree 1. Hence each connected component of $G$ contains a source. By part (\ref{p3}), we can tighten any source to get $\psi_{(0)}$. Hence the set of tight descending abacus configurations in $G$ contains exactly one source, namely $\psi_{(0)}$. By Lemma \ref{gglemma} and part (\ref{p2}), the set of tight descending abacus configurations is closed under the operators $f_i$, so it is a complete connected component of $G$. 

(\ref{p4}): This follows immediately from part (\ref{p3}), and the observation from the proof of (\ref{p3b}) that every connected component of $G$ contains a source.
\end{proof}

\begin{Definition} \label{Lambdaa}
Let $\psi_{(0)}$ be a compact, descending abacus configuration. Define a dominant integral weight of $\asl_n$ by
\begin{equation}
\Lambda(\psi_{(0)}):= \sum_{i=0}^{n-1} m_i \Lambda_i,
\end{equation}
where $m_i$ is the number of $0 \leq j \leq \ell-1$ such that the last black bead of $\psi_{(0)j}$ is in position $i-1/2$ modulo n. Equivalently, 
\begin{equation}
m_i = \max \{ m : f_i^m (\psi_{(0)}) \neq 0 \}
\end{equation}
Note that $\Lambda(\psi_{(0)})$ uniquely determines $\psi_{(0)}$, up to a transformation of the form $\psi_{(0) i}^j \rightarrow \psi_{(0) i+m}^j$ for some $m \in \bz$ (see Definition \ref{moveup}). That is, up to a series of moves as in Figure \ref{cyclic_move}.
\end{Definition}

\begin{Theorem}  \label{hipart} Fix $n \geq 2$ and $\ell \geq 1$, and let $\psi_{(0)}$ be a compact, descending configuration of beads on an $\ell$-strand abacus colored with $c_0, \ldots c_{n-1}$ as shown in Figure \ref{abacuscrystal}. Let $B$ be the colored directed graph whose vertices are all tight, descending abacus configurations with compactification $\psi_{(0)}$, and there is a $c_i$ colored edge from $\psi$ to $\phi$ if $f_i(\psi)=\phi$. Then $B$ is a realization of the crystal graph for the $\asl_n$ representation $V_{\Lambda(\psi_{(0)})}$.
\end{Theorem}

\begin{proof}
Lemma \ref{structure_lemma} part (\ref{p3b}) shows that $B$ is a connected graph. For $n \geq 3$, Theorem \ref{crystal_proof_th} shows that this is in fact the crystal graph of an irreducible $\asl_n$ representation. For $n=2$, it is also the crystal of an irreducible $\asl_2$ representation, but we delay the proof until Section \ref{like_kashiwara}.

It is clear that $\psi_{(0)}$ is the highest weight element. It's weight is $ \sum_{i=0}^{n-1} m_i \Lambda_i$, where $m_i = \max \{ m : f_i^m (\psi_{(0)}) \neq 0 \}$. This is $\Lambda (\psi_{(0)})$ by Definition \ref{Lambdaa}. 
\end{proof}

This implies we can realize any integral irreducible highest weight crystal for $\asl_n$ using the abacus model by choosing the appropriate $\psi_{(0)}$ (see Figure \ref{compact1} for a typical highest weight element $\psi_{(0)}$, and Figure \ref{good_abacus} for a typical element of the highest irreducible part). However, we do not get every irreducible representation from the partition model. To see this, recall that the charge of a row of beads $r$ is the integer $c$ such that, when we push all the black  beads of $r$ to the left, the last black bead is in position $c-1/2$.  Define the charge of an abacus configuration $\psi$ to be the sum of the charges of the rows. The charge of a compact abacus configuration of highest weight $\Lambda$ is then well defined modulo $n$. Only those $\Lambda$ which correspond to compact configurations of charge zero $\mbox{mod} (n)$ can be realized using the crystal structure on partitions, as we define it. However, one can shift the colors in Figure \ref{ribboncrystal} to realize the other irreducible crystals.

\section{Relation to cylindric plane partitions} \label{cylindric_partitions_and_abacus}

We construct a bijection between  the set of descending abacus configurations with a given compactification and the set of cylindric plane partitions on a certain cylinder. This allows us to define a crystal structure on cylindric plane partitions, which is explicitly described in Section \ref{cpp_crystal}. The crystal we obtain is reducible, and carries an action of $\gli$ which commutes with $e_i$ and $f_i$. This allows us to show that the generating function for cylindric plane partitions with a fixed boundary is given by the principally graded character ($q$-character) of a certain $\agl_n$ representation. Borodin \cite{Borodin:2006} has previously calculated this generating function, and we show directly that the two results are equivalent. We also observe a form of rank-level duality, which comes from reflecting the cylindric plane partition in a vertical axis. 

\begin{figure} 
\begin{center}
\input{period_nl}
\end{center}
\caption[Labeled cylindric plane partition]{ A cylindric plane partition with period 9. Here $n=3$ and $\ell=6$, because there are 3 lines on the boundary sloping down and 6 sloping up in each period. $\pi_{ij}$ from Definition \ref{cppd} is the number at the intersection of diagonals $\pi_i$ and $c_j$ (where the empty squares are filled in with $0$). By Theorem \ref{weights} part (\ref{ccpp}), this cylindric plane partition corresponds to a basis element of $V_{2\Lambda_0+3 \Lambda_1 +\Lambda_2} \otimes F$ (a representation of $\asl_3 \oplus \gli$).  We can reflect the picture in a vertical axis, interchanging the $\pi$ diagonals with the $c$ diagonals. In this way the same cylindric plane partition corresponds to a basis element of $V_{\Lambda_0+\Lambda_1+\Lambda_4} \otimes F$ (a representation of $\asl_6 \oplus \gli$). \label{period_nl}}
\end{figure}
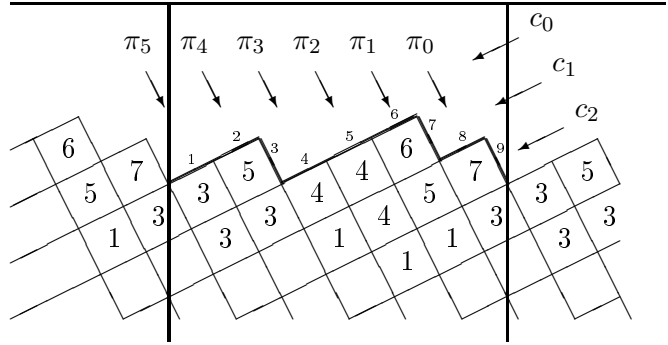

\subsection{A bijection between descending abacus configurations and cylindric plane partitions} \label{bpp}
In this section we define our bijection (Definition \ref{pidef}). Lemma \ref{piisgood} and Theorem \ref{acppb} show that this is a bijection, and relate the compactification of an abacus configuration to the boundary of the corresponding cylindric plane partition. Our definition of a cylindric plane partition (\ref{cppd}) is essentially the same as that used by Borodin in \cite{Borodin:2006}. See Figure \ref{period_nl} for an example. 
  
\begin{Definition} \label{cppd_body}
A cylindric plane partition of type $(n, \ell)$ is an array of non-negative integers $(\pi_{ij})$, defined for all sufficiently large $i, j \in \bz$, and satisfying:
\begin{enumerate}
\item \label{cppd1_body} If $\pi_{ij}$ is defined, then so is $\pi_{km}$ whenever $k \geq i$ and $m \geq j$.

\item \label{cppd2_body} $(\pi_{ij})$ is weakly decreasing in both $i$ and $j$. Furthermore, for all $i$, $\displaystyle \lim_{j \rightarrow \infty} \pi_{ij}=0$, and, for all $j$, $\displaystyle \lim_{i \rightarrow \infty} \pi_{ij}=0$. 

\item \label{cppd3_body} If $\pi_{ij}$ is defined, then $\pi_{ij} = \pi_{i+\ell, j-n}$.
\end{enumerate}

Notice that condition (\ref{cppd3_body}) implies that $\pi$ can be represented on a cylinder as in Figure \ref{period_nl}. The period of the cylinder, and the directions of the grid lines, are determined by $n$ and $\ell$. The information of which $\pi_{ij}$ are defined will be called the boundary of the cylinder.
\end{Definition}

\begin{Definition} \label{onemoredef}
A charged partition of charge $k \in \bz$ is a sequence of non-negative integers $(\lambda_{k}  \geq  \lambda_{k+1} \geq \lambda_{k+2} \geq \cdots )$ such that $\lambda_j=0$ for sufficiently large $j$. That is, it is a partition, but with the parts indexed starting at $k$. 
\end{Definition}

\begin{Definition} \label{pidef}
Let $\psi$ be an abacus configuration. For each $i \in \bz$, let $p_i$ be the integer such that the last black bead in the compactification of $\psi_i$ is in position $p_i-1/2$. For $j \geq p_i$, define $\pi_{ij}$ to be the number of black beads to the right of the $j-p_i+1$ st white bead of $\psi_i$, counting from the left.  Define $\pi(\psi):= (\pi_{ij}).$ As shown below (Lemma \ref{piisgood}), $\pi(\psi)$ is a cylindric plane partition. \end{Definition}

\begin{Comment}  \label{alongcomment}
Let $\pi=\pi(\psi)$. It should be clear that each $\pi_{i}:= (\pi_{ip_i}, \pi_{i,p_i+1}, \ldots)$ is a charged partition with charge $c_i$, as defined in Definition \ref{onemoredef}. Also, note that the boundary of $\pi (\psi)$ is determined by the charges of the $\ell$ rows of $\psi$, since $\pi_{ij}$ is well defined exactly when $j$ is at least the charge of $\psi_i$. In particular, the boundary will only depend on the compactification $\psi_{(0)}$ of $\psi$. \end{Comment}

\begin{Lemma} \label{piisgood}
For any descending abacus configuration $\psi$, $\pi(\psi)$ is a cylindric plane partition.
\end{Lemma}

\begin{proof}
Fix a descending abacus configurations $\psi$. 
We first show that part (\ref{cppd2_body}) of Definition \ref{cppd_body} holds for $\pi(\psi)$.
As in Section \ref{lotsastuff}, translate each row $\psi_i$ of $\psi$ into a diagram by going down and to the right one step for each black bead, and up and to the right one step for each white bead. Place this diagram at a height so that, far to the right, the diagrams of each $\psi_i$ lie along the same axis, as shown in Figure \ref{cyclic}. Since the $k^{th}$ black bead of $\psi_i$ is always to the right of the $k^{th}$ black bead of $\psi_{i+1}$, one can see that the sequence of diagrams is weakly decreasing by containment (i.e. the diagram of $\psi_{i+1}$ is never above the diagram for $\psi_i$). As in Definition \ref{pidef}, $\pi_{ij}$ is defined to be the number of black beads to the right of the $j-p_i+1$ st white bead of $\psi_i$. In the diagram, this is the distance of the diagram for $\psi_i$ in the $y$-direction from the interval $[j, j+1]$ on the $x$ axis (see Figure \ref{cyclic}). So $\pi_{ij}$ is weakly decreasing in $i$ because the sequence of diagrams is decreasing. Also $\pi_{ij}$ is weakly decreasing in $j$, since $\pi_i$ is a charged partition. So Definition \ref{cppd_body} part (\ref{cppd2_body}) holds. In fact, this argument also shows that Definition \ref{cppd_body} part (\ref{cppd1_body}) holds.

Now, by Definition \ref{moveup}, the row $\psi_{i+\ell}$ is just the row $\psi_{i},$ but shifted $n$ steps to the left. By definition \ref{onemoredef}, that means $\pi_{i+\ell}$ will be the same partition as $\pi_i$, but with the charge shifted by $n$. That is, $\pi_{ij}= \pi_{i+\ell, j-n}$, as required in part (\ref{cppd3_body}) of Definition \ref{cppd_body}.
\end{proof}

\begin{Definition} \label{Lambdab}
Let $\pi$ be a cylindric plane partition satisfying $\pi_{ij}= \pi_{i+\ell, j-n}$. Label the diagonals $c_0, c_1, c_2, \ldots$, as shown in Figure \ref{cpp}. Let
\begin{equation}
\Lambda(\pi):= \sum_{i=0}^{n-1} m_i \Lambda_i,
\end{equation}
where $m_i$ is the number of $0 \leq j \leq \ell-1$ such that the first entry of $\pi_j$ is in diagonal $c_k$ with $k \equiv i$  modulo n. 
\end{Definition}

\begin{Comment}
$\Lambda(\pi)$ is only determined by an (unlabeled) cylindric plane partition $\pi$ up to cyclically relabeling the fundamental weights $\Lambda_i$. This corresponds to the fact that there is a diagram automorphism of $\asl_n$ which cyclically permutes these weights.
\end{Comment}

\setlength{\unitlength}{0.17cm}
\begin{figure}[t]

\begin{center}
\begin{picture}(90,30)
\input{cyclic}
\end{picture}
\caption[The diagram of an abacus configuration]{The diagram of an abacus configuration.
This shows the diagram for $\psi_i$ ($0 \leq i \leq 5$), as described in the proof of Lemma \ref{piisgood}, for the example from Figure \ref{good_abacus}. \label{cyclic}}
\end{center}
\end{figure}
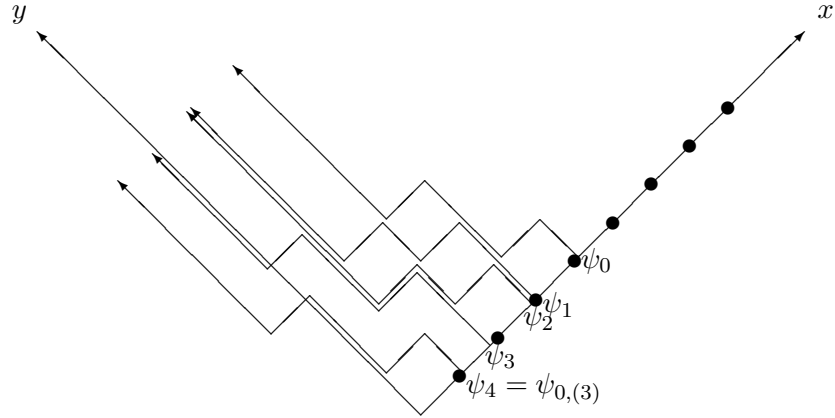

\setlength{\unitlength}{0.3cm}
\begin{figure}
\begin{center}
\input{ccp2}
\end{center}

\caption[The cylindric plane partition associated to an abacus configuration]{The cylindric plane partition associated to an abacus configuration. This figure shows $\pi (\psi)$ for the configuration $\psi$ shown in Figures \ref{good_abacus} and \ref{cyclic}. $\pi_{ij}$ is the intersection of the diagonals labeled $\pi_i$ and $c_j$. Notice that $\pi_4$ is just a shift of $\pi_0$. So we can cut on the lines shown, and wrap the diagram around a cylinder, to get a cylindric plane partition. We can also describe Definition \ref{Lambdab}:  $\Lambda (\pi) = \Lambda_0 + 2 \Lambda_1 + \Lambda_2$, since the coefficient of $\Lambda_j$ is the number of $0 \leq i \leq n-1$ for which $\pi_i$ has its first entry at a position $c_k$ with $k \equiv j$ modulo $n$. Note that the choice of labels $c_i$ can be shifted by $c_i \rightarrow c_{i+k}$. This corresponds to rotating the dynkin diagram of $\asl_n$. \label{cpp}}
\end{figure}
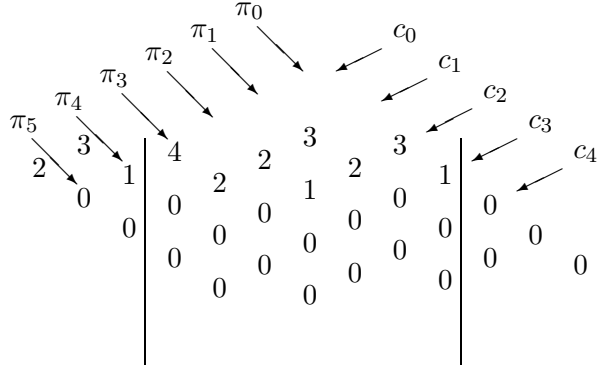
\setlength{\unitlength}{0.5cm}

\begin{Theorem} \label{acppb}
The  map $S: \psi \rightarrow \pi (\psi )$ is a bijection between the set of descending abacus configurations with compactification $\psi_{(0)}$, and the set of cylindric plane partitions with a given boundary. Furthermore, the boundary is determined by $\Lambda(\pi)= \Lambda(\psi_{(0)})$.
\end{Theorem}

\begin{proof}
Lemma \ref{piisgood} shows that $\pi(\psi)$ is always a cylindric plane partition, and, by Comment \ref{alongcomment}, the boundary of $\pi$ depends only on the compactification of $\psi$. Next, notice that, given a cylindric plane partition $\pi$ with the correct boundary, we can construct a diagram as in Figure \ref{cyclic} and then an abacus configuration, simply by reversing the procedure in Definition \ref{pidef}. This is clearly an inverse for $S$, so $S$ is a  bijection. It remains to check that $\Lambda(\pi)= \Lambda(\psi_{(0)})$. The charge of $\pi_i(\psi)$ is both the integer $k$ such that the last black bead of $\psi_i$ is in position $k-1/2$, and the integer $k$ such that the first entry of $\pi_i(\psi)$ is labeled with color $c_k$, so the identity follows from Definitions \ref{Lambdaa} and \ref{Lambdab}.
\end{proof}

\subsection{The crystal structure on cylindric plane partitions} \label{cpp_crystal}

We already have a crystal structure on abacus configurations (for $n \geq 3$) which, by Lemma \ref{gglemma}, preserves the set of descending configurations. Hence Theorem \ref{acppb} implies we have a crystal structure on the set of cylindric plane partitions. At the moment, we need to translate a cylindric plane partition into an abacus configuration to calculate $e_i$ and $f_i$. We now describe how to calculate these operators directly on the cylindric plane partitions. This section is largely independent of the rest of the paper, since later on we find it simpler to work with descending abacus configurations.

\begin{figure}
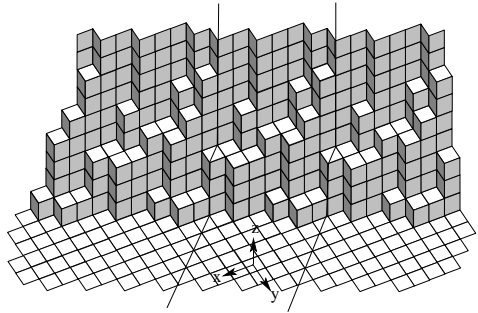


$$\mathfig{0.41}{3Dcpp.eps}$$

\caption[The three dimensional representation of a cylindric plane partition]{The three dimensional representation of a cylindric plane partition. This figure shows the cylindric plane partition from Figure \ref{period_nl}. The picture is periodic, with one period shown between the dark lines. The first layer of boxes should be colored as in Figure \ref{period_nl}, and higher levels according to the rule that color is constant along $(x_0, y_0+s, z_0+s)$ as $s$ varies. The planes $t(x,y,z)=C$ intersect the "floor" of the picture in a line which is horizontal in the projection shown, and are angled so that $(x,y,z)$ and $(x,y+1, z+1)$ are always on the same plane. For any $i \in I$, each such plane intersects the center of at most one $c_i$ colored box that could be added or removed  to/from $\pi$ (in each period). To calculate $f_i(\pi)$, place a $``("$  for each box in $A_i(\pi)$ and a  $``)"$  for each $c_i$ box in $R_i(\pi)$, ordered by $t$ calculated on the coordinates of the center of the box. $f_i$ adds the box corresponding to the first uncanceled $``("$, if there is one, and sends the element to $0$ otherwise.   \label{cpp3}} 
\end{figure}

It will be convenient to view a cylindric plane partition as a 3-dimensional picture, where $\pi_{ij}$ is the height of a pile of boxes placed at position $(i,j)$. See Figure \ref{cpp3}. Each box will be labeled by the coordinates of it's center in the $x$, $y$ and $z$ directions, as shown in Figure \ref{cpp3} (only relative positions matter, so the origin is arbitrary). Each box will also be labeled by a color $c_i$ for some residue $i$ modulo $n$. For the first layer of boxes, this will be determined as in Figure \ref{period_nl}. For higher boxes, one uses the rule that color is constant along lines of the form $\{ (x, y+k, z+k) : k \in \bz \}$. Note that due to the periodicity, $(x,y,z)$ labels the same box as $(x+\ell, y-n, z)$, so the coordinates are only well defined up to this type of transformation. 

For any cylindric plane partition $\pi$ and any residue $i$ modulo $n$, define the following two sets:

\begin{Definition}  $A_i(\pi) $ is the set of all $c_i$ colored boxes that can be added to $\pi$ so that each slice $\pi_i$ (see Figure \ref{cpp}) is still a valid partition.

 $R_i(\pi)$ is the set of all $c_i$ colored boxes that can be removed from $\pi$ so that each slice $\pi_i$ is still a valid partition.
 
\begin{Comment}Adding (removing) a box in $A_i(\pi)$ ($R_i(\pi)$) can result in something which is not a cylindric plane partition, since the slices $c_i$ (see Figure \ref{cpp}) may no longer be partitions.
 \end{Comment}
 
\end{Definition}

\begin{Definition} \label{defoft}
Define $t(x,y,z)=nx/\ell+y -z.$ Note that $t(x,y,z)= t(x+\ell, y-n, z)$, so $t$ is well defined as a function on boxes in a cylindric plane partition. For a box $n$ as in Figure \ref{cpp3}, define $t(n)$ to be $t$ calculated on the coordinates of the center of $n$.
\end{Definition}

\begin{Lemma} \label{total_order}
Fix a cylindric plane partition $\pi$. Let $n_1, n_2 \in A_i(\pi) \bigcup R_i(\pi)$. Then $t(n_1) \neq t(n_2)$ unless $n_1=n_2$.  
\end{Lemma}

\begin{proof}
Choose $n_0 \in A_i(\pi) \bigcup R_i(\pi)$, and consider the equation $ t(x,y,z)= t(n_0)$. It should be clear from Figure \ref{cpp3} that this plane intersects the center of at most one box in $A_i(\pi) \bigcup R_i(\pi)$ over each period. 
\end{proof}

We are now ready a define our crystal structure. Let $\pi$ be a cylindric plane partition. Define $S_i (\pi)$ to be the string of brackets formed by placing a $``("$  for every box $x$ in $A_i(\pi)$, and  $``)"$ for every box in $R_i(\pi)$. These are ordered with the brackets corresponding to $n_1$ coming before the brackets corresponding to $n_2$ if and only if $t(n_1) < t(n_2)$ (this is possible because of Lemma \ref{total_order}). Then $f_i(\pi)$ is obtained by adding the box corresponding to the first uncanceled $``("$, if there is one, and is $0$ otherwise. Similarly, $e_i(\pi)$ is obtained by removing the box corresponding to the first uncanceled $``)"$, if there is one, and is $0$ otherwise. It follows from that fact that the set of descending abacus configurations is closed under the crystal operators that both $f_i(\pi)$ and $e_i(\pi)$ will always be cylindric plane partitions. See Figure \ref{cpp3}.

To see this agrees with the crystal structure coming from descending abacus configurations, one should notice that adding (removing) a box $n$ colored $c_i$ with $t(n)=C$ always corresponds to moving a bead across the same gap on the abacus, which we will call $g_i(C)$. One must check that, for any $c_i$ colored boxes $n_1$ and $n_2$ in $A_i(\pi) \cup R_i(\pi)$, and $g_i(t(n_1))$ comes to the left of $g_i(t(n_2))$ when we calculate the crystal moves on the abacus if and only if $t(n_1)< t(n_2)$. This is straightforward.

\subsection{More structure on cylindric plane partitions}
In this section we give a bijection between descending abacus configurations with compactification $\psi_{(0)}$ and $B_\Lambda \times P$, where $B_\Lambda$ is the crystal associated to an irreducible $\asl_n$ representation $V_\Lambda$, and $P$ is the set of all partitions. It follows that there is also a bijection between cylindric plane partitions with a given boundary and $B_\Lambda \times P$. These bijections are useful because they preserve an appropriate notion of weight.

\begin{Definition}
Let $\psi$ be an abacus configuration. Notice that you can always transform $\psi$ to it's compactification by a series of moves, each of which moves one black bead exactly one step to the left. One can easily see that the number of such moves needed is well defined. This number is the weight of $\psi$, denoted $| \psi |$.
\end{Definition}

\begin{Definition}
If $\pi$ is a cylindric plane partition, the weight of $\pi$, denoted $|\pi |$, is the sum of the entries of $\pi$ over one period. 
\end{Definition}

\begin{Definition} \label{zdef} To each descending abacus configuration $\psi$, assign:
\begin{enumerate}

\item \label{z1} A tight descending abacus configuration $\gamma (\psi)$.

\item \label{z2} A partition $\lambda (\psi)$.

\item \label{z4} A canonical basis vector $v (\psi) \in V_\Lambda$.

\end{enumerate}

\noindent as follows:

(\ref{z1}): Moving from left to right, act by $T_k$ (see definition \ref{dtk}) until $\psi_\bigdot^k$ is tight with respect to $\psi_\bigdot^{k+1}$ for every $k$. The result is $\gamma (\psi)$. 

(\ref{z2}): Construct a row of beads by putting a black bead for each $\psi^i_\bigdot$, and, between the beads corresponding to $\psi^i_\bigdot$ and $\psi^{i+1}_\bigdot$, put a white bead for each time you can apply $T_i$ to $\psi^i_\bigdot$ and stay strictly to the right of $\psi^{i+1}_\bigdot$. Then use the correspondence between rows of beads and partitions to get $\lambda (\psi)$ (shifting to get charge 0 if necessary). 

(\ref{z4}): By Theorem \ref{hipart},  we know that the set of tight descending abacus configurations with compactification $\psi_{(0)}$ are the vertices of the crystal graph for $V_\Lambda$. These in turn correspond to canonical basis vectors in $V_\Lambda$ (see, for example, \cite{Hong&Kang:2000} for a discussion of canonical bases). Let $v(\psi)$ be the canonical basis element associated to $\gamma(\psi)$.
\end{Definition}

\begin{Comment} Since we have a bijection between cylindric plane partitions and abacus configurations, we can also associate this data to a cylindric plane partition $\pi$. In that case, we use the notation $\gamma(\pi), \lambda (\pi)$ and $v(\pi)$. 
\end{Comment}

Let $\gli$ be the Lie algebra of $(\bz+1/2) \times (\bz+1/2)$ matrices with finitely many non-zero entries. As in \cite{Kac:1990}, Chapter 14 (among others), the space $F$ spanned by all partitions is an irreducible $\gli$ module. The generators $E_{p, p+1}$ and $E_{p+1, p}$ of $\gli$ act on a partition as follows: Let $r \in F$, and consider the row of beads corresponding to $r$. 

$\bullet$  If position $p+1$ of $r$ is full and position $p$ is empty then $E_{p, p+1}$ moves the bead in position $p+1$ to position $p$. Otherwise $E_{p, p+1} (r)=0$.

$\bullet$  If position $p$ of $r$ is full and position $p+1$ is empty, then $E_{p+1, p}$ moves the bead in position $p$ to position $p+1$. Otherwise $E_{p+1, p} (r)=0$.

We then get a $\gli$ action on the space of abacus configurations, where $\gli$ acts on $\lambda (\psi)$ without affecting $\gamma(\psi)$. The following theorem will be our key tool in this section:

\begin{Theorem} \label{weights} 
\begin{enumerate}
\item \label{ccommutes} The $\asl_n$ crystal structure on the space of descending abacus configurations with compactification $\psi_{(0)}$ commutes with the $\gli$ action. Furthermore, $ \psi \rightarrow (\gamma (\psi),  \lambda (\psi))$ is a bijection with $B_\Lambda \times P$. Here $\Lambda = \Lambda(\psi_{(0)})$ (see Definition \ref{Lambdaa}), $B_\Lambda$ is the crystal graph of the $\asl_n$ representation $V_\Lambda$ realized according to Theorem \ref{hipart}, and $P$ is the set of all partitions.

\item \label{cabacus} $\{ v(\psi) \otimes \lambda (\psi) \},$ where $\psi$ ranges over all descending abacus configurations with a given compactification $\psi_{(0)}$, is a basis  $V_\Lambda \otimes F$. Here $F$ is the space spanned by all partitions and $\Lambda = \Lambda (\psi_{(0)})$ (see Definition \ref{Lambdaa}). Furthermore, the weight of $\psi$ is the sum of the principle graded weight of $v (\psi)$ and $n |\lambda(\psi)|$.

\item \label{ccpp} $\{ v(\pi) \otimes \lambda(\pi) \}$, where $\pi$ ranges over the set of all cylindric plane partitions $\pi$ with given boundary, is a basis for $V_\Lambda  \otimes F$. Here $F$ is the space of all partitions and $\Lambda = \Lambda (\pi)$ (see Definition \ref{Lambdab}). Furthermore, the weight of such a configuration is the sum of the principle graded weight of $v(\pi)$ and $n |\lambda(\pi)|$.

\end{enumerate}
\end{Theorem}

\begin{proof}
(\ref{ccommutes}): Notice that $\psi$ is uniquely defined by the pair $(\gamma (\psi), \lambda (\psi))$. Lemma \ref{structure_lemma} part (\ref{p2}) implies that $e_i$ and$f_i$ act on $\gamma (\psi)$, ignoring $\lambda (\psi)$. Also, by definition, $\gli$ acts only on $\lambda (\psi)$, ignoring $\gamma (\psi)$. Together, this implies that all $f_i$ commute with the $\gli$ action. Similarly, the $e_i$ commute with the $\gli$ action. Theorem \ref{hipart} says that the $\asl_n$ crystal generated by $\psi_{(0)}$ is $B_\Lambda$. Hence Lemma \ref{structure_lemma} part (\ref{p4}) implies that this map is a bijection to $B_\Lambda \otimes F$.

In order to prove (\ref{cabacus}) and (\ref{ccpp}), recall that the principally graded weight of a canonical basis vector $v \in V_\Lambda$ can be calculated from the crystal graph: set the weight of the highest weight element to be $0$, and let $f_i$ have degree $1$. This is a well defined grading on $B_\Lambda$ by standard results (see \cite{Hong&Kang:2000}). The principle grading on $V_\Lambda$ is obtained by letting $v(b)$ have the same weight as $b$. In particular, acting on $\psi$ by $f_i$ increases the principally graded weight of $v(\psi)$ by one.

(\ref{cabacus}): It follows from part (\ref{ccommutes}) that $\{ v(\psi) \otimes \lambda (\psi) \}$ is a basis for $V_\Lambda \otimes F$. It remains to show that the weights are correct. It is clear from the definitions that the operators $f_i$ add one to the weight of an abacus configuration, and $T_k^*$ adds $n$. As above, acting by $f_i$ increases the principally graded weight of $v(\psi)$ by $1$, without affecting $\lambda(\psi)$. Acting by $T_k^*$ adds one box to $\lambda(\psi)$ without affecting $v(\psi)$. The result follows.

(\ref{ccpp}): This follows from (\ref{cabacus}), once one notices that moving a bead one position to the right on the abacus (in such a way that it remains decreasing) always corresponds to adding $1$ to one entry of the corresponding cylindric plane partition.
\end{proof}

\subsection{The generating function for cylindric plane partitions} \label{rcpp}

In \cite{Borodin:2006}, Borodin studied the expected behavior of large random cylindric plane partitions under the distribution where the probability of a cylindric plane partition $\pi$ is proportional to $q^{|\pi |}$. In particular, he found the generating function for cylindric plane partitions on a given cylinder, which is the partition function for this system. In Corollary \ref{mypart_redo} below, we show that this generating function is the principally graded character of a certain representation of $\agl_n$.  We then directly show that this agrees with Borodin's formula.

Recall that, by Theorem \ref{weights}, we can  identify the set of cylindric plane partitions on a given cylinder with a basis for $V_\Lambda \otimes F$, and that this preserves an appropriate notion of weight.  Let $E_i$ and $F_i$ be the Chevalley generators of $\g$. Recall that the principle grading of a highest weight representation $V$ of $\g$ is the $\bz_{\geq 0}$ grading induced by putting the highest weight vector in degree $0$, letting all $F_i$ have degree 1, and letting all $E_i$ have degree $-1$.  Let $V_k$ denote the degree $k$ part of $V$ in this grading. The principally graded character or $q$-character of $V$ is the formal sum:
\begin{equation} \label{dimq_eq}
\dim_q (V) = \sum_{k \geq 0} \dim (V_k) q^k.
\end{equation}
Also, recall that $F$ is the space spanned by all partitions, and define the $q^n$ character of $F$ by:
\begin{equation} \label{dimq}
\dim_{q^n} (F) := \sum_{\mbox{partitions} \hspace{1.5pt} \lambda} q^{n | \lambda |}.
\end{equation}

\begin{Theorem} \label{mypart}
The generating function for cylindric plane partitions on a cylinder whose boundary satisfies $\Lambda (\pi) = \Lambda $ (see Definition \ref{Lambdab}) is given by:
\begin{align} 
 \label{dimqs} Z := \sum_{\begin{array}{c} \mbox{cylindric partitions} \quad \pi \\ \mbox{on a given cylinder} \end{array}}  q^{| \pi |} &= \dim_q (V_\Lambda)  \dim_{q^n} (F)
 \end{align}
\end{Theorem}

\begin{proof} This follows immediately from Theorem \ref{weights} part (\ref{ccpp}).
\end{proof}

A formula for $\dim_q(V_\Lambda)$ can be found in \cite[Proposition 10.10]{Kac:1990}. The $q^n$-character of $F$ is the sum over all partitions $\lambda$ of  $q^{n | \lambda |}$, and is well known (see for example \cite{Andrews:1976}, Theorem 1.1). Substituting these results into Equation (\ref{dimqs}) we obtain an equivalent formula for $Z$:
\begin{align} \label{myZ}
 Z=
\prod_{\alpha \in \Delta_+^\vee} \left( \frac{1- q^{ \langle \Lambda + \rho, \alpha \rangle}}{1-q^{\langle \rho, \alpha \rangle}} \right)^{\mbox{mult} (\alpha)}  \prod_{k =1}^\infty \frac{1}{1-q^{kn}}.
\end{align} 

Theorem \ref{mypart} can be simplified by wording it in terms of $\agl_n$ instead of $\asl_n$. In order to do this, we need to review some notation. 
 Our definitions agree with those in \cite{IFrenkel:1982}, although we have reworded some parts.

We use the following realization of $\hgl_n$:
\begin{equation}
\hgl_n = \mbox{gl}_n \otimes \bc[t, t^{-1}] + \bc c + \bc d.
\end{equation}

\noindent Here $c$ is central and $d$ acts by the derivation $t \frac{d}{dt}$. The bracket on the remaining part is defined by
\begin{equation}
[x \otimes t^p, y \otimes t^q] = [x,y] \otimes t^{p+q} + p \delta_{p+q,0} \langle x, y \rangle c,
\end{equation}
where $\langle x, y \rangle= \mbox{Tr} ( \mbox{ad}(x) \mbox{ad}(y)).$

Notice that $\asl_n$ is naturally contained in $\agl_n$. For a dominant integral weight $\Lambda$ of $\asl_n$, we say an irreducible representation $W_\Lambda$ of $\hgl_n$ has highest weight $\Lambda$ if there is some $w_\Lambda \in W_\Lambda$ such that:
\begin{enumerate}

\item{  $w_\Lambda$ is a highest weight vector of $W_\Lambda$ as an $\asl_n$ representation, and has weight $\Lambda$.}

\item{$\mbox{Id} \otimes t^n \cdot w_\Lambda=0$ for all $n >0$.}

\end{enumerate}

The principle grading on such a $W_\Lambda$ can be defined as follows: Put $w_\Lambda$ in degree zero. For each $0 \leq i \leq n-1$, let $E_i$ have degree $-1$ and $F_i$ have degree 1. Finally, let $c$ and $d$ have degree $0$, and $\mbox{Id} \otimes t^k$ have degree $-nk$. Denote by $W_{\Lambda,k}$ the degree $k$ weight space in this grading. The $q$-character of $W_\Lambda$ is the formal sum:
\begin{equation} \label{qqchar}
\dim_q(W_\Lambda)= \sum_{k=0}^\infty \dim(W_{\Lambda,k})q^k.
\end{equation}

\begin{Comment} \label{abcd}
It is a straightforward exercise to see that this definition agrees with that given in \cite{IFrenkel:1982}, and that the restriction of this grading to the $\asl_n$ representation $\asl_n \cdot w_\Lambda$ agrees with the principle grading defined in Equation \ref{dimq_eq}.
For any given $\Lambda$ there are many possibilities for $W_\Lambda$, since  $\mbox{Id} \otimes 1$ can act by any scalar. However, these all have the same $q$-character, so the $q$-character is well defined as a function of $\Lambda$.
\end{Comment}

\begin{Lemma} \label{aformula}
Fix a dominant integral weight $\Lambda$ for $\asl_n$. Let $V_\Lambda$ be the irreducible $\asl_n$ representation of highest weight $\Lambda$, and $W_\Lambda$ be an irreducible highest weight representation of $\hgl_n$ of highest weight $\Lambda$. Then:
\begin{equation} 
\dim_q(W_\Lambda)= \dim_q(V_\Lambda) \prod_{k=1}^\infty \frac{1}{1-q^k}.
\end{equation}
\end{Lemma}

\begin{proof}

Let $\mathcal{W}$ be the $\bc$ span of $c$ and $\mbox{Id} \otimes t^p$ for $p \in \bz$. Let $\mathcal{S}$ be the $\bc$ space of $c$, $d$, and $x \otimes t^p$ for $x \in \mbox{sl}_n$ and $p \in \bz$.
One can check that $\mathcal{W}$ forms a copy of the infinite dimensional Weyl algebra and $\mathcal{S}$ forms a copy of $\asl_n$. Furthermore, $[\mathcal{W},\mathcal{S}]=0$, $\hgl_n = \mathcal{W} + \mathcal{S}$, and $\mathcal{W} \cap \mathcal{S} = \bc c$. It follows that
\begin{equation} \label{stufff}
W_\Lambda= V \boxtimes F,
\end{equation}
where $V$ is an irreducible representation of $\asl_n$ and $F$ is an irreducible representation of $\mathcal{W}$. By Comment \ref{abcd}, we must have $V=V_\Lambda$. There are many irreducible highest weight representations of $\mathcal{W}$ (parameterized by the actions of $c$ and $\Id \otimes t^0$), but each has the same $q$-character with respect to the grading where $c$ has degree $0$, $\Id \otimes t^k$ has weight $-n k $, and the highest weight space is placed in degree $0$. This is given by:
\begin{equation}
\dim_q(F) = \prod_{k=1}^\infty \frac{1}{1-q^{nk}}.
\end{equation}

The Lemma follows by taking the $q$-character of each side of Equation (\ref{stufff}).
\end{proof}

We can now reword Theorem \ref{mypart} as follows:

\begin{Corollary} \label{mypart_redo}
Fix a dominant integral weight $\Lambda$ for $\asl_n$, and let $W_n$ be an irreducible representation of $\agl_n$ of highest weight $\Lambda$. 
The generating function for cylindric plane partitions on a cylinder whose boundary satisfies $\Lambda (\pi) = \Lambda $ (see Definition \ref{Lambdab}) is given by:
\begin{align} 
 \label{dimqs_redo} Z := \sum_{\begin{array}{c} \mbox{cylindric partitions} \quad \pi \\ \mbox{on a given cylinder} \end{array}}  q^{| \pi |} &= \dim_q (W_\Lambda)  
 \end{align}
\end{Corollary}

\begin{proof}
Follows immediately from Theorem \ref{mypart} and Lemma \ref{aformula}.
\end{proof}

We now compare this with Borodin's result. Following the conventions from \cite{Borodin:2006}, label the steps of the boundary of a cylindric plane partition with residues $i$ modulo $n + \ell$, as in Figure \ref{period_nl}. Define:

$\bullet$ $N= n+\ell$. This is the period of the boundary of the cylindric plane partition.

$\bullet$ For any $k \in \bz$, $k (N)$ is the smallest non-negative integer congruent to $k$ mod $N$.

$\bullet$ $\overline{1,N}$ is the set of integers modulo $N$.

$\displaystyle \bullet A[i] = \begin{cases}
1 \quad \mbox{if the boundary is sloping up and to the right on diagonal}  \quad i \\
0 \quad \mbox{otherwise}
\end{cases}$

$\displaystyle \bullet B[i] = \begin{cases}
1 \quad \mbox{if the boundary is sloping down and to the right on diagonal}  \quad i \\
0 \quad \mbox{otherwise}
\end{cases}$

Note: We have reversed the definitions of $A[i]$ and $B[i]$ from those used by Borodin. This corresponds to reflecting the cylindric plane partition about a vertical axis, so by symmetry does not change the partition function.

\begin{Theorem} (Borodin 2006) The generating function for cylindric plane partitions is given by:
\begin{equation} \label{BorodinZ} Z  := \hspace{-0.2in} \sum_{\begin{array}{c} \mbox{cylindric partitions} \quad \pi \\ \mbox{on a given cylinder} \end{array}} \hspace{-0.3in}  q^{| \pi |} = 
\prod_{k \geq 1} \frac{1}{1-q^{kN}} \prod_{\begin{array}{c} i \in \overline{1,N}: A[i]=1 \\ j \in \overline{1,N}: B[j]=1 \end{array} } \hspace{-0.2in} \frac{1}{1 - q^{(i-j)(N)+(k-1)N}}.
\end{equation}
\end{Theorem}

We will now directly show that Equations (\ref{myZ}) and (\ref{BorodinZ}) are equivalent. We already have a proof of Theorem \ref{mypart}, so some readers may wish to skip to Section \ref{rld}. 

We need the following fact about affine Lie algebras, which can be found in \cite[Chapter 14.2]{Kac:1990}. Let $\g$ be a finite dimensional simple complex Lie algebra. Let $\widehat{\g}$ be the associated untwisted affine algebra, and $\widehat{\g}'$ the corresponding derived algebra. The principle grading on $\widehat{\g}'$ is the grading determined by setting $\deg(H_i)=\deg(c)=0$, $\deg(E_i)=1$ and $\deg(F_i)= -1$, for each $i$. Let $r = \mbox{rank} (\g)$.

\begin{Lemma} (\cite{Kac:1990}, Chapter 14.2)
The dimension of the $k^{th}$ principally graded component of $\widehat{\g}'$ is $r+ s(k)$, where $s(k)$ is the number of exponents of $\g$ congruent to $k$ modulo the dual Coxeter number.
\end{Lemma}

\noindent In the case of $\asl_n$, this reduces to:

\begin{Corollary} \label{principle_dimension}
The dimension of the $k^{th}$ principally graded component of $\asl_n'$ is $n$, unless $k \equiv 0$ mod n, in which case the dimension is $n-1$.
\end{Corollary}

We start with Equation (\ref{myZ}) and manipulate it to reach Equation (\ref{BorodinZ}). First, notice that
$q^{\langle \rho, \alpha \rangle}$ is just $q$ to the principally graded weight of $\alpha$. Using this, and identifying $\Delta_+^\vee$ with $\Delta_+$ (since $\asl_n$ is self dual),  Corollary \ref{principle_dimension} allows us to simplify Equation (\ref{myZ}) as follows:

\vspace{-0.05in}

\begin{align}
Z & = \prod_{\alpha \in \Delta_+^\vee} (1- q^{ \langle \Lambda + \rho, \alpha \rangle} )^{mult ( \alpha )}
\prod_{k=1}^\infty \left( \frac{1}{1-q^k} \right)^{\begin{array}{l} \mbox{dim of the } \hspace{0.2pt} k^{th} \hspace{0.2pt} \mbox{principally} \\ \mbox{graded component of} \hspace{2pt} \asl_n' \end{array}}
\prod_{k=1}^{\infty} \frac{1}{1-q^{kn}} \label{aZeq}
\\
\label{twoparts} &=\prod_{\alpha \in \Delta_+^\vee} (1- q^{ \langle \Lambda + \rho, \alpha \rangle} )^{mult ( \alpha )}
\prod_{k=1}^\infty \left( \frac{1}{1-q^k} \right)^n.
\end{align}

We will deal with the second factor first. One can readily see that there are exactly $n$ residues $j \in\overline{1,N}$ such that $B[j]=1$ (see Figure \ref{period_nl}). Hence:

\begin{align}
\label{apart}
\prod_{k=1}^\infty \left( \frac{1}{1-q^k} \right)^n  & =  \prod_{\begin{array}{c}  j \in \{ 1, \ldots N \} , B[j]=1\\ i > j \end{array}}  \frac{1}{1-q^{i-j}}\\
\nonumber & = \prod_{k=1}^\infty \prod_{i,j \in \overline{1,N}: B[j]=1} \frac{1}{1 - q^{(i-j)(N)+(k-1)N}}
\end{align}

Now we consider the first factor of Equation (\ref{twoparts}) . Write:
\begin{equation}
\Lambda (\pi) = \sum_{i=0}^{n-1} m_i \Lambda_i.
\end{equation}
Notice that $m_i$ is the length of the $i^{th}$ downward sloping piece of the boundary of the cylindric plane partition, counting from the left over one period in Figure \ref{period_nl}, and starting with $0$ (see Definition \ref{Lambdab}). For each $i \in I$, set
$s_i := 1 + m_i.$
Let $\g_k (s)$ be the $k^{th}$ degree piece of $\asl_n'$ with the grading induced by letting $E_i$ have degree $s_i$, and $F_i$ have degree $-s_i$. This will be called the $s$-grading of $\asl_n'$. Then $q^{\langle \Lambda+\rho, \alpha \rangle}$ is just $q$ to the $s$-graded degree of $\alpha$. Again identifying $\Delta_+^\vee$ with $\Delta_+$, we see that: 

\begin{equation} \label{aroots}
 \prod_{\alpha \in \Delta_+^\vee} \left(1- q^{ \langle \Lambda + \rho, \alpha \rangle} \right)^{mult ( \alpha )} =
 \prod_{k =1}^\infty \left( 1-q^{\mbox{dim} ( \g_k (s) )} \right).
\end{equation}

\noindent $\mbox{dim} ( \g_k (s) )$ is the number of positive roots of $\asl_n$ of degree $k$ in the $s$ grading, counted with multiplicity. These are of two forms, which we must consider separately: 

\begin{enumerate} 

\item Real roots of the form $\beta+ k \delta$. Here $\beta$ is a root of $\mbox{sl}_n$, $k \geq 0$, and $k > 0$ if $\beta$ is negative. These each have multiplicity 1, and $s$ graded weight
$\langle \Lambda+ \rho, \beta \rangle + kN$ (since $\langle \Lambda + \rho, \delta \rangle = N$). 

\item Imaginary roots of the form $k \delta$. These each have multiplicity $n-1$, and $s$-graded weight $kN$.

\end{enumerate}

The imaginary roots contribute a factor of 
\begin{equation} \label{broots}
\prod_{k=1}^\infty \left( 1- q^{kN} \right)^{n-1} 
\end{equation}
to Equation (\ref{aroots}).

Next we find the contribution from the real roots. The positive roots of $\mbox{sl}_n$ are $\alpha_a + \alpha_{a+1} + \ldots + \alpha_b$ for all $1 \leq a <b \leq n-1$, and the negative roots are their negatives. A straightforward calculation shows that
\begin{equation}
\langle \rho+\Lambda , \sum_{i=a}^b \alpha_i \rangle=\left( b + \sum_{k=0}^{b-1} m_k \right)- \left(a + \sum_{k=0}^{a-1} m_k \right).\end{equation} 

One can see from Figure \ref{period_nl} that $B[a + m_0 + m_1 + \ldots + m_{a-1}]=1$ for all $a$. By a similar argument, if $\beta$ is a negative root of $\mbox{sl}_n$, then $\langle \rho+ \Lambda, \beta+\delta \rangle$ will be the difference of two integers that have $B[\cdot]=1$. It should be clear from Figure \ref{period_nl} that each pair $i,j \in \overline{1,N}$ with  $B[i]=B[j]=1$  corresponds to a root $\beta$ of $\mbox{sl}_n$ in one of these two ways. Therefore the real roots contribute a factor of
\begin{equation} \label{croots}
\prod_{k=1}^\infty \prod_{\begin{array}{c} i \in \overline{1,N}: B[i]=1 \\ j \in \overline{1,N}: B[j]=1  \\ i \neq j\end{array} } (1 - q^{(i-j)(N)+(k-1)N})
\end{equation}
to Equation (\ref{aroots}). Using Equations (\ref{broots}) and (\ref{croots}), Equations (\ref{aroots}) is equivalent to:
\begin{equation} \label{bpart}
 \prod_{\alpha \in \Delta_+^\vee} \left(1- q^{ \langle \Lambda + \rho, \alpha \rangle} \right)^{mult ( \alpha )} =
\prod_{k=1}^\infty  \left( 1- q^{kN} \right)^{n-1} \hspace{-0.3in}  \prod_{\begin{array}{c} i \in \overline{1,N}: B[i]=1 \\ j \in \overline{1,N}: B[j]=1  \\ i \neq j \end{array} }  \hspace{-0.4in}   \left(1 - q^{(i-j)(N)+(k-1)N}\right).
\end{equation}
Next, substitute Equations (\ref{apart}) and (\ref{bpart}) into Equation (\ref{aZeq}) to get:
\begin{align} 
Z & = \prod_{k=1}^\infty \left( 1- q^{kN} \right)^{n-1}  \hspace{-0.2in}  \prod_{\begin{array}{l} i, j \in \overline{1,N} \\ B[j]=1 \end{array}}  \hspace{-0.2in}  \frac{1}{1 - q^{(i-j)(N)+(k-1)N}}
\hspace{-0.8cm}  \prod_{\begin{array}{c} i \in \overline{1,N}: B[i]=1 \\ j \in \overline{1,N}: B[j]=1  \\ i \neq j \end{array} }  \hspace{-0.5in}   \left( 1 - q^{(i-j)(N)+(k-1)N} \right) \\
& \label{almost_done} = 
\prod_{k=1}^\infty \frac{1}{1-q^{kN}}
 \hspace{-0.2in}  \prod_{\begin{array}{l} i, j \in \overline{1,N} \\ B[j]=1 \end{array}}  \hspace{-0.2in}   \frac{1}{1 - q^{(i-j)(N)+(k-1)N}}
  \hspace{-0.1in}   \prod_{\begin{array}{c} i \in \overline{1,N}: B[i]=1 \\ j \in \overline{1,N}: B[j]=1  \end{array} } \hspace{-0.5in}  \left( 1 - q^{(i-j)(N)+(k-1)N} \right).
\end{align}
Equation (\ref{almost_done}) follows because there are exactly $n$ residues $i \in \overline{1,N}$ that have $B[i]=1$. This then simplifies to Equation (\ref{BorodinZ}), using the fact that $A[i]=1$ if and only if $B[i] \neq 1$. Hence we have directly shown that our equation for the partition function agrees with Borodin's.

\subsection{A rank-level duality result} \label{rld}
Fix a cylinder of type $(n, \ell)$. By Corollary \ref{mypart_redo}, we obtain the generating functions for cylindric plane partitions on this cylinder as the $q$-character of a certain $\agl_n$ representation of level $\ell$. We can reflect Figure \ref{period_nl} about a vertical axis to see that this generating function is also given by the $q$-character of a certain $\agl_\ell$ representation of level $n$. Thus we see that these two $q$-characters are equal. We now state this precisely, recovering a result of I. Frenkel \cite[Theorem 2.3]{IFrenkel:1982}.

For a residue $i$ modulo $n$, let  $\Lambda_i^{(n)}$ denote the $i^{th}$ fundamental weight of $\asl_n$.
For any level $\ell$ dominant integral weight $\Lambda = \sum_{i=0}^{n-1} c_i \Lambda_i^{(n)}$ of $\asl_n$, define a corresponding level $n$ dominant integral weight of $\asl_\ell$ by:
\begin{equation} \label{Lambda'}
\Lambda' := \sum_{i=0}^{n-1} \Lambda^{(\ell)}_{c_i+ c_{i+1}+ \cdots + c_{n-1}}.
\end{equation}
This should be understood using Figure \ref{period_nl}: If $\Lambda$ is defined as in the caption, then $\Lambda'$ is defined in the same way, but after first reflecting in the vertical axis, interchanging diagonal $\pi_k$ and diagonal $c_k$.

\begin{Theorem} (Frenkel 1982) \label{wldtheorem}
Let $\Lambda$ be a level $\ell$ dominant integral weight of $\asl_n$, with $n, \ell \geq 2$. Let $\Lambda'$ be the corresponding dominant integral weight of $\asl_\ell$ defined by Equation (\ref{Lambda'}). Let $W_\Lambda$ be an irreducible representation of $\agl_n$ with highest weight $\Lambda$, and $W_{\Lambda'}$ an irreducible representation of $\agl_\ell$ with highest weight $\Lambda'$. Then:
\begin{equation} \label{wldequation}
\dim_q (W_\Lambda)= \dim_q (W_{\Lambda'}). 
\end{equation}
\end{Theorem}

\begin{proof}[New proof of theorem \ref{wldtheorem}]
Let $S$ be the set of all cylindric plane partitions with a fixed boundary,  such that $\Lambda(\pi)= \Lambda$ for any $\pi \in S$.  By Corollary \ref{mypart_redo}, the generating function for $S$ is given by:
\begin{equation} \label{aa1}
Z := \sum_{\pi \in S}  q^{|\pi|} = \dim_q (W_\Lambda). 
\end{equation}
Now let $S'$ be the set of cylindric plane partitions with a fixed boundary,  such that $\Lambda(\pi)= \Lambda'$ for any $\pi \in S'$. Once again, by Corollary \ref{mypart_redo}, the generating function for $S'$ is given by:
\begin{equation} \label{aa2}
Z' := \sum_{\pi \in S'}  q^{|\pi|} = \dim_q (W_{\Lambda'}). 
\end{equation}
As discussed in the caption of Figure \ref{period_nl}, we can reflect an element of $S$ about a vertical axis to get an element of $S'$. This is a weight preserving bijection between $S$ and $S'$, so the two generating functions $Z$ and $Z'$ must be equal. 
\end{proof}
For completeness, we also include the corresponding result when $\ell=1$:

\begin{Theorem} \label{n2}
Let $\Lambda$ be a level $1$ dominant integral weight of $\asl_n$, with $n \geq 2$. Let $W_\Lambda$ be an irreducible representation of $\agl_n$ of highest weight $\Lambda$. Then:
\begin{equation}
\label{n22} \dim_q (W_\Lambda) = \prod_{k=0}^\infty \frac{1}{1-q^k}.
\end{equation}
\end{Theorem}

\begin{proof}
By Corollary \ref{mypart_redo}, 
\begin{equation}
\dim_q(W_\Lambda)= \dim_q(V_\Lambda) \prod_{k=0}^\infty \frac{1}{1-q^{nk}} .
\end{equation}
By Theorems \ref{mypart} and \ref{acppb}, this is $\sum q^{|\psi|}$, where the sum is over all $1$-strand (descending) abacus configurations $\psi$ with a given compactification $\psi_{(0)}$. Such a configuration is just one row of beads. We can simultaneously shift each $\psi$ until $\psi_{(0)}$ has its last black bead in position $-1/2$, without changing any $|\psi|$. As in Section \ref{lotsastuff}, there is a bijection between partitions and rows of beads whose compactification is this $\psi_{(0)}$, and one can easily see that this bijection preserves the weights. Hence the left side of Equation (\ref{n22}) is equal to $\sum_{\mbox{partitions } \lambda} q^{|\lambda|},$ which is well known to be given by the right side. \end{proof}

\section{Relation to the Kyoto path model} \label{like_kashiwara}

We now present an explicit bijection between the set of tight descending abacus configurations with a given compactification $\psi_{(0)}$, and the Kyoto path model for an integrable irreducible level $\ell$ representation of $\asl_n$, and show that the images of our operators $e_i$ and $f_i$ are the crystal operators $e_i$ and $f_i$. This proof does not assume that the abacus model is an $\asl_n$ crystal. Since the Kyoto path model is an $\asl_n$ crystal, we get a new proof of Theorem \ref{hipart}, which works for $\asl_2$ (where our previous proof failed). Since Lemma \ref{structure_lemma} parts (\ref{p2}) and (\ref{p4}) also hold in the case of $\asl_2$, the set of all descending abacus configurations is still an $\asl_2$ crystal. Hence the results in Section \ref{cylindric_partitions_and_abacus} still hold in the case $n=2$.

We refer the reader to (\cite{Hong&Kang:2000}, Chapter 10) for the details of the Kyoto path model. We will only use the path models corresponding  to one family $B_\ell$ of perfect crystals for $\asl_n$, one of each level $\ell > 0$. These can be found in (\cite{KKMMNN2} Theorem 1.2.2), where they are called $B^{1,\ell}$. Here $B_\ell$ consists of semi standard fillings of the partition $( \ell )$ with $0.5, 1.5, \ldots (2n-1)/2$. Below is an example for $n=3$ and $\ell=4$:

\begin{center}
\begin{picture}(4,1)
\put(0,0){\line(1,0){4}}
\put(0,1){\line(1,0){4}}
\put(0,0){\line(0,1){1}}
\put(1,0){\line(0,1){1}}
\put(2,0){\line(0,1){1}}
\put(3,0){\line(0,1){1}}
\put(4,0){\line(0,1){1}} 
\put(0.1, 0.32){\small{0.5}}
\put(1.1, 0.32){\small{1.5}}
\put(2.1, 0.32){\small{2.5}}
\put(3.1, 0.32){\small{2.5}}
\end{picture}

\end{center}

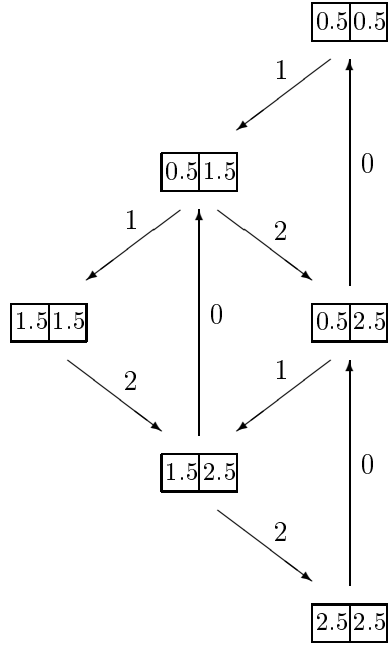
\begin{figure}
\begin{center}
\input{perfect_crystal}
\end{center}
\caption[A perfect crystal]{A perfect crystal of level 2 for $\asl_3$. The edges show the actions of $f_0, f_1,f_2$. The $e_i$ just follow the arrows backwards. If there is no $i$ arrow coming out of a vertex, $f_i$ sends that element to $0$. \label{perfect_crystal}}
\end{figure}

The operator $f_i$ changes the rightmost $i - 1/2$ to $i + 1/2$, if possible, and sends the element to 0 if there is no  $i - 1/2$. $f_0$ changes one $n - 1/2$ to $1/2$, if possible, then shuffles that element to the front. If there is no $n-1/2$ then $f_0$ sends that element to $0$. $e_i$ inverts $f_i$ if possible, and send the element to $0$ otherwise. This should be clear from Figure \ref{perfect_crystal}. Notice that we are using half integers instead of integers. The usual conventions are obtained by adding $1/2$ to everything.

Given a level $\ell$ perfect crystal $B_\ell$ for an affine algebra $\g$, and a level $\ell$ integral weight $\lambda$ of $\g$, 
Kashiwara et. al. define a {\it ground state path} of weight $\lambda$ to be a sequence 
\begin{equation} {\bf p}_\lambda =
\cdots \otimes b_{3}  \otimes  b_2 \otimes b_1 \in \cdots \otimes B_\ell \otimes B_\ell \otimes B_\ell \end{equation}
such that:
\begin{enumerate}
\item  $\varphi ( b_1) = \lambda$
\item For each $k \geq 1$, $\varepsilon ( b_k) = \varphi (b_{k+1})$.
\end{enumerate}
They show that, once $B_\ell$ is fixed, there is a unique ground state path of weight $\lambda$ for any level $\ell$ dominant integral weight $\lambda$. 

The irreducible crystal $P(\lambda)$ of highest weight $\lambda$ is the set of all all paths in $\cdots \otimes B_\ell \otimes B_\ell \otimes B_\ell$ which differ from the ground state path in only finitely many places. The crystal operators $e_i$ and $f_i$ are obtained from $B_\ell$ using the tensor product rule.
We use the notation $P_{n, \ell} (\Lambda)$ to denote the Kyoto path model for the $\asl_n$ crystal $B_\Lambda$ using the perfect crystal $B_\ell$ defined above (where $\ell$ is the level of $\Lambda$). 

Given an abacus configuration $\psi$, we define a sequence $\overline{\psi_\bigdot^k}$ of elements in $B_\ell$ (where $\ell$ is the number of rows of the abacus) by letting $\overline{\psi_\bigdot^k}$ be the unique semi-standard filling of $(\ell)$ with $\{ \psi_0^k(n), \ldots, \psi_{\ell-1}^k(n) \}$. Recall that $\psi_j^k$ is the position of the $k^{th}$ black bead from the right on row $j$ of $\psi$, and, as in Section \ref{rcpp}, $\psi_j^k(n)$ is the unique number $0 \leq \psi_j^k(n) < n$ that is congruent to $\psi_j^k$ modulo $n$.

\begin{Theorem} \label{Like_Kashiwara} Consider an $\ell$ row abacus colored with $c_0, \ldots c_{n-1}$ as in figure \ref{abacuscrystal}. Let $B$ be the set of tight, descending abacus configurations with compactification $\psi_{(0)}$. The map $J: B \rightarrow P_{n,\ell} (\Lambda (\psi_{(0)}))$ given by $\psi \rightarrow \cdots \otimes \overline{\psi_\bigdot^3} \otimes \overline{\psi_\bigdot^2} \otimes \overline{\psi_\bigdot^1}$ is a bijection. Furthermore, the images of our operators $e_i$ and $f_i$ are the crystal operators on $P_{n,\ell} (\Lambda(\psi_{(0)}))$.
\end{Theorem}

\begin{Comment}
This theorem really says $J$ is a crystal isomorphism. It gives a new proof of Theorem \ref{hipart}, which works even in the case $n=2$.
\end{Comment}

\begin{proof}

First note that the image of $\psi_{(0)}$ will satisfy the conditions to be a ground state path of weight $\Lambda(\psi_{(0)})$. For any $\psi \in B$, $\psi_j^k={\psi_{(0)}}_j^k$ for all but finitely many $0 \leq j \leq \ell-1$ and $k \geq 1$. Hence $J(\psi)$ will always lie in $P_{n,\ell} (\Lambda (\psi_{(0)}))$.

Use the string of brackets $S_i$ from Lemma \ref{gglemma} to calculate $f_i$ acting on $B$, and use the string of brackets $S_i'$ from Corollary \ref{mycrystaldef} to calculate $f_i$  acting on  $P_{n,\ell}(\Lambda(\psi)).$ One can easily see that $S_i'(J(\psi))=S_i(\psi)$. It follows from the two rules that the image of our operator $f_i$ is the crystal operator $f_i$ on $P_{n,\ell} (\Lambda(\psi_{(0)}))$. The proof for $e_i$ is similar.

Our operators $e_i$ and $f_i$ act transitively on $B$ by Lemma \ref{structure_lemma} part (\ref{p3b}), and the crystal operators $e_i$ and $f_i$ act transitively on  $P_{n,\ell} (\Lambda(\psi_{(0)}))$ since it is an irreducible crystal. $J$ preserves these operators, so it must be a bijection.
\end{proof}

\section{Future directions} \label{abacus_questions}

We finish with a short discussion of two questions which we feel are natural next steps for the work presented in this chapter.

\begin{Question}
Can one lift our crystal structures to get representations of $U_q(\asl_n)$? In particular, do cylindric plane partitions parameterize a basis for a representation of $U_q(\asl_n)$ in any natural way?
\end{Question}

In \cite{KMPY:1996}, Kashiwara, Miwa, Petersen and Yung study a $q$-deformed Fock space. This space has an action of $U_q(\asl_n)$, and a commuting action by the Bosons. In our picture, the space spanned by cylindric plane partitions has a $\asl_n$ crystal structure and a commuting action of $\gli$. Using the Boson-Fermion correspondence, one can realize $\gli$ inside a completion of the Bosonic algebra. We hope there is a natural action of $U_q(\asl_n)$ on the space spanned by cylindric plane partitions, coming from an embedding of this space into the $q$-deformed Fock space. 

This should have interesting consequences because of our rank-level duality: The space spanned by cylindric plane partitions on a given cylinder would simultaneously carry actions of $U_q(\asl_n)$, $U_q(\asl_\ell)$, and two copies of $\gli$. One copy of $\gli$ would commute with the $U_q(\asl_n)$ action and the other would commute with the $U_q(\asl_\ell)$ action. It would be interesting to see how the other actions interact. 

There are also several papers of Denis Uglov which may shed some light on this question (see \cite{Uglov:1999}, and references therein). He studies an action of $\asl_n$ on the space spanned by $\ell$-tuples of charged partitions. The crystal structure coming from his action is very similar to ours, and we believe it can be made to agree exactly by changing some conventions. He also studies commuting actions of $\asl_\ell$ and the Heisenberg algebra. It would be nice to fit these into our picture as well. However, his action does not appear to preserve the space spanned by descending abacus configurations, so it will not descend directly to the desired action on the space spanned by cylindric plane partitions.

\begin{Question}
Our Theorem \ref{mypart} gives a new way to calculate the partition function for a system of random cylindric plane partitions, using the representation theory of $\asl_n$. Can one use similar methods to do other calculations on such a statistical system? In particular, does this give a new way to calculate the correlation functions found by Borodin in \cite{Borodin:2006}? Alternatively, can one use results about random cylindric plane partitions to get interesting results about the representation theory of $\asl_n$?

\end{Question}

We have not pursued this line of inquiry in any depth, but it seems we have a new link between cylindric plane partitions and $\asl_n$ representation theory. One can certainly hope that results from one area will transfer in some way to the other.

\begin{Question}
What is the relationship between the combinatorics discussed in this chapter and the various geometric realizations of $\asl_n$ representations and their crystals? In particular, what is the relationship with Nakajima quiver varieties?
\end{Question}

There is a natural action of $\asl_n$ on the union of the cohomology rings of various Nakajima quiver varieties (see \cite{Nakajima:1994}). This is reducible, but by considering only the top non-zero cohomology of each space one obtains an irreducible representation. In fact, all irreducible integrable highest weight modules for $\asl_n$ occur in this way. 

After choosing a generic torus action on the Nakajima quiver varieties, one can index the cohomology by $\ell$-tuples of partitions. These in turn correspond to $\ell$-strand abacus configurations. However, depending on the choice of torus action, the set of $\ell$-tuples corresponding to the top non-zero cohomology changes. We believe this should give an explanation for an oddity of the abacus model: The same irreducible representation in fact occurs in many ways, not just as the set of tight descending abacus configurations.

\bibliographystyle{plain}
\bibliography{mybib}

\end{document}

%% file: tensor_figure.tex
\setlength{\unitlength}{0.6cm}
\begin{picture}(7.5,7)

\put(-1.5,1.8){$C$}
\put(4.7,6.8){$B$}

\put(0,0){\circle*{0.4}}
\put(0,2){\circle*{0.4}}
\put(0,4){\circle*{0.4}}

\put(2,0){\circle*{0.4}}
\put(2,2){\circle*{0.4}}
\put(2,4){\circle*{0.4}}
\put(2,6){\circle*{0.4}}

\put(4,0){\circle*{0.4}}
\put(4,2){\circle*{0.4}}
\put(4,4){\circle*{0.4}}
\put(4,6){\circle*{0.4}}

\put(6,0){\circle*{0.4}}
\put(6,2){\circle*{0.4}}
\put(6,4){\circle*{0.4}}
\put(6,6){\circle*{0.4}}

\put(8,0){\circle*{0.4}}
\put(8,2){\circle*{0.4}}
\put(8,4){\circle*{0.4}}
\put(8,6){\circle*{0.4}}

\put(0,2){\vector(0,-1){1.8}}
\put(0,4){\vector(0,-1){1.8}}

\put(2,6){\vector(1,0){1.8}}
\put(4,6){\vector(1,0){1.8}}
\put(6,6){\vector(1,0){1.8}}

\put(2,4){\vector(1,0){1.8}}
\put(4,4){\vector(1,0){1.8}}
\put(6,4){\vector(1,0){1.8}}

\put(2,2){\vector(1,0){1.8}}
\put(4,2){\vector(1,0){1.8}}

\put(2,0){\vector(1,0){1.8}}

\put(8,4){\vector(0,-1){1.8}}
\put(8,2){\vector(0,-1){1.8}}

\put(6,2){\vector(0,-1){1.8}}

\end{picture}
\setlength{\unitlength}{0.5cm}

%% file: ribboncrystal.tex
\begin{picture}(30,14.8)

\put(4.7,15.5){\Large{(}}
\put(14.7,15.5){\Large{(}}

\put(16.7,15.5){\Large{)}}
\put(17.7,15.5){\Large{)}}
\put(27.7,15.5){\Large{(}}
\put(29.7,15.5){\Large{(}}

\put(5, 15){\line(0,-1){4}}
\put(15, 15){\line(0,-1){2}}
\put(17, 15){\line(0,-1){2}}
\put(18, 15){\line(0,-1){3}}
\put(28, 15){\line(0,-1){1}}

\put(29,14){\line(-1,1){1}}
\put(28,15){\line(-1,-1){1}}
\put(27,14){\line(-1,1){1}}
\put(26,15){\line(-1,-1){2}}
\put(-0.3,13.9){$c_2$}
\put(0.7,12.9){$c_2$}
\put(1.7,11.9){$c_2$}
\put(2.7,10.9){$c_0$}
\put(3.7,9.9){$c_0$}
\put(4.7,8.9){$c_0$}
\put(5.7,7.9){$c_0$}
\put(6.7,6.9){$c_1$}
\put(7.7,5.9){$c_1$}
\put(8.7,4.9){$c_1$}
\put(9.7,3.9){$c_1$}
\put(10.7,2.9){$c_2$}
\put(11.7,1.9){$c_2$}
\put(12.7,0.9){$c_2$}
\put(13.7,-0.1){$c_2$}

\put(15.7,-0.1){$c_0$}
\put(16.7,0.9){$c_0$}
\put(17.7,1.9){$c_0$}
\put(18.7,2.9){$c_1$}
\put(19.7,3.9){$c_1$}
\put(20.7,4.9){$c_1$}
\put(21.7,5.9){$c_1$}
\put(22.7,6.9){$c_2$}
\put(23.7,7.9){$c_2$}
\put(24.7,8.9){$c_2$}
\put(25.7,9.9){$c_2$}
\put(26.7,10.9){$c_0$}
\put(27.7,11.9){$c_0$}
\put(28.7,12.9){$c_0$}
\put(29.7,13.9){$c_0$}

\put(5.7,9.9){$c_0$}
\put(6.7,8.9){$c_1$}
\put(7.7,7.9){$c_1$}
\put(8.7,6.9){$c_1$}
\put(9.7,5.9){$c_1$}
\put(10.7,4.9){$c_2$}
\put(11.7,3.9){$c_2$}
\put(12.7,2.9){$c_2$}
\put(13.7,1.9){$c_2$}
\put(14.7,0.9){$c_0$}
\put(15.7,1.9){$c_0$}
\put(16.7,2.9){$c_0$}
\put(17.7,3.9){$c_0$}
\put(18.7,4.9){$c_1$}
\put(19.7,5.9){$c_1$}
\put(20.7,6.9){$c_1$}
\put(21.7,7.9){$c_1$}
\put(22.7,8.9){$c_2$}
\put(23.7,9.9){$c_2$}
\put(24.7,10.9){$c_2$}
\put(25.7,11.9){$c_2$}

\put(7.7,9.9){$c_1$}
\put(8.7,8.9){$c_1$}
\put(9.7,7.9){$c_1$}
\put(10.7,6.9){$c_2$}
\put(11.7,5.9){$c_2$}
\put(12.7,4.9){$c_2$}
\put(13.7,3.9){$c_2$}
\put(14.7,2.9){$c_0$}
\put(15.7,3.9){$c_0$}
\put(16.7,4.9){$c_0$}
\put(17.7,5.9){$c_0$}
\put(18.7,6.9){$c_1$}
\put(19.7,7.9){$c_1$}
\put(20.7,8.9){$c_1$}
\put(21.7,9.9){$c_1$}
\put(22.7,10.9){$c_2$}
\put(23.7,11.9){$c_2$}

\put(8.7,10.9){$c_1$}
\put(9.7,9.9){$c_1$}
\put(10.7,8.9){$c_2$}
\put(11.7,7.9){$c_2$}
\put(12.7,6.9){$c_2$}
\put(13.7,5.9){$c_2$}
\put(14.7,4.9){$c_0$}
\put(15.7,5.9){$c_0$}
\put(16.7,6.9){$c_0$}
\put(17.7,7.9){$c_0$}
\put(18.7,8.9){$c_1$}
\put(19.7,9.9){$c_1$}
\put(20.7,10.9){$c_1$}
\put(21.7,11.9){$c_1$}

\put(12.7,8.9){$c_2$}
\put(13.7,7.9){$c_2$}
\put(14.7,6.9){$c_0$}
\put(15.7,7.9){$c_0$}
\put(16.7,8.9){$c_0$}
\put(17.7,9.9){$c_0$}
\put(18.7,10.9){$c_1$}
\put(19.7,11.9){$c_1$}

\put(13.7,9.9){$c_2$}
\put(14.7,8.9){$c_0$}
\put(15.7,9.9){$c_0$}
\put(16.7,10.9){$c_0$}
\put(17.7,11.9){$c_0$}

\put(14.7,10.9){$c_0$}
\put(15.7,11.9){$c_0$}
\put(16.7,12.9){$c_0$}

\put(15,0){\line(1,1){15}}
\put(14,1){\line(1,1){12}}
\put(13,2){\line(1,1){11}}
\put(12,3){\line(1,1){10}}
\put(11,4){\line(1,1){9}}
\put(10,5){\line(1,1){8}}
\put(9,6){\line(1,1){8}}
\put(8,7){\line(1,1){3}}
\put(7,8){\line(1,1){3}}
\put(6,9){\line(1,1){3}}
\put(5,10){\line(1,1){1}}

\put(15,0){\line(-1,1){15}}
\put(16,1){\line(-1,1){10}}
\put(17,2){\line(-1,1){9}}
\put(18,3){\line(-1,1){9}}
\put(19,4){\line(-1,1){6}}
\put(20,5){\line(-1,1){6}}
\put(21,6){\line(-1,1){6}}
\put(22,7){\line(-1,1){6}}
\put(23,8){\line(-1,1){6}}
\put(24,9){\line(-1,1){4}}
\put(25,10){\line(-1,1){3}}
\put(26,11){\line(-1,1){2}}
\put(27,12){\line(-1,1){1}}

\end{picture}

%% file: abacuscrystal.tex
\vspace{0.2in}

\begin{picture}(20,4)

\put(3,0){\ldots}
\put(3,1){\ldots}
\put(3,2){\ldots}
\put(3,3){\ldots}

\put(16,0){\ldots}
\put(16,1){\ldots}
\put(16,2){\ldots}
\put(16,3){\ldots}

\put(4.5,3){\circle*{0.5}}
\put(4.5,2){\circle*{0.5}}
\put(4.5,1){\circle*{0.5}}
\put(4.5,0){\circle*{0.5}}

\put(4.8, -1){$c_1$}

\put(5.5,3){\circle*{0.5}}
\put(5.5,2){\circle*{0.5}}
\put(5.5,1){\circle*{0.5}}
\put(5.5,0){\circle*{0.5}}

\put(5.8, -1){$c_2$}

\put(6.5,3){\circle*{0.5}}
\put(6.5,2){\circle*{0.5}}
\put(6.5,1){\circle*{0.5}}
\put(6.5,0){\circle*{0.5}}

\put(6.8, -1){$c_0$}

\put(6.8, 3.7){\big{(}}

\put(7.5,3){\circle*{0.5}}
\put(7.5,2){\circle{0.5}}
\put(7.5,1){\circle*{0.5}}
\put(7.5,0){\circle*{0.5}}

\put(7.8, -1){$c_1$}

\put(8.5,3){\circle*{0.5}}
\put(8.5,2){\circle*{0.5}}
\put(8.5,1){\circle{0.5}}
\put(8.5,0){\circle{0.5}}

\put(8.8, -1){$c_2$}

\put(9.5,3){\circle{0.5}}
\put(9.5,2){\circle{0.5}}
\put(9.5,1){\circle{0.5}}
\put(9.5,0){\circle*{0.5}}

\put(9.8, -1){$c_0$}

\put(9.3, 3.7){\big{(}}
\put(9.8, 3.7){\big{)}}
\put(10.3, 3.7){\big{)}}

\put(10,-0.5){\line(0,1){4}}

\put(10.5,3){\circle*{0.5}}
\put(10.5,2){\circle*{0.5}}
\put(10.5,1){\circle{0.5}}
\put(10.5,0){\circle{0.5}}

\put(10.8, -1){$c_1$}

\put(11.5,3){\circle*{0.5}}
\put(11.5,2){\circle{0.5}}
\put(11.5,1){\circle*{0.5}}
\put(11.5,0){\circle{0.5}}

\put(11.8, -1){$c_2$}

\put(12.5,3){\circle*{0.5}}
\put(12.5,2){\circle{0.5}}
\put(12.5,1){\circle*{0.5}}
\put(12.5,0){\circle{0.5}}

\put(12.5,0.8){$\longrightarrow$}

\put(12.8, -1){$c_0$}

\put(12.55, 3.7){\big{(}}
\put(13.05, 3.7){\big{(}}

\put(13.5,3){\circle{0.5}}
\put(13.5,2){\circle{0.5}}
\put(13.5,1){\circle{0.5}}
\put(13.5,0){\circle{0.5}}

\put(13.8, -1){$c_1$}

\put(14.5,3){\circle{0.5}}
\put(14.5,2){\circle{0.5}}
\put(14.5,1){\circle{0.5}}
\put(14.5,0){\circle{0.5}}

\put(14.8, -1){$c_2$}

\put(15.5,3){\circle{0.5}}
\put(15.5,2){\circle{0.5}}
\put(15.5,1){\circle{0.5}}
\put(15.5,0){\circle{0.5}}

\end{picture}

\vspace{0.1in}

%% file: compactification_intro.tex
\put(3,0){\ldots}
\put(3,1){\ldots}
\put(3,2){\ldots}
\put(3,3){\ldots}

\put(16,0){\ldots}
\put(16,1){\ldots}
\put(16,2){\ldots}
\put(16,3){\ldots}

\put(4.5,3){\circle*{0.5}}
\put(4.5,2){\circle*{0.5}}
\put(4.5,1){\circle*{0.5}}
\put(4.5,0){\circle*{0.5}}

\put(4.8, -1){$c_1$}

\put(5.5,3){\circle*{0.5}}
\put(5.5,2){\circle*{0.5}}
\put(5.5,1){\circle*{0.5}}
\put(5.5,0){\circle*{0.5}}

\put(5.8, -1){$c_2$}

\put(6.5,3){\circle*{0.5}}
\put(6.5,2){\circle*{0.5}}
\put(6.5,1){\circle*{0.5}}
\put(6.5,0){\circle*{0.5}}

\put(6.8, -1){$c_0$}

\put(7.5,3){\circle*{0.5}}
\put(7.5,2){\circle*{0.5}}
\put(7.5,1){\circle*{0.5}}
\put(7.5,0){\circle*{0.5}}

\put(7.8, -1){$c_1$}

\put(8.5,3){\circle*{0.5}}
\put(8.5,2){\circle*{0.5}}
\put(8.5,1){\circle*{0.5}}
\put(8.5,0){\circle*{0.5}}

\put(8.8, -1){$c_2$}

\put(9.5,3){\circle*{0.5}}
\put(9.5,2){\circle{0.5}}
\put(9.5,1){\circle*{0.5}}
\put(9.5,0){\circle{0.5}}

\put(9.8, -1){$c_0$}

\put(10,-0.5){\line(0,1){4}}

\put(10.5,3){\circle*{0.5}}
\put(10.5,2){\circle{0.5}}
\put(10.5,1){\circle{0.5}}
\put(10.5,0){\circle{0.5}}

\put(10.8, -1){$c_1$}

\put(11.5,3){\circle*{0.5}}
\put(11.5,2){\circle{0.5}}
\put(11.5,1){\circle{0.5}}
\put(11.5,0){\circle{0.5}}

\put(11.8, -1){$c_2$}

\put(12.5,3){\circle{0.5}}
\put(12.5,2){\circle{0.5}}
\put(12.5,1){\circle{0.5}}
\put(12.5,0){\circle{0.5}}

\put(12.8, -1){$c_0$}

\put(13.5,3){\circle{0.5}}
\put(13.5,2){\circle{0.5}}
\put(13.5,1){\circle{0.5}}
\put(13.5,0){\circle{0.5}}

\put(13.8, -1){$c_1$}

\put(14.5,3){\circle{0.5}}
\put(14.5,2){\circle{0.5}}
\put(14.5,1){\circle{0.5}}
\put(14.5,0){\circle{0.5}}

\put(14.8, -1){$c_2$}

\put(15.5,3){\circle{0.5}}
\put(15.5,2){\circle{0.5}}
\put(15.5,1){\circle{0.5}}
\put(15.5,0){\circle{0.5}}

%% file: descendingabacus.tex
\begin{picture}(20,5)

%\put(-1,-0.2){\small{(0)}}
%\put(-1,0.8){\small{(1)}}
%\put(-1,1.8){\small{(2)}}
%\put(-1,2.8){\small{(3)}}
%\put(-1,3.8){\small{(4)}}

\put(2.5,1){\line(1,-1){1}}

\put(1.5,4){\line(2,-1){2}}
\put(3.5,3){\line(0,-1){2}}
\put(3.5,1){\line(1,-1){1}}

\put(2.5,4){\line(2,-1){2}}
\put(4.5,3){\line(0,-1){2}}
\put(4.5,1){\line(1,-1){1}}

\put(4.5,4){\line(1,-1){1}}
\put(5.5,3){\line(0,-1){1}}
\put(5.5,2){\line(1,-1){2}}

\put(6.5,4){\line(0,-1){1}}
\put(6.5,3){\line(1,-1){1}}
\put(7.5,2){\line(1,-1){1}}
\put(8.5,1){\line(1,-1){1}}

\put(7.5,4){\line(1,-1){1}}
\put(8.5,3){\line(0,-1){1}}
\put(8.5,2){\line(2,-1){2}}
\put(10.5,1){\line(0,-1){1}}

\put(9.5,4){\line(2,-1){2}}
\put(11.5,3){\line(0,-1){1}}
\put(11.5,2){\line(1,-1){1}}
\put(12.5,1){\line(0,-1){1}}

\put(12.5,4){\line(1,-1){1}}
\put(13.5,3){\line(0,-1){1}}
\put(13.5,2){\line(1,-1){2}}

\put(2,0){\ldots}
\put(2,1){\ldots}
\put(2,2){\ldots}
\put(2,3){\ldots}
\put(-,4){\ldots}

\put(17,0){\ldots}
\put(17,1){\ldots}
\put(17,2){\ldots}
\put(17,3){\ldots}
\put(14,4){\ldots}

\put(0.5,4){\circle*{0.5}}
\put(0.5,4){\line(2,-1){2}}

\put(3.5,3){\circle*{0.5}}
\put(3.5,2){\circle*{0.5}}
\put(3.5,1){\circle*{0.5}}
\put(3.5,0){\circle*{0.5}}

\put(3.8, -1){$c_0$}

\put(4.5,3){\circle*{0.5}}
\put(4.5,2){\circle*{0.5}}
\put(4.5,1){\circle*{0.5}}
\put(4.5,0){\circle*{0.5}}
\put(1.5,4){\circle*{0.5}}

\put(4.8, -1){$c_1$}

\put(5.5,3){\circle*{0.5}}
\put(5.5,2){\circle*{0.5}}
\put(5.5,1){\circle{0.5}}
\put(5.5,0){\circle*{0.5}}
\put(2.5,4){\circle*{0.5}}

\put(5.8, -1){$c_2$}

\put(6.5,3){\circle*{0.5}}
\put(6.5,2){\circle{0.5}}
\put(6.5,1){\circle*{0.5}}
\put(6.5,0){\circle{0.5}}
\put(3.5,4){\circle{0.5}}

\put(6.8, -1){$c_0$}

\put(7.5,3){\circle{0.5}}
\put(7.5,2){\circle*{0.5}}
\put(7.5,1){\circle{0.5}}
\put(7.5,0){\circle*{0.5}}
\put(4.5,4){\circle*{0.5}}

\put(7.8, -1){$c_1$}

\put(8.5,3){\circle*{0.5}}
\put(8.5,2){\circle*{0.5}}
\put(8.5,1){\circle*{0.5}}
\put(8.5,0){\circle{0.5}}
\put(5.5,4){\circle{0.5}}

\put(8.8, -1){$c_2$}

\put(9.5,3){\circle{0.5}}
\put(9.5,2){\circle{0.5}}
\put(9.5,1){\circle{0.5}}
\put(9.5,0){\circle*{0.5}}
\put(6.5,4){\circle*{0.5}}

\put(9.8, -1){$c_0$}

\put(10,-0.5){\line(0,1){4}}

\put(10.5,3){\circle{0.5}}
\put(10.5,2){\circle{0.5}}
\put(10.5,1){\circle*{0.5}}
\put(10.5,0){\circle*{0.5}}
\put(7.5,4){\circle*{0.5}}

\put(10.8, -1){$c_1$}

\put(11.5,3){\circle*{0.5}}
\put(11.5,2){\circle*{0.5}}
\put(11.5,1){\circle{0.5}}
\put(11.5,0){\circle{0.5}}
\put(8.5,4){\circle{0.5}}

\put(11.8, -1){$c_2$}

\put(12.5,3){\circle{0.5}}
\put(12.5,2){\circle{0.5}}
\put(12.5,1){\circle*{0.5}}
\put(12.5,0){\circle*{0.5}}
\put(9.5,4){\circle*{0.5}}

\put(12.8, -1){$c_0$}

\put(13.5,3){\circle*{0.5}}
\put(13.5,2){\circle*{0.5}}
\put(13.5,1){\circle{0.5}}
\put(13.5,0){\circle{0.5}}
\put(10.5,4){\circle{0.5}}

\put(13.8, -1){$c_1$}

\put(14.5,3){\circle{0.5}}
\put(14.5,2){\circle{0.5}}
\put(14.5,1){\circle*{0.5}}
\put(14.5,0){\circle{0.5}}
\put(11.5,4){\circle{0.5}}

\put(14.8, -1){$c_2$}

\put(15.5,3){\circle{0.5}}
\put(15.5,2){\circle{0.5}}
\put(15.5,1){\circle{0.5}}
\put(15.5,0){\circle*{0.5}}
\put(12.5,4){\circle*{0.5}}

\put(15.8, -1){$c_0$}

\put(16.5,3){\circle{0.5}}
\put(16.5,2){\circle{0.5}}
\put(16.5,1){\circle{0.5}}
\put(16.5,0){\circle{0.5}}
\put(13.5,4){\circle{0.5}}

\end{picture}
\vspace{0.1in}

%% file: aplaincpp.tex
\setlength{\unitlength}{0.3cm}

\begin{center}
\begin{picture}(28,12)

\put(0,-17){ \begin{picture}(32,28)

\linethickness{0.3mm}
\put(-1,28){\line(1,0){29}}
\put(6,28){\line(0,-1){12}}
\put(21,28){\line(0,-1){12}}

\put(6,24.1){\line(2,1){4}}
\put(6,24.05){\line(2,1){4}}
\put(10.05,26){\line(1,-2){1}}
\put(10.1,26){\line(1,-2){1}}
\put(11,24.1){\line(2,1){6}}
\put(11,24.05){\line(2,1){6}}
\put(17.05,27){\line(1,-2){1}}
\put(17.1,27){\line(1,-2){1}}
\put(18,25.1){\line(2,1){2}}
\put(18,25.05){\line(2,1){2}}
\put(20.05,26){\line(1,-2){1}}
\put(20.1,26){\line(1,-2){1}}

\linethickness{0.15mm}

\put(1.25,25.2){6}
\put(2.25,23.2){5}
\put(3.25,21.2){1}

\put(4.25,24.2){7}
\put(5.25,22.2){3}

\put(7.25,23.2){3}
\put(8.25,21.2){3}

\put(9.25,24.2){5}
\put(10.23,22.2){3}

\put(12.25,23.2){4}
\put(13.25,21.2){1}

\put(14.25,24.2){4}
\put(15.25,22.2){4}
\put(16.25,20.2){1}

\put(16.25,25.2){6}
\put(17.25,23.2){5}
\put(18.25,21.2){1}

\put(19.25,24.2){7}
\put(20.25,22.2){3}

\put(22.25,23.2){3}
\put(23.25,21.2){3}

\put(24.25,24.2){5}
\put(25.23,22.2){3}

\put(0,26){\line(1,-2){4}}
\put(2,27){\line(1,-2){4.5}}
\put(5,26){\line(1,-2){4}}
\put(8,25){\line(1,-2){3.5}}
\put(10,26){\line(1,-2){4}}
\put(13,25){\line(1,-2){3.5}}
\put(15,26){\line(1,-2){4}}
\put(17,27){\line(1,-2){4.5}}
\put(20,26){\line(1,-2){4}}
\put(23,25){\line(1,-2){3.5}}
\put(25,26){\line(1,-2){1}}

\put(2,27){\line(-2,-1){3}}
\put(5,26){\line(-2,-1){6}}
\put(10,26){\line(-2,-1){11}}
\put(17,27){\line(-2,-1){18}}
\put(20,26){\line(-2,-1){16}}
\put(25,26){\line(-2,-1){16}}
\put(26,24){\line(-2,-1){12}}
\put(26,21.5){\line(-2,-1){7}}
\put(26,19){\line(-2,-1){2}}

 \end{picture}
}
\end{picture}

\end{center}

\setlength{\unitlength}{0.5cm}

%% file: partition.tex
\put(15,0){\line(1,1){15}}
\put(14,1){\line(1,1){12}}
\put(13,2){\line(1,1){11}}
\put(12,3){\line(1,1){10}}
\put(11,4){\line(1,1){9}}
\put(10,5){\line(1,1){7}}
\put(9,6){\line(1,1){5}}
\put(8,7){\line(1,1){3}}
\put(7,8){\line(1,1){3}}
\put(6,9){\line(1,1){3}}
\put(5,10){\line(1,1){1}}

\put(15,0){\line(-1,1){15}}
\put(16,1){\line(-1,1){10}}
\put(17,2){\line(-1,1){9}}
\put(18,3){\line(-1,1){9}}
\put(19,4){\line(-1,1){6}}
\put(20,5){\line(-1,1){6}}
\put(21,6){\line(-1,1){5}}
\put(22,7){\line(-1,1){5}}
\put(23,8){\line(-1,1){4}}
\put(24,9){\line(-1,1){4}}
\put(25,10){\line(-1,1){3}}
\put(26,11){\line(-1,1){2}}
\put(27,12){\line(-1,1){1}}

\put(2.5,0){\line(0,1){12.5}}
\put(3.5,0){\line(0,1){11.5}}
\put(4.5,0){\line(0,1){10.5}}
\put(6.5,0){\line(0,1){10.5}}
\put(9.5,0){\line(0,1){11.5}}
\put(10.5,0){\line(0,1){10.5}}
\put(11.5,0){\line(0,1){9.5}}
\put(14.5,0){\line(0,1){10.5}}
\put(17.5,0){\line(0,1){11.5}}
\put(20.5,0){\line(0,1){12.5}}
\put(22.5,0){\line(0,1){12.5}}
\put(24.5,0){\line(0,1){12.5}}
\put(26.5,0){\line(0,1){12.5}}

\put(2.5,0){\circle*{0.5}}
\put(3.5,0){\circle*{0.5}}
\put(4.5,0){\circle*{0.5}}
\put(6.5,0){\circle*{0.5}}
\put(9.5,0){\circle*{0.5}}
\put(10.5,0){\circle*{0.5}}
\put(11.5,0){\circle*{0.5}}
\put(14.5,0){\circle*{0.5}}
\put(17.5,0){\circle*{0.5}}
\put(20.5,0){\circle*{0.5}}
\put(22.5,0){\circle*{0.5}}
\put(24.5,0){\circle*{0.5}}
\put(26.5,0){\circle*{0.5}}

\put(1.5,0){\circle*{0.5}}
\put(0.5,0){\circle*{0.5}}
\put(-0.5,0){\circle*{0.5}}

\put(-2.2,-0.05){\ldots}

%% file: addribbon.tex
\begin{picture}(30,15)

\put(15,0){\line(1,1){15}}
\put(14,1){\line(1,1){12}}
\put(13,2){\line(1,1){11}}
\put(12,3){\line(1,1){10}}
\put(11,4){\line(1,1){9}}
\put(10,5){\line(1,1){7}}
\put(9,6){\line(1,1){5}}
\put(8,7){\line(1,1){3}}
\put(7,8){\line(1,1){3}}
\put(6,9){\line(1,1){3}}
\put(5,10){\line(1,1){1}}

\put(15,0){\line(-1,1){15}}
\put(16,1){\line(-1,1){10}}
\put(17,2){\line(-1,1){9}}
\put(18,3){\line(-1,1){9}}
\put(19,4){\line(-1,1){6}}
\put(20,5){\line(-1,1){6}}
\put(21,6){\line(-1,1){5}}
\put(22,7){\line(-1,1){5}}
\put(23,8){\line(-1,1){4}}
\put(24,9){\line(-1,1){4}}
\put(25,10){\line(-1,1){3}}
\put(26,11){\line(-1,1){2}}
\put(27,12){\line(-1,1){1}}

\put(2.5,0){\line(0,1){12.5}}
\put(3.5,0){\line(0,1){11.5}}
\put(4.5,0){\line(0,1){10.5}}
\put(6.5,0){\line(0,1){10.5}}
\put(9.5,0){\line(0,1){11.5}}
\put(10.5,0){\line(0,1){10.5}}
\put(11.5,0){\line(0,1){9.5}}
\put(18.5,0){\line(0,1){12.5}}
\put(17.5,0){\line(0,1){13.5}}
\put(20.5,0){\line(0,1){12.5}}
\put(22.5,0){\line(0,1){12.5}}
\put(24.5,0){\line(0,1){12.5}}
\put(26.5,0){\line(0,1){12.5}}

\put(2.5,0){\circle*{0.5}}
\put(3.5,0){\circle*{0.5}}
\put(4.5,0){\circle*{0.5}}
\put(6.5,0){\circle*{0.5}}
\put(9.5,0){\circle*{0.5}}
\put(10.5,0){\circle*{0.5}}
\put(11.5,0){\circle*{0.5}}
\put(14.5,0){\circle{0.5}}
\put(17.5,0){\circle*{0.5}}
\put(20.5,0){\circle*{0.5}}
\put(22.5,0){\circle*{0.5}}
\put(24.5,0){\circle*{0.5}}
\put(26.5,0){\circle*{0.5}}
\put(18.5,0){\circle*{0.5}}

\put(1.5,0){\circle*{0.5}}
\put(0.5,0){\circle*{0.5}}
\put(-0.5,0){\circle*{0.5}}

\put(-2.2,-0.05){\ldots}

\put(14,11){\line(1,1){3}}
\put(19,12){\line(-1,1){2}}

\end{picture}

%% file: bead_strand.tex
\put(5,0){
\begin{picture}(12,0.5)

\put(-5,0){\ldots}
\put(24,0){\ldots}

\put(-3.5,0){\circle*{0.5}}
\put(-2.5,0){\circle*{0.5}}
\put(-1.5,0){\circle*{0.5}}
\put(-0.5,0){\circle*{0.5}}
\put(0.5,0){\circle{0.5}}
\put(1.5,0){\circle*{0.5}}
\put(2.5,0){\circle{0.5}}
\put(3.5,0){\circle{0.5}}
\put(4.5,0){\circle*{0.5}}
\put(5.5,0){\circle*{0.5}}
\put(6.5,0){\circle*{0.5}}
\put(7.5,0){\circle{0.5}}
\put(8.5,0){\circle{0.5}}
\put(9.5,0){\circle*{0.5}}

\put(10,-0.5){\line(0,1){1}}

\put(10.5,0){\circle{0.5}}
\put(11.5,0){\circle{0.5}}
\put(12.5,0){\circle*{0.5}}
\put(13.5,0){\circle{0.5}}
\put(14.5,0){\circle{0.5}}
\put(15.5,0){\circle*{0.5}}
\put(16.5,0){\circle{0.5}}
\put(17.5,0){\circle*{0.5}}
\put(18.5,0){\circle{0.5}}
\put(19.5,0){\circle*{0.5}}
\put(20.5,0){\circle{0.5}}
\put(21.5,0){\circle*{0.5}}
\put(22.5,0){\circle{0.5}}
\put(23.5,0){\circle{0.5}}

\end{picture}

}

%% file: grouped_beads.tex
\put(5,0){
\begin{picture}(12,0.5)

\put(-5,0){\ldots}
\put(24,0){\ldots}

\put(-3.5,0){\circle*{0.5}}
\put(-2.5,0){\circle*{0.5}}

\put(-2,-0.3){\huge(}
\put(-2.4,-0.3){\huge)}

\put(-1.5,0){\circle*{0.5}}
\put(-0.5,0){\circle*{0.5}}
\put(0.5,0){\circle{0.5}}
\put(1.5,0){\circle*{0.5}}

\put(2,-0.3){\huge(}
\put(1.6,-0.3){\huge)}

\put(2.5,0){\circle{0.5}}
\put(3.5,0){\circle{0.5}}
\put(4.5,0){\circle*{0.5}}
\put(5.5,0){\circle*{0.5}}

\put(6,-0.3){\huge(}
\put(5.6,-0.3){\huge)}

\put(6.5,0){\circle*{0.5}}
\put(7.5,0){\circle{0.5}}
\put(8.5,0){\circle{0.5}}
\put(9.5,0){\circle*{0.5}}

\put(10,-0.3){\huge(}
\put(9.6,-0.3){\huge)}

\put(10,-0.5){\line(0,1){1}}

\put(10.5,0){\circle{0.5}}
\put(11.5,0){\circle{0.5}}
\put(12.5,0){\circle*{0.5}}
\put(13.5,0){\circle{0.5}}

\put(14,-0.3){\huge(}
\put(13.6,-0.3){\huge)}

\put(14.5,0){\circle{0.5}}
\put(15.5,0){\circle*{0.5}}
\put(16.5,0){\circle{0.5}}
\put(17.5,0){\circle*{0.5}}

\put(18,-0.3){\huge(}
\put(17.6,-0.3){\huge)}

\put(18.5,0){\circle{0.5}}
\put(19.5,0){\circle*{0.5}}
\put(20.5,0){\circle{0.5}}
\put(21.5,0){\circle*{0.5}}

\put(22,-0.3){\huge(}
\put(21.6,-0.3){\huge)}

\put(22.5,0){\circle{0.5}}
\put(23.5,0){\circle{0.5}}

\end{picture}

}

%% file: first_abacus.tex
\begin{picture}(20,3.7)

\put(0,0){
\begin{picture}(12,2)

\put(3,0){\ldots}
\put(3,1){\ldots}
\put(3,2){\ldots}
\put(3,3){\ldots}

\put(16,0){\ldots}
\put(16,1){\ldots}
\put(16,2){\ldots}
\put(16,3){\ldots}

\put(4.5,3){\circle*{0.5}}
\put(4.5,2){\circle*{0.5}}
\put(4.5,1){\circle*{0.5}}
\put(4.5,0){\circle*{0.5}}

\put(5.5,3){\circle*{0.5}}
\put(5.5,2){\circle*{0.5}}
\put(5.5,1){\circle*{0.5}}
\put(5.5,0){\circle*{0.5}}

\put(6.5,3){\circle*{0.5}}
\put(6.5,2){\circle*{0.5}}
\put(6.5,1){\circle*{0.5}}
\put(6.5,0){\circle*{0.5}}

\put(7.5,3){\circle*{0.5}}
\put(7.5,2){\circle{0.5}}
\put(7.5,1){\circle*{0.5}}
\put(7.5,0){\circle*{0.5}}

\put(8.5,3){\circle*{0.5}}
\put(8.5,2){\circle*{0.5}}
\put(8.5,1){\circle{0.5}}
\put(8.5,0){\circle{0.5}}

\put(9.5,3){\circle*{0.5}}
\put(9.5,2){\circle{0.5}}
\put(9.5,1){\circle{0.5}}
\put(9.5,0){\circle*{0.5}}

\put(10,-0.5){\line(0,1){4}}

\put(10.5,3){\circle{0.5}}
\put(10.5,2){\circle*{0.5}}
\put(10.5,1){\circle{0.5}}
\put(10.5,0){\circle{0.5}}

\put(11.5,3){\circle*{0.5}}
\put(11.5,2){\circle{0.5}}
\put(11.5,1){\circle*{0.5}}
\put(11.5,0){\circle{0.5}}

\put(12.5,3){\circle*{0.5}}
\put(12.5,2){\circle{0.5}}
\put(12.5,1){\circle*{0.5}}
\put(12.5,0){\circle{0.5}}

\put(13.5,3){\circle{0.5}}
\put(13.5,2){\circle{0.5}}
\put(13.5,1){\circle{0.5}}
\put(13.5,0){\circle{0.5}}

\put(14.5,3){\circle{0.5}}
\put(14.5,2){\circle{0.5}}
\put(14.5,1){\circle{0.5}}
\put(14.5,0){\circle{0.5}}

\put(15.5,3){\circle{0.5}}
\put(15.5,2){\circle{0.5}}
\put(15.5,1){\circle{0.5}}
\put(15.5,0){\circle{0.5}}

\end{picture}

}

\end{picture}

\vspace{0.1in}

%% file: second_abacus.tex
\begin{picture}(20,4)
\put(5,0){
\begin{picture}(12,2)

\put(3,0){\ldots}
\put(3,1){\ldots}
\put(3,2){\ldots}
\put(3,3){\ldots}

\put(16,0){\ldots}
\put(16,1){\ldots}
\put(16,2){\ldots}
\put(16,3){\ldots}

\put(4.5,3){\circle*{0.5}}
\put(4.5,2){\circle*{0.5}}
\put(4.5,1){\circle*{0.5}}
\put(4.5,0){\circle*{0.5}}

\put(5.5,3){\circle*{0.5}}
\put(5.5,2){\circle*{0.5}}
\put(5.5,1){\circle*{0.5}}
\put(5.5,0){\circle*{0.5}}

\put(6.5,3){\circle*{0.5}}
\put(6.5,2){\circle*{0.5}}
\put(6.5,1){\circle*{0.5}}
\put(6.5,0){\circle*{0.5}}

\put(7.5,3){\circle*{0.5}}
\put(7.5,2){\circle{0.5}}
\put(7.5,1){\circle*{0.5}}
\put(7.5,0){\circle*{0.5}}

\put(8.5,3){\circle*{0.5}}
\put(8.5,2){\circle*{0.5}}
\put(8.5,1){\circle{0.5}}
\put(8.5,0){\circle{0.5}}

\put(9.5,3){\circle{0.5}}
\put(9.5,2){\circle{0.5}}
\put(9.5,1){\circle{0.5}}
\put(9.5,0){\circle*{0.5}}

\put(10,-0.5){\line(0,1){4}}

\put(10.5,3){\circle*{0.5}}
\put(10.5,2){\circle*{0.5}}
\put(10.5,1){\circle{0.5}}
\put(10.5,0){\circle{0.5}}

\put(11.5,3){\circle*{0.5}}
\put(11.5,2){\circle{0.5}}
\put(11.5,1){\circle*{0.5}}
\put(11.5,0){\circle{0.5}}

\put(12.5,3){\circle*{0.5}}
\put(12.5,2){\circle{0.5}}
\put(12.5,1){\circle*{0.5}}
\put(12.5,0){\circle{0.5}}

\put(13.5,3){\circle{0.5}}
\put(13.5,2){\circle{0.5}}
\put(13.5,1){\circle{0.5}}
\put(13.5,0){\circle{0.5}}

\put(14.5,3){\circle{0.5}}
\put(14.5,2){\circle{0.5}}
\put(14.5,1){\circle{0.5}}
\put(14.5,0){\circle{0.5}}

\put(15.5,3){\circle{0.5}}
\put(15.5,2){\circle{0.5}}
\put(15.5,1){\circle{0.5}}
\put(15.5,0){\circle{0.5}}

\end{picture}

}
\end{picture}

\vspace{0.1in}

%% file: crystal_proof_diagram.tex
\begin{picture}(30,19)

\put(0,-1){\begin{picture}(30,20)

\put(17.5,8.5){\vector(1,0){1.5}}

\put(20.9,18.5){$\cdot$}
\put(20.9,18.8){$\cdot$}
\put(20.9,19.1){$\cdot$}

\put(21.9,18.5){$\cdot$}
\put(21.9,18.8){$\cdot$}
\put(21.9,19.1){$\cdot$}

\put(22.9,18.5){$\cdot$}
\put(22.9,18.8){$\cdot$}
\put(22.9,19.1){$\cdot$}

\put(20.9,1.5){$\cdot$}
\put(20.9,1.8){$\cdot$}
\put(20.9,2.1){$\cdot$}

\put(21.9,1.5){$\cdot$}
\put(21.9,1.8){$\cdot$}
\put(21.9,2.1){$\cdot$}

\put(22.9,1.5){$\cdot$}
\put(22.9,1.8){$\cdot$}
\put(22.9,2.1){$\cdot$}

\put(21, 18){\circle{0.5}}
\put(22, 18){\circle{0.5}}
\put(23, 18){\circle{0.5}}
\put(21, 17){\circle{0.5}}
\put(22, 17){\circle{0.5}}
\put(23, 17){\circle{0.5}}
\put(21, 16){\circle{0.5}}
\put(22, 16){\circle{0.5}}
\put(23, 16){\circle{0.5}}
\put(21, 15){\circle{0.5}}
\put(22, 15){\circle{0.5}}
\put(23, 15){\circle{0.5}}

\put(20.5,14.5){\line(1,0){3}}

\put(21, 14){\circle*{0.5}}
\put(22, 14){\circle{0.5}}
\put(23, 14){\circle{0.5}}
\put(21, 13){\circle*{0.5}}
\put(22, 13){\circle{0.5}}
\put(23, 13){\circle{0.5}}
\put(21, 12){\circle{0.5}}
\put(22, 12){\circle*{0.5}}
\put(23, 12){\circle{0.5}}
\put(21, 11){\circle{0.5}}
\put(22, 11){\circle*{0.5}}
\put(23, 11){\circle*{0.5}}

\put(20.5,10.5){\line(1,0){3}}

\put(21, 10){\circle*{0.5}}
\put(22, 10){\circle*{0.5}}
\put(23, 10){\circle{0.5}}
\put(21, 9){\circle*{0.5}}
\put(22, 9){\circle{0.5}}
\put(23, 9){\circle*{0.5}}
\put(21, 8){\circle{0.5}}
\put(22, 8){\circle*{0.5}}
\put(23, 8){\circle*{0.5}}
\put(21, 7){\circle*{0.5}}
\put(22, 7){\circle*{0.5}}
\put(23, 7){\circle{0.5}}

\put(20.5,6.5){\line(1,0){3}}

\put(21, 6){\circle*{0.5}}
\put(22, 6){\circle*{0.5}}
\put(23, 6){\circle*{0.5}}
\put(21, 5){\circle*{0.5}}
\put(22, 5){\circle*{0.5}}
\put(23, 5){\circle*{0.5}}
\put(21, 4){\circle*{0.5}}
\put(22, 4){\circle*{0.5}}
\put(23, 4){\circle*{0.5}}
\put(21, 3){\circle*{0.5}}
\put(22, 3){\circle*{0.5}}
\put(23, 3){\circle*{0.5}}

\put(21.2,2){$c_1$}
\put(22.2,2){$c_2$}

\put(19.5,4.2){$b_{-1}$}
\put(19.5,8.2){$b_{0}$}
\put(19.5,12.2){$b_{1}$}
\put(19.5,16.2){$b_{2}$}

\put(0,7){\begin{picture}(20,4)
\put(3,0){\ldots}
\put(3,1){\ldots}
\put(3,2){\ldots}
\put(3,3){\ldots}

\put(16,0){\ldots}
\put(16,1){\ldots}
\put(16,2){\ldots}
\put(16,3){\ldots}

\put(4.5,3){\circle*{0.5}}
\put(4.5,2){\circle*{0.5}}
\put(4.5,1){\circle*{0.5}}
\put(4.5,0){\circle*{0.5}}

\put(4.8, -1){$c_1$}

\put(5.5,3){\circle*{0.5}}
\put(5.5,2){\circle*{0.5}}
\put(5.5,1){\circle*{0.5}}
\put(5.5,0){\circle*{0.5}}

\put(5.8, -1){$c_2$}

\put(6.5,3){\circle*{0.5}}
\put(6.5,2){\circle*{0.5}}
\put(6.5,1){\circle*{0.5}}
\put(6.5,0){\circle*{0.5}}

\put(6.8, -1){$c_0$}

\put(7.5,3){\circle*{0.5}}
\put(7.5,2){\circle*{0.5}}
\put(7.5,1){\circle{0.5}}
\put(7.5,0){\circle*{0.5}}

\put(7.8, -1){$c_1$}

\put(8.5,3){\circle*{0.5}}
\put(8.5,2){\circle{0.5}}
\put(8.5,1){\circle*{0.5}}
\put(8.5,0){\circle*{0.5}}

\put(8.8, -1){$c_2$}

\put(9.5,3){\circle{0.5}}
\put(9.5,2){\circle*{0.5}}
\put(9.5,1){\circle*{0.5}}
\put(9.5,0){\circle{0.5}}

\put(9.8, -1){$c_0$}

\put(10,-0.5){\line(0,1){4}}
\put(7,-0.5){\line(0,1){4}}
\put(13,-0.5){\line(0,1){4}}
\put(5,4){$b_{-1}$}
\put(8,4){$b_{0}$}
\put(11,4){$b_{1}$}
\put(14,4){$b_{2}$}

\put(10.5,3){\circle*{0.5}}
\put(10.5,2){\circle*{0.5}}
\put(10.5,1){\circle{0.5}}
\put(10.5,0){\circle{0.5}}

\put(10.8, -1){$c_1$}

\put(11.5,3){\circle{0.5}}
\put(11.5,2){\circle{0.5}}
\put(11.5,1){\circle*{0.5}}
\put(11.5,0){\circle*{0.5}}

\put(11.8, -1){$c_2$}

\put(12.5,3){\circle{0.5}}
\put(12.5,2){\circle{0.5}}
\put(12.5,1){\circle{0.5}}
\put(12.5,0){\circle*{0.5}}

\put(12.8, -1){$c_0$}

\put(13.5,3){\circle{0.5}}
\put(13.5,2){\circle{0.5}}
\put(13.5,1){\circle{0.5}}
\put(13.5,0){\circle{0.5}}

\put(13.8, -1){$c_1$}

\put(14.5,3){\circle{0.5}}
\put(14.5,2){\circle{0.5}}
\put(14.5,1){\circle{0.5}}
\put(14.5,0){\circle{0.5}}

\put(14.8, -1){$c_2$}

\put(15.5,3){\circle{0.5}}
\put(15.5,2){\circle{0.5}}
\put(15.5,1){\circle{0.5}}
\put(15.5,0){\circle{0.5}}

\end{picture}}
\end{picture}}
\end{picture}

%% file: nottight.tex
\begin{picture}(20,4)

\put(4.5,3){\line(0,-1){2}}
\put(5.5,3){\line(0,-1){2}}
\put(6.5,3){\line(0,-1){2}}

\put(3.5,1){\line(1,-1){1}}
\put(4.5,1){\line(1,-1){1}}
\put(5.5,1){\line(1,-1){1}}
\put(6.5,1){\line(1,-1){1}}

\put(7.5,2){\line(1,-1){1}}
\put(8.5,3){\line(1,-1){1}}
\put(10.5,2){\line(1,-1){1}}
\put(11.5,1){\line(1,-1){1}}

\put(7.5,3){\line(0,-1){1}}
\put(8.5,1){\line(0,-1){1}}
\put(9.5,2){\line(0,-1){1}}
\put(10.5,3){\line(0,-1){1}}

\put(9.5,1){\line(2,-1){2}}

\put(9.33,0.4){$\downarrow$}
\put(9.33,1.4){$\downarrow$}
\put(8.33,2.4){$\downarrow$}
\put(11.33,-0.6){$\downarrow$}
\put(8.33,3.4){$\downarrow$}

\put(3,0){\ldots}
\put(3,1){\ldots}
\put(3,2){\ldots}
\put(3,3){\ldots}

\put(16,0){\ldots}
\put(16,1){\ldots}
\put(16,2){\ldots}
\put(16,3){\ldots}

\put(4.5,3){\circle*{0.5}}
\put(4.5,2){\circle*{0.5}}
\put(4.5,1){\circle*{0.5}}
\put(4.5,0){\circle*{0.5}}

\put(4.8, -1){$c_1$}

\put(5.5,3){\circle*{0.5}}
\put(5.5,2){\circle*{0.5}}
\put(5.5,1){\circle*{0.5}}
\put(5.5,0){\circle*{0.5}}

\put(5.8, -1){$c_2$}

\put(6.5,3){\circle*{0.5}}
\put(6.5,2){\circle*{0.5}}
\put(6.5,1){\circle*{0.5}}
\put(6.5,0){\circle*{0.5}}

\put(6.8, -1){$c_0$}

\put(7.5,3){\circle*{0.5}}
\put(7.5,2){\circle*{0.5}}
\put(7.5,1){\circle{0.5}}
\put(7.5,0){\circle*{0.5}}

\put(7.8, -1){$c_1$}

\put(8.5,3){\circle*{0.5}}
\put(8.5,2){\circle{0.5}}
\put(8.5,1){\circle*{0.5}}
\put(8.5,0){\circle*{0.5}}

\put(8.8, -1){$c_2$}

\put(9.5,3){\circle{0.5}}
\put(9.5,2){\circle*{0.5}}
\put(9.5,1){\circle*{0.5}}
\put(9.5,0){\circle{0.5}}

\put(9.8, -1){$c_0$}

\put(10,-0.5){\line(0,1){4}}

\put(10.5,3){\circle*{0.5}}
\put(10.5,2){\circle*{0.5}}
\put(10.5,1){\circle{0.5}}
\put(10.5,0){\circle{0.5}}

\put(10.8, -1){$c_1$}

\put(11.5,3){\circle{0.5}}
\put(11.5,2){\circle{0.5}}
\put(11.5,1){\circle*{0.5}}
\put(11.5,0){\circle*{0.5}}

\put(11.8, -1){$c_2$}

\put(12.5,3){\circle{0.5}}
\put(12.5,2){\circle{0.5}}
\put(12.5,1){\circle{0.5}}
\put(12.5,0){\circle*{0.5}}

\put(12.8, -1){$c_0$}

\put(13.5,3){\circle{0.5}}
\put(13.5,2){\circle{0.5}}
\put(13.5,1){\circle{0.5}}
\put(13.5,0){\circle{0.5}}

\put(13.8, -1){$c_1$}

\put(14.5,3){\circle{0.5}}
\put(14.5,2){\circle{0.5}}
\put(14.5,1){\circle{0.5}}
\put(14.5,0){\circle{0.5}}

\put(14.8, -1){$c_2$}

\put(15.5,3){\circle{0.5}}
\put(15.5,2){\circle{0.5}}
\put(15.5,1){\circle{0.5}}
\put(15.5,0){\circle{0.5}}

\end{picture}
\vspace{0.1in}

%% file: cyclic_move.tex
\put(-2,0){\begin{picture}(30,4)

\put(4,0){\ldots}
\put(4,1){\ldots}
\put(4,2){\ldots}
\put(4,3){\ldots}

\put(15,0){\ldots}
\put(15,1){\ldots}
\put(15,2){\ldots}
\put(15,3){\ldots}

\put(4.5,3){\line(1,-1){1}}
\put(5.5,2){\line(0,-1){1}}
\put(5.5,1){\line(1,-1){1}}

\put(5.5,3){\line(1,-1){1}}
\put(6.5,2){\line(0,-1){1}}
\put(6.5,1){\line(1,-1){1}}

\put(7.5,3){\line(0,-1){1}}
\put(7.5,2){\line(0,-1){1}}
\put(7.5,1){\line(1,-1){1}}

\put(8.5,3){\line(0,-1){1}}
\put(8.5,2){\line(1,-1){1}}
\put(9.5,1){\line(1,-1){1}}

\put(10.5,3){\line(0,-1){1}}
\put(10.5,2){\line(1,-1){1}}
\put(11.5,1){\line(0,-1){1}}

\put(11.5,3){\line(1,-1){1}}
\put(12.5,2){\line(0,-1){1}}
\put(12.5,1){\line(1,-1){1}}

\put(4.5,1){\line(1,-1){1}}

\put(5.5,3){\circle*{0.5}}
\put(5.5,2){\circle*{0.5}}
\put(5.5,1){\circle*{0.5}}
\put(5.5,0){\circle*{0.5}}

\put(5.8, -1){$c_2$}

\put(6.5,3){\circle{0.5}}
\put(6.5,2){\circle*{0.5}}
\put(6.5,1){\circle*{0.5}}
\put(6.5,0){\circle*{0.5}}

\put(6.8, -1){$c_0$}

\put(7.5,3){\circle*{0.5}}
\put(7.5,2){\circle*{0.5}}
\put(7.5,1){\circle*{0.5}}
\put(7.5,0){\circle*{0.5}}

\put(7.8, -1){$c_1$}

\put(8.5,3){\circle*{0.5}}
\put(8.5,2){\circle*{0.5}}
\put(8.5,1){\circle{0.5}}
\put(8.5,0){\circle*{0.5}}

\put(8.8, -1){$c_2$}

\put(9.5,3){\circle{0.5}}
\put(9.5,2){\circle{0.5}}
\put(9.5,1){\circle*{0.5}}
\put(9.5,0){\circle{0.5}}

\put(9.8, -1){$c_0$}

\put(10,-0.5){\line(0,1){4}}

\put(10.5,3){\circle*{0.5}}
\put(10.5,2){\circle*{0.5}}
\put(10.5,1){\circle{0.5}}
\put(10.5,0){\circle*{0.5}}

\put(10.8, -1){$c_1$}

\put(11.5,3){\circle*{0.5}}
\put(11.5,2){\circle{0.5}}
\put(11.5,1){\circle*{0.5}}
\put(11.5,0){\circle*{0.5}}

\put(11.8, -1){$c_2$}

\put(12.5,3){\circle{0.5}}
\put(12.5,2){\circle*{0.5}}
\put(12.5,1){\circle*{0.5}}
\put(12.5,0){\circle{0.5}}

\put(12.8, -1){$c_0$}

\put(13.5,3){\circle{0.5}}
\put(13.5,2){\circle{0.5}}
\put(13.5,1){\circle{0.5}}
\put(13.5,0){\circle*{0.5}}

\put(13.8, -1){$c_1$}

\put(14.5,3){\circle{0.5}}
\put(14.5,2){\circle{0.5}}
\put(14.5,1){\circle{0.5}}
\put(14.5,0){\circle{0.5}}

\end{picture}}

\put(14.5, 1.5){\vector(1,0){2}}

\put(13,0){\begin{picture}(30,4)

\put(4,0){\ldots}
\put(4,1){\ldots}
\put(4,2){\ldots}
\put(4,3){\ldots}

\put(15,0){\ldots}
\put(15,1){\ldots}
\put(15,2){\ldots}
\put(15,3){\ldots}

\put(4.5,2){\line(1,-1){1}}
\put(5.5,1){\line(0,-1){1}}

\put(5.5,2){\line(1,-1){1}}
\put(6.5,1){\line(0,-1){1}}

\put(7.5,2){\line(0,-1){1}}
\put(7.5,1){\line(0,-1){1}}

\put(8.5,2){\line(0,-1){1}}
\put(8.5,1){\line(1,-1){1}}

\put(10.5,2){\line(0,-1){1}}
\put(10.5,1){\line(1,-1){1}}

\put(11.5,2){\line(1,-1){1}}
\put(12.5,1){\line(0,-1){1}}

\put(4.5,3){\line(1,-1){1}}
\put(5.5,3){\line(2,-1){2}}
\put(7.5,3){\line(1,-1){1}}
\put(8.5,3){\line(2,-1){2}}
\put(10.5,3){\line(1,-1){1}}

\put(5.5,3){\circle*{0.5}}
\put(5.5,2){\circle*{0.5}}
\put(5.5,1){\circle*{0.5}}
\put(5.5,0){\circle*{0.5}}

\put(5.8, -1){$c_2$}

\put(6.5,3){\circle{0.5}}
\put(6.5,2){\circle{0.5}}
\put(6.5,1){\circle*{0.5}}
\put(6.5,0){\circle*{0.5}}

\put(6.8, -1){$c_0$}

\put(7.5,3){\circle*{0.5}}
\put(7.5,2){\circle*{0.5}}
\put(7.5,1){\circle*{0.5}}
\put(7.5,0){\circle*{0.5}}

\put(7.8, -1){$c_1$}

\put(8.5,3){\circle*{0.5}}
\put(8.5,2){\circle*{0.5}}
\put(8.5,1){\circle*{0.5}}
\put(8.5,0){\circle{0.5}}

\put(8.8, -1){$c_2$}

\put(9.5,3){\circle{0.5}}
\put(9.5,2){\circle{0.5}}
\put(9.5,1){\circle{0.5}}
\put(9.5,0){\circle*{0.5}}

\put(9.8, -1){$c_0$}

\put(10,-0.5){\line(0,1){4}}

\put(10.5,3){\circle*{0.5}}
\put(10.5,2){\circle*{0.5}}
\put(10.5,1){\circle*{0.5}}
\put(10.5,0){\circle{0.5}}

\put(10.8, -1){$c_1$}

\put(11.5,3){\circle{0.5}}
\put(11.5,2){\circle*{0.5}}
\put(11.5,1){\circle{0.5}}
\put(11.5,0){\circle*{0.5}}

\put(11.8, -1){$c_2$}

\put(12.5,3){\circle{0.5}}
\put(12.5,2){\circle{0.5}}
\put(12.5,1){\circle*{0.5}}
\put(12.5,0){\circle*{0.5}}

\put(12.8, -1){$c_0$}

\put(13.5,3){\circle{0.5}}
\put(13.5,2){\circle{0.5}}
\put(13.5,1){\circle{0.5}}
\put(13.5,0){\circle{0.5}}

\put(13.8, -1){$c_1$}

\put(14.5,3){\circle{0.5}}
\put(14.5,2){\circle{0.5}}
\put(14.5,1){\circle{0.5}}
\put(14.5,0){\circle{0.5}}

\end{picture}}

%% file: compact1.tex
\put(3,0){\ldots}
\put(3,1){\ldots}
\put(3,2){\ldots}
\put(3,3){\ldots}

\put(16,0){\ldots}
\put(16,1){\ldots}
\put(16,2){\ldots}
\put(16,3){\ldots}

\put(4.5,3){\line(1,-1){1}}
\put(5.5,2){\line(0,-1){1}}
\put(5.5,1){\line(1,-1){1}}

\put(5.5,3){\line(1,-1){1}}
\put(6.5,2){\line(0,-1){1}}
\put(6.5,1){\line(1,-1){1}}
\put(6.5,3){\line(1,-1){1}}
\put(7.5,2){\line(0,-1){1}}
\put(7.5,1){\line(1,-1){1}}
\put(7.5,3){\line(1,-1){1}}
\put(8.5,2){\line(0,-1){1}}
\put(8.5,1){\line(1,-1){1}}
\put(8.5,3){\line(1,-1){1}}
\put(9.5,2){\line(0,-1){1}}
\put(9.5,1){\line(1,-1){1}}
\put(9.5,3){\line(1,-1){1}}
\put(10.5,2){\line(0,-1){1}}
\put(10.5,1){\line(1,-1){1}}
\put(3.5,3){\line(1,-1){1}}
\put(4.5,2){\line(0,-1){1}}
\put(4.5,1){\line(1,-1){1}}
\put(3.5,1){\line(1,-1){1}}

\put(4.5,3){\circle*{0.5}}
\put(4.5,2){\circle*{0.5}}
\put(4.5,1){\circle*{0.5}}
\put(4.5,0){\circle*{0.5}}

\put(4.8, -1){$c_1$}

\put(5.5,3){\circle*{0.5}}
\put(5.5,2){\circle*{0.5}}
\put(5.5,1){\circle*{0.5}}
\put(5.5,0){\circle*{0.5}}

\put(5.8, -1){$c_2$}

\put(6.5,3){\circle*{0.5}}
\put(6.5,2){\circle*{0.5}}
\put(6.5,1){\circle*{0.5}}
\put(6.5,0){\circle*{0.5}}

\put(6.8, -1){$c_0$}

\put(7.5,3){\circle*{0.5}}
\put(7.5,2){\circle*{0.5}}
\put(7.5,1){\circle*{0.5}}
\put(7.5,0){\circle*{0.5}}

\put(7.8, -1){$c_1$}

\put(8.5,3){\circle*{0.5}}
\put(8.5,2){\circle*{0.5}}
\put(8.5,1){\circle*{0.5}}
\put(8.5,0){\circle*{0.5}}

\put(8.8, -1){$c_2$}

\put(9.5,3){\circle*{0.5}}
\put(9.5,2){\circle*{0.5}}
\put(9.5,1){\circle*{0.5}}
\put(9.5,0){\circle*{0.5}}

\put(9.8, -1){$c_0$}

\put(10,-0.5){\line(0,1){4}}

\put(10.5,3){\circle{0.5}}
\put(10.5,2){\circle*{0.5}}
\put(10.5,1){\circle*{0.5}}
\put(10.5,0){\circle*{0.5}}

\put(10.8, -1){$c_1$}

\put(11.5,3){\circle{0.5}}
\put(11.5,2){\circle{0.5}}
\put(11.5,1){\circle{0.5}}
\put(11.5,0){\circle*{0.5}}

\put(11.8, -1){$c_2$}

\put(12.5,3){\circle{0.5}}
\put(12.5,2){\circle{0.5}}
\put(12.5,1){\circle{0.5}}
\put(12.5,0){\circle{0.5}}

\put(12.8, -1){$c_0$}

\put(13.5,3){\circle{0.5}}
\put(13.5,2){\circle{0.5}}
\put(13.5,1){\circle{0.5}}
\put(13.5,0){\circle{0.5}}

\put(13.8, -1){$c_1$}

\put(14.5,3){\circle{0.5}}
\put(14.5,2){\circle{0.5}}
\put(14.5,1){\circle{0.5}}
\put(14.5,0){\circle{0.5}}

\put(14.8, -1){$c_2$}

\put(15.5,3){\circle{0.5}}
\put(15.5,2){\circle{0.5}}
\put(15.5,1){\circle{0.5}}
\put(15.5,0){\circle{0.5}}

%% file: good_abacus.tex
\put(3,0){\ldots}
\put(3,1){\ldots}
\put(3,2){\ldots}
\put(3,3){\ldots}

\put(16,0){\ldots}
\put(16,1){\ldots}
\put(16,2){\ldots}
\put(16,3){\ldots}

\put(3.5,3){\line(1,-1){1}}
\put(4.5,2){\line(0,-1){1}}
\put(4.5,1){\line(1,-1){1}}

\put(4.5,3){\line(1,-1){1}}
\put(5.5,2){\line(0,-1){1}}
\put(5.5,1){\line(1,-1){1}}

\put(5.5,3){\line(1,-1){1}}
\put(6.5,2){\line(0,-1){1}}
\put(6.5,1){\line(1,-1){1}}

\put(7.5,3){\line(0,-1){1}}
\put(7.5,2){\line(0,-1){1}}
\put(7.5,1){\line(1,-1){1}}

\put(8.5,3){\line(0,-1){1}}
\put(8.5,2){\line(1,-1){1}}
\put(9.5,1){\line(1,-1){1}}

\put(10.5,3){\line(0,-1){1}}
\put(10.5,2){\line(1,-1){1}}
\put(11.5,1){\line(0,-1){1}}

\put(11.5,3){\line(1,-1){1}}
\put(12.5,2){\line(0,-1){1}}
\put(12.5,1){\line(1,-1){1}}

\put(3.5,1){\line(1,-1){1}}

\put(4.5,3){\circle*{0.5}}
\put(4.5,2){\circle*{0.5}}
\put(4.5,1){\circle*{0.5}}
\put(4.5,0){\circle*{0.5}}

\put(4.8, -1){$c_1$}

\put(5.5,3){\circle*{0.5}}
\put(5.5,2){\circle*{0.5}}
\put(5.5,1){\circle*{0.5}}
\put(5.5,0){\circle*{0.5}}

\put(5.8, -1){$c_2$}

\put(6.5,3){\circle{0.5}}
\put(6.5,2){\circle*{0.5}}
\put(6.5,1){\circle*{0.5}}
\put(6.5,0){\circle*{0.5}}

\put(6.8, -1){$c_0$}

\put(7.5,3){\circle*{0.5}}
\put(7.5,2){\circle*{0.5}}
\put(7.5,1){\circle*{0.5}}
\put(7.5,0){\circle*{0.5}}

\put(7.8, -1){$c_1$}

\put(8.5,3){\circle*{0.5}}
\put(8.5,2){\circle*{0.5}}
\put(8.5,1){\circle{0.5}}
\put(8.5,0){\circle*{0.5}}

\put(8.8, -1){$c_2$}

\put(9.5,3){\circle{0.5}}
\put(9.5,2){\circle{0.5}}
\put(9.5,1){\circle*{0.5}}
\put(9.5,0){\circle{0.5}}

\put(9.8, -1){$c_0$}

\put(10,-0.5){\line(0,1){4}}

\put(10.5,3){\circle*{0.5}}
\put(10.5,2){\circle*{0.5}}
\put(10.5,1){\circle{0.5}}
\put(10.5,0){\circle*{0.5}}

\put(10.8, -1){$c_1$}

\put(11.5,3){\circle*{0.5}}
\put(11.5,2){\circle{0.5}}
\put(11.5,1){\circle*{0.5}}
\put(11.5,0){\circle*{0.5}}

\put(11.8, -1){$c_2$}

\put(12.5,3){\circle{0.5}}
\put(12.5,2){\circle*{0.5}}
\put(12.5,1){\circle*{0.5}}
\put(12.5,0){\circle{0.5}}

\put(12.8, -1){$c_0$}

\put(13.5,3){\circle{0.5}}
\put(13.5,2){\circle{0.5}}
\put(13.5,1){\circle{0.5}}
\put(13.5,0){\circle*{0.5}}

\put(13.8, -1){$c_1$}

\put(14.5,3){\circle{0.5}}
\put(14.5,2){\circle{0.5}}
\put(14.5,1){\circle{0.5}}
\put(14.5,0){\circle{0.5}}

\put(14.8, -1){$c_2$}

\put(15.5,3){\circle{0.5}}
\put(15.5,2){\circle{0.5}}
\put(15.5,1){\circle{0.5}}
\put(15.5,0){\circle{0.5}}

%% file: period_nl.tex
\setlength{\unitlength}{0.3cm}

\begin{center}
\begin{picture}(32,15)

\put(0,-17){ \begin{picture}(32,32)

\put(4,30){$\pi_5$}
\put(6.5,30){$\pi_4$}
\put(9,30){$\pi_3$}
\put(11.5,30){$\pi_2$}
\put(14,30){$\pi_1$}
\put(16.5,30){$\pi_0$}

\put(5,29){\vector(1,-2){0.8}}
\put(7.5,29){\vector(1,-2){0.8}}
\put(10,29){\vector(1,-2){0.8}}
\put(12.5,29){\vector(1,-2){0.8}}
\put(15,29){\vector(1,-2){0.8}}
\put(17.5,29){\vector(1,-2){0.8}}

\put(24,27){$c_2$}
\put(23,29){$c_1$}
\put(22,31){$c_0$}

\put(23.5,26.5){\vector(-2,-1){2}}
\put(22.5,28.5){\vector(-2,-1){2}}
\put(21.5,30.5){\vector(-2,-1){2}}

\put(6.8,24.9){\tiny{1}}
\put(8.8,25.9){\tiny{2}}
\put(11.8,24.9){\tiny{4}}
\put(13.8,25.9){\tiny{5}}
\put(15.8,26.9){\tiny{6}}
\put(18.8,25.9){\tiny{8}}

\put(20.5,25.4){\tiny{9}}
\put(17.5,26.4){\tiny{7}}
\put(10.5,25.4){\tiny{3}}

\linethickness{0.3mm}
\put(-1,32){\line(1,0){29}}
\put(6,32){\line(0,-1){15}}
\put(21,32){\line(0,-1){15}}

\put(6,24.1){\line(2,1){4}}
\put(6,24.05){\line(2,1){4}}
\put(10.05,26){\line(1,-2){1}}
\put(10.1,26){\line(1,-2){1}}
\put(11,24.1){\line(2,1){6}}
\put(11,24.05){\line(2,1){6}}
\put(17.05,27){\line(1,-2){1}}
\put(17.1,27){\line(1,-2){1}}
\put(18,25.1){\line(2,1){2}}
\put(18,25.05){\line(2,1){2}}
\put(20.05,26){\line(1,-2){1}}
\put(20.1,26){\line(1,-2){1}}

\linethickness{0.15mm}

\put(1.25,25.2){6}
\put(2.25,23.2){5}
\put(3.25,21.2){1}

\put(4.25,24.2){7}
\put(5.25,22.2){3}

\put(7.25,23.2){3}
\put(8.25,21.2){3}

\put(9.25,24.2){5}
\put(10.23,22.2){3}

\put(12.25,23.2){4}
\put(13.25,21.2){1}

\put(14.25,24.2){4}
\put(15.25,22.2){4}
\put(16.25,20.2){1}

\put(16.25,25.2){6}
\put(17.25,23.2){5}
\put(18.25,21.2){1}

\put(19.25,24.2){7}
\put(20.25,22.2){3}

\put(22.25,23.2){3}
\put(23.25,21.2){3}

\put(24.25,24.2){5}
\put(25.23,22.2){3}

\put(0,26){\line(1,-2){4}}
\put(2,27){\line(1,-2){4.5}}
\put(5,26){\line(1,-2){4}}
\put(8,25){\line(1,-2){3.5}}
\put(10,26){\line(1,-2){4}}
\put(13,25){\line(1,-2){3.5}}
\put(15,26){\line(1,-2){4}}
\put(17,27){\line(1,-2){4.5}}
\put(20,26){\line(1,-2){4}}
\put(23,25){\line(1,-2){3.5}}
\put(25,26){\line(1,-2){1}}

\put(2,27){\line(-2,-1){3}}
\put(5,26){\line(-2,-1){6}}
\put(10,26){\line(-2,-1){11}}
\put(17,27){\line(-2,-1){18}}
\put(20,26){\line(-2,-1){16}}
\put(25,26){\line(-2,-1){16}}
\put(26,24){\line(-2,-1){12}}
\put(26,21.5){\line(-2,-1){7}}
\put(26,19){\line(-2,-1){2}}

 \end{picture}
}
\end{picture}

\end{center}

\setlength{\unitlength}{0.5cm}

%% file: cyclic.tex
\put(13, 31){$y$}
\put(76,31){$x$}
\put(45,0){\vector(1,1){30}}

\put(45,0){\vector(-1,1){30}}

\put(48,3){\circle*{1}}
\put(51,6){\circle*{1}}
\put(54,9){\circle*{1}}
\put(57,12){\circle*{1}}
\put(60,15){\circle*{1}}
\put(63,18){\circle*{1}}
\put(66,21){\circle*{1}}
\put(69,24){\circle*{1}}

\put(48.5,1.7){$\psi_4 = \psi_{0, (3)}$}
\put(50,4){$\psi_3$}
\put(53,7){$\psi_2$}
\put(54.5,8){$\psi_1$}
\put(57.5,11.5){$\psi_0$}

\put(57.3,12.3){\line(-1,1){3}}
\put(54.3,15.3){\line(-1,-1){3}}
\put(51.3,12.3){\line(-1,1){6}}
\put(45.3,18.3){\line(-1,-1){3}}
\put(42.3,15.3){\vector(-1,1){12}}

\put(54,9){\line(-1,1){6}}
\put(48,15){\line(-1,-1){3}}
\put(45,12){\line(-1,1){3}}
\put(42,15){\line(-1,-1){3}}
\put(39,12){\vector(-1,1){12}}

\put(53.7,8.7){\line(-1,1){3}}
\put(50.7,11.7){\line(-1,-1){3}}
\put(47.7,8.7){\line(-1,1){3}}
\put(44.7,11.7){\line(-1,-1){3}}
\put(41.7,8.7){\vector(-1,1){15}}

\put(50.4,5.4){\line(-1,1){5.7}}
\put(44.7,11.1){\line(-1,-1){3}}
\put(41.7,8.1){\line(-1,1){6}}
\put(35.7,14.1){\line(-1,-1){2.7}}
\put(33,11.4){\vector(-1,1){9}}

\put(48.3,3.3){\line(-1,1){3}}
\put(45.3,6.3){\line(-1,-1){3}}
\put(42.3,3.3){\line(-1,1){6}}
\put(36.3,9.3){\line(-1,-1){3}}
\put(33.3,6.3){\vector(-1,1){12}}

%% file: ccp2.tex
\begin{picture}(40, 16)

\put(13, 10){\line(0,-1){10}}
\put(27, 10){\line(0,-1){10}}

\put(17,15.5){$\pi_0$}
\put(18,15){\vector(1,-1){2}}

\put(15,14.5){$\pi_1$}
\put(16,14){\vector(1,-1){2}}

\put(13,13.5){$\pi_2$}
\put(14,13){\vector(1,-1){2}}

\put(11,12.5){$\pi_3$}
\put(12,12){\vector(1,-1){2}}

\put(9,11.5){$\pi_4$}
\put(10,11){\vector(1,-1){2}}

\put(7,10.5){$\pi_5$}
\put(8,10){\vector(1,-1){2}}

\put(24,14.5){$c_{0}$}
\put(26,13.17){$c_{1}$}
\put(28,11.83){$c_{2}$}
\put(30,10.5){$c_{3}$}
\put(32,9.17){$c_{4}$}

\put(23.5,14){\vector(-2,-1){2}}
\put(25.5,12.67){\vector(-2,-1){2}}
\put(27.5,11.33){\vector(-2,-1){2}}
\put(29.5,10){\vector(-2,-1){2}}
\put(31.5,8.67){\vector(-2,-1){2}}

\put(0,-8){\begin{picture}(40,20)
\put(20,20){}
\put(22,18.67){}
\put(24,17.33){$3$}
\put(26,16){$1$}
\put(28,14.67){$0$}
\put(30,13.33){$0$}
\put(32,12){$0$}

\put(18,19){}
\put(20,17.67){$3$}
\put(22,16.33){$2$}
\put(24,15){$0$}
\put(26,13.67){$0$}
\put(28,12.33){$0$}

\put(16,18){}
\put(18,16.67){$2$}
\put(20,15.33){$1$}
\put(22,14){$0$}
\put(24,12.67){$0$}
\put(26,11.33){$0$}

\put(14,17){$4$}
\put(16,15.67){$2$}
\put(18,14.33){$0$}
\put(20,13){$0$}
\put(22,11.67){$0$}

\put(10,17.33){$3$}
\put(12,16){$1$}
\put(14,14.67){$0$}
\put(16,13.33){$0$}
\put(18,12){$0$}
\put(20,10.67){$0$}

\put(8,16.33){$2$}
\put(10,15){$0$}
\put(12,13.67){$0$}
\put(14,12.33){$0$}
\put(16,11){$0$}

\end{picture} }

\end{picture}

%% file: perfect_crystal.tex
\begin{picture}(30,18)

\put(17.5, 15.5){\vector(-4,-3){2.5}}
\put(13.5, 11.5){\vector(-4,-3){2.5}}
\put(17.5, 7.5){\vector(-4,-3){2.5}}
\put(16,15){1}
\put(12,11){1}
\put(16,7){1}
\put(14.5, 11.5){\vector(4,-3){2.5}}
\put(10.5, 7.5){\vector(4,-3){2.5}}
\put(14.5, 3.5){\vector(4,-3){2.5}}
\put(16,10.7){2}
\put(12,6.7){2}
\put(16,2.7){2}
\put(18,9.5){\vector(0,1){6}}
\put(18,1.5){\vector(0,1){6}}
\put(14,5.5){\vector(0,1){6}}
\put(18.3,12.5){0}
\put(18.3,4.5){0}
\put(14.3,8.5){0}

\put(17,16){\begin{picture}(2,1)
\put(0,0){\line(1,0){2}}
\put(0,1){\line(1,0){2}}
\put(0,0){\line(0,1){1}}
\put(1,0){\line(0,1){1}}
\put(2,0){\line(0,1){1}} 
\put(0.1, 0.32){\small{0.5}}
\put(1.1, 0.32){\small{0.5}}
\end{picture}}

\put(13,12){\begin{picture}(2,1)
\put(0,0){\line(1,0){2}}
\put(0,1){\line(1,0){2}}
\put(0,0){\line(0,1){1}}
\put(1,0){\line(0,1){1}}
\put(2,0){\line(0,1){1}} 
\put(0.1, 0.32){\small{0.5}}
\put(1.1, 0.32){\small{1.5}}
\end{picture}}

\put(9,8){\begin{picture}(2,1)
\put(0,0){\line(1,0){2}}
\put(0,1){\line(1,0){2}}
\put(0,0){\line(0,1){1}}
\put(1,0){\line(0,1){1}}
\put(2,0){\line(0,1){1}} 
\put(0.1, 0.32){\small{1.5}}
\put(1.1, 0.32){\small{1.5}}
\end{picture}}

\put(17,8){\begin{picture}(2,1)
\put(0,0){\line(1,0){2}}
\put(0,1){\line(1,0){2}}
\put(0,0){\line(0,1){1}}
\put(1,0){\line(0,1){1}}
\put(2,0){\line(0,1){1}} 
\put(0.1, 0.32){\small{0.5}}
\put(1.1, 0.32){\small{2.5}}
\end{picture}}

\put(13,4){\begin{picture}(2,1)
\put(0,0){\line(1,0){2}}
\put(0,1){\line(1,0){2}}
\put(0,0){\line(0,1){1}}
\put(1,0){\line(0,1){1}}
\put(2,0){\line(0,1){1}} 
\put(0.1, 0.32){\small{1.5}}
\put(1.1, 0.32){\small{2.5}}
\end{picture}}

\put(17,0){\begin{picture}(2,1)
\put(0,0){\line(1,0){2}}
\put(0,1){\line(1,0){2}}
\put(0,0){\line(0,1){1}}
\put(1,0){\line(0,1){1}}
\put(2,0){\line(0,1){1}} 
\put(0.1, 0.32){\small{2.5}}
\put(1.1, 0.32){\small{2.5}}
\end{picture}}

\end{picture}